\documentclass[a4paper,10pt,reqno]{amsart}
\usepackage{amsmath, amssymb, amsthm, amsfonts, amsgen, stmaryrd}
\usepackage{indentfirst}
\usepackage{graphics}
\usepackage{psfig}

 \numberwithin{equation}{section}

\makeatletter
\def\rightharpoonupfill@{\arrowfill@\relbar\relbar\rightharpoonup}
\newcommand{\xrightharpoonup}[2][]{\ext@arrow
0359\rightharpoonupfill@{#1}{#2}} \makeatother

\newcommand{\ds}{\displaystyle}

\newcommand{\medi}{- \hskip-0.9em \int}
\newcommand{\med}{- \hskip-1.0em \int}

\newcommand{\C}{{\mathcal{C}}}

\newcommand{\Q}{{\mathcal{Q}}}
\newcommand{\Nb}{{\mathbb{N}}}

\newcommand{\Rb}{{\mathbb{R}}}
\newcommand{\Zb}{{\mathbb{Z}}}

\let\e=\varepsilon
\let\a=\alpha
\let\d=\delta

\let\g=\gamma

\let\O=\Omega
\let\o=\omega
\let\G=\Gamma

\newtheorem{thm}{Theorem}[section]
\newtheorem{defi}[thm]{Definition}
\newtheorem{rmk}[thm]{Remark}
\newtheorem{lemma}[thm]{Lemma}
\newtheorem{proposition}[thm]{Proposition}

\newcommand{\loc}{{\rm loc}}
\newcommand{\dx}{\,dx}
\newcommand{\dxa}{\dx_{\alpha}}
\newcommand{\dy}{\,dy}

\newcommand{\ie}{; {\it i.e., }}
\newcommand{\HH}{{\mathcal H}}

\def\rr{\Bbb R}

\begin{document}
\title[The Neumann sieve problem and dimensional reduction]
{The Neumann sieve problem and dimensional reduction: a multiscale
approach}
\author[N. Ansini, J.-F. Babadjian \& C. I. Zeppieri]
{Nadia Ansini, Jean-Fran\c cois Babadjian \& Caterina Ida Zeppieri}
\address[Nadia Ansini\footnote{Present Address: Section de Math\'ematiques,
EPFL, 1015 Lausanne, Switzerland. E-mail:
nadia.ansini@epfl.ch}]{Dipartimento di Matematica `G. Castelnuovo'
Universit\`a di Roma `La Sapienza', Piazzale Aldo Moro, 2, 00185
Rome, Italy} \email{ansini@mat.uniroma1.it}

\address[Jean-Fran\c cois Babadjian]{SISSA, Via Beirut 2-4, 34014
Trieste, Italy} \email{babadjia@sissa.it}

\address [Caterina Ida Zeppieri]{Dipartimento di Matematica `G.
Castelnuovo' Universit\`a di Roma `La Sapienza', Piazzale Aldo Moro,
2, 00185 Rome, Italy} \email{zeppieri@mat.uniroma1.it}

\maketitle

\begin{abstract} We perform a multiscale
analysis for the elastic energy of a $n$-dimensional bilayer thin
film of thickness $2\d$ whose layers are connected through an
$\e$-periodically distributed contact zone. Describing the contact
zone as a union of $(n-1)$-dimensional balls of radius $r\ll \e$
(the holes of the sieve) and assuming that $\d \ll \e$,  we show
that the asymptotic memory of the sieve (as $\e \to 0$) is
witnessed by the presence of an extra interfacial energy term.
Moreover we find three different limit behaviors (or regimes)
depending on the mutual vanishing rate of $\d$ and $r$. We also
give an explicit nonlinear capacitary-type formula for the
interfacial energy density in each regime.
\end{abstract}
\begin{center}
\begin{minipage}{13.5cm}
\vspace{0.7cm}

\small{ \noindent \textsc{Keywords}: $\Gamma$-convergence,
dimension reduction, sieve problem, nonlinear capacity.

\vspace{6pt} \noindent {\it 2000 Mathematics Subject
Classification:} 49J45, 74K35, 74K15, 74B20, 74G65, 35B40. }
\end{minipage}
\end{center}

\vspace{0.5cm}


\section{Introduction}
\noindent For an ever increasing variety of applications, an
interesting problem to be explored is to model the debonding of a
thin film from a substrate.

If we consider a stretched film bonded to an infinite rigid
substrate, the elastic energy of this film scales as its thickness.
If the film debonds from the substrate, on one hand its elastic
energy tends to zero, while on the other hand this creates a new
surface and then an interfacial energy independent of the thickness.

In \cite{BhFoFra} Bhattacharya, Fonseca and Francfort examine,
among other, the asymptotic behavior of a bilayer thin film
allowing for the possibility of a debonding at the interface, but
penalizing it postulating an interfacial energy which scales as
the overall thickness of the film to some exponent. Thus the
energy they consider consists of the elastic energy of the two
layers and the interfacial energy with penalized debonding.

In this paper we deal with thin films connected by a hyperplane
(sieve plane) through a periodically distributed contact zone.
Thus we see the debonding as the effect of the {\it weak
interaction} of the two thin films through this contact zone and
we recover the interfacial energy term by a limit procedure.

Since we are mainly interested in describing the interaction
phenomenon due to the presence of the sieve, we make a
simplification choosing two thin films having the same elastic
properties (for a generalization to the case of two different
materials interacting, we refer the reader to \cite{Ans}).

\bigskip Consider a nonlinear elastic
$n$-dimensional bilayer thin film of thickness $2\d$ with layers
connected through $(n-1)$-dimensional balls $B^{n-1}_{r}(x_i^\e)$
centered in $x_i^\e:=i\e$, $i\in \Zb^{n-1}$ and with radius $r>0$.
Thus the investigated elastic body occupies the reference
configuration parametrized as
$$
\O^{\d}_{\e,r}:=\o^{+\d} \cup \o^{-\d} \cup \big( \o_{\e, r}
\times \{0\} \big)
$$
where $\o$ is a bounded open subset of $\Rb^{n-1}$, $\o^{+\d}:=\o
\times(0,\d)$, $\o^{-\d}:=\o \times (-\d,0)$ and $\o_{\e, r}:=
\bigcup_{i \in\Zb^{n-1}}B^{n-1}_{r}(x_i^\e) \cap \o$ (see Figure
1).

\begin{figure}[th]\label{fig1}
\centerline{\psfig{file=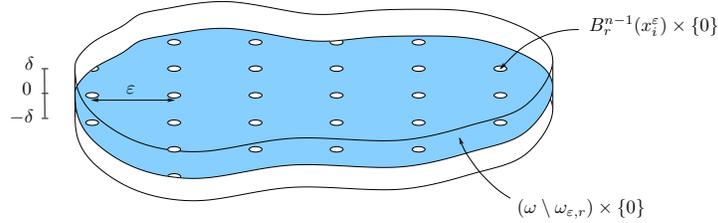,width=4in}} \vspace*{8pt}
\caption{The domain $\Omega^{\delta}_{\varepsilon, r}$.}
\end{figure}

In the nonlinear membrane theory setting the (scaled) elastic
energy associated to the material modelled by $\O^{\d}_{\e,r}$ is
given by
\begin{equation}\label{energia}
\frac{1}{\d} \int_{\O^{\d}_{\e,r}} W(D u )\, dx,
\end{equation}
where $u:\O^{\d}_{\e,r} \to \Rb^m$ is the deformation field and
$W$ is the stored energy density.

The $\G$-convergence approach has been used successfully in recent
years to rigorously obtain limit models for various dimensional
reductional problems (see for example \cite{BhJ, BraFo, BFF, LDR,
Shu}).

In this paper we study the multiscale asymptotic behavior of
(\ref{energia}) via $\Gamma$-convergence, as $\e$, $\d$ and $r$ tend
to zero, assuming that $\d= \d(\e)$, $r =r(\e,\d)$ and where $W:
\Rb^{m\times n}\to [0,+\infty)$ is a Borel function satisfying a
growth condition of order $p$, with $1<p<n-1$. The case $p=n-1$
requires a further appropriate analysis and it cannot be easily
derived from $p<n-1$ by slight changes. Unfortunately, three
dimensional linearized elasticity falls into this framework.

Since the sieve $(\o\setminus\o_{\e,r})\times \{0\}$ is not a part
of the domain $\O^{\d}_{\e,r}$, for any fixed $\e, \d, r>0$ we
have no information on the admissible deformation across part of
the mid-section $\o\times\{0\}$. This possible lack of regularity
might produce, in the limit, the above mentioned debonding and
correspondingly an interfacial energy depending on the jump of the
limit deformation. Moreover, we expect that this interfacial
energy will depend on the scaling of the radius of the connecting
zones with respect to the period of their distribution and the
thickness of the thin film.

\bigskip

The cases $\d=1$ and $\d=\e$ have been studied by Ansini
\cite{Ans} who proved that, to recover a non trivial limit
model\ie to obtain a limit model remembering the presence of the
sieve, the meaningful radius (or critical size) of the contact
zones must be of order $\e^{(n-1)/(n-p)}$ and $\e^{n/(n-p)}$,
respectively. In fact a different choice should lead in the limit
to two decoupled problems (if $r$ tends to zero faster than the
critical size) or to the same result that is obtained without the
presence of connecting zones in the mid-section (if $r$ tends to
zero more slowly than the critical size).

The proofs of the $\Gamma$-convergence results in \cite{Ans} (see
Theorems 3.2 and 8.2 therein) are based on a technical lemma
(\cite{Ans}, Lemma 3.4) that allows to modify a sequence of
deformations $u_\e$ with equi-bounded energy, on a suitable
$n$-dimensional spherical annuli surrounding the balls
$B^{n-1}_{r}(x_i^\e)$ without essentially changing their energies,
and to study the behavior of the energies along the new modified
sequence. Both in the case $\d=1$ and $\d=\e$ the $\Gamma$-limits
consist of three terms. The first two terms represent the
contribution of the new sequence far from the balls
$B^{n-1}_{r}(x_i^\e)$; more precisely, they are the
$\Gamma$-limits of two problems defined separately on the upper
and lower part (with respect to the `sieve plane') of the
considered domain. The third term describes the contribution near
the balls $B^{n-1}_{r}(x_i^\e)$ through a  nonlinear
capacitary-type formula that is the same for both $\d=1$ and
$\d=\e$. The equality of the two formulas is due to the fact that
the radii of the annuli suitably chosen to separate the two
contributions are less than $c\,\e$, with $c$ an arbitrary small
positive constant. In fact as a consequence, all constructions can
be performed in the interior of the domain, and the same procedure
yielding the nonlinear capacitary-type formula, applies for $\d=1$
and for $\d=\e$ as well. The cases $\e \sim \d$ and $\e \ll \d$
can be treated in the same way.

This approach follows the method introduced by Ansini-Braides in
\cite{AnsB,erratum} where the asymptotic behavior of periodically
perforated nonlinear domains has been studied; in particular,
Lemma 3.4 in \cite{Ans} is a suitable variant, for the sieve
problem, of Lemma 3.1 in \cite{AnsB}.

For other problems related to this subject, we refer the reader to
Attouch-Damlamian-Murat-Picard \cite{ADMP}, \cite{Mur},
\cite{Pic}, Attouch-Picard \cite{AP}, Conca \cite{Con1, Con2,
Con3}, Del Vecchio \cite{Del} and Sanchez-Palencia \cite{S.P.,
S.P.1, S.P.2}, among others.

\bigskip
In this paper we focus our attention on the case $\d=\d(\e)\ll\e$.
As in \cite{Ans}, we expect the existence of a critical radius
$r=r(\e,\d) \ll \e$ for which the limit model is nontrivial but now
we expect also to find different limit regimes depending on the
mutual vanishing rate of $r$ and $\d$. Moreover Lemma 3.4 in
\cite{Ans} cannot be directly applied to our setting since the
spherical annuli surrounding the connecting zones
$B^{n-1}_{r}(x_i^\e)$ as above, are well contained in a strip of
thickness $c\,\e$ but not in $\O^\d_{\e,r}$ ($\d\ll\e$). However, we
are able to modify Lemma 3.4 in \cite{Ans} by considering, instead
of spherical annuli, suitable cylindrical annuli of thickness of
order $\d$ (see Lemma \ref{important} and Lemma
\ref{important+equiint}).

As a consequence, also in this case the asymptotic analysis of
(\ref{energia}) as $\e$, $\d$ and $r$ tend to zero can be carried
on studying separately the energy contributions far from and close
to $B^{n-1}_{r}(x_i^\e)$; we get three terms in the limit. The
first two terms still describe the contribution `far' from the
connecting zones\ie they are the $\Gamma$-limits of the two
dimensional reduction problems defined by
$$
\frac{1}{\d} \int_{\o^{+\d}} W(D u )\, dx\,,\qquad \frac{1}{\d}
\int_{\o^{-\d}} W(D u )\, dx\,;
$$
while the third term, arising in the limit from the energy
contribution close to the connecting zones, represents the
asymptotic memory of the sieve: it is the above mentioned
interfacial energy.

\bigskip

The main results of this paper are stated in Theorem \ref{ABZ} and
Theorem \ref{ABZ2}. In Theorem \ref{ABZ} we prove a
$\Gamma$-convergence result for the sequence of functionals
(\ref{energia}) while in Theorem \ref{ABZ2} we give an explicit
characterization of the interfacial energy term occurring in the
$\Gamma$-limit. More precisely, for every sequence $(\e_j)$
converging to zero, we set $\d_j:=\d(\e_j)$, $r_j:= r(\e_j,\d_j)$,
$\O_j:=\O^{\d_j}_{\e_j,r_j}$ and
$$
{\mathcal F}_j(u):=\left\{
\begin{array}{ll}
\ds \frac{1}{\d_j} \int_{\O_j} W(D
u)\, dx  & \text{if }u \in W^{1,p}(\O_j;\Rb^m)\\
&\\
\ds +\infty & \text{otherwise\,.}
\end{array}
\right.
$$
Up to subsequence we can define
$$
\ell:= \lim_{j \to +\infty}\frac{r_j}{\d_j} \quad \text{and} \quad
g(F):=\lim_{j\to +\infty} r_j^p\, \Q_n W (r_j^{-1} F).
$$
where $\Q_n W$ is the $n$-quasiconvexification of $W$.

If $\ell \in (0,+\infty]$ and
$$
0< R^{(\ell)}:=\lim_{j \to
+\infty}\frac{r_j^{n-1-p}}{\e_j^{n-1}}<+\infty,
$$
then $({\mathcal F}_j)$ $\Gamma$-converges to
$$
{\mathcal F}^{(\ell)}(u^+,u^-)=\int_\o \Q_{n-1}\overline W(D_\a
u^+)\, dx_\a + \int_\o \Q_{n-1}\overline W(D_\a u^-)\, dx_\a +
R^{(\ell)} \int_\o \varphi^{(\ell)} (u^+ - u^-)\, dx_\a
$$
on $W^{1,p}(\o;\Rb^m) \times W^{1,p}(\o;\Rb^m)$ with respect to
the convergence introduced in Definition {\rm\ref{conv}}, where
$\overline W(\overline F):=\inf \{ W(\overline F|z) : z \in
\Rb^m\}$, $\Q_{n-1}\overline W$ is the
$(n-1)$-quasiconvexification of $\overline W$ and
$\varphi^{(\ell)}:\Rb^m \to [0,+\infty)$ is a locally Lipschitz
continuous function for any $\ell\in [0,+\infty]$. Similarly, if
$\ell =0$ and
$$
0< R^{(0)}:=\lim_{j \to +\infty}\frac{r_j^{n-p}}{\d_j\,
\e_j^{n-1}}<+\infty,
$$
then we still have $\Gamma$-convergence, as above, to
$$
{\mathcal F}^{(0)}(u^+,u^-)=\int_\o \Q_{n-1}\overline W(D_\a
u^+)\, dx_\a + \int_\o \Q_{n-1}\overline W(D_\a u^-)\, dx_\a +
R^{(0)} \int_\o \varphi^{(0)} (u^+ - u^-)\, dx_\a
$$
on $W^{1,p}(\o;\Rb^m) \times W^{1,p}(\o;\Rb^m)$.
\smallskip

For any $\ell\in [0,+\infty]$, $\varphi^{(\ell)}$ is described by
the following nonlinear capacitary-type formulas:

(1) if $\ell=+\infty$, then
\begin{eqnarray*}
\varphi^{(\infty)}(z) & = & \inf \Bigg\{\int_{\Rb^{n-1}}\Bigl(
\Q_{n-1} \,\overline g (D_\a \zeta^+ ) + \Q_{n-1} \, \overline g
(D_\a \zeta^- ) \Bigr)\, dx_\a : \; \zeta^\pm \in
W^{1,p}_{\rm loc}(\Rb^{n-1};\Rb^m),\\
&&\hspace{3.5 cm}\zeta^+=\zeta^- \text{ in }B_1^{n-1}(0), \quad
D_\a \zeta^\pm
\in L^p(\Rb^{n-1} ;\Rb^{m \times (n-1)}),\\
&&\hspace{7.5cm} (\zeta^+ - z) \,, \zeta^- \in
L^{p^*}(\Rb^{n-1};\Rb^m) \Bigg\},
\end{eqnarray*}
where again, $\bar g (\overline F):=\inf \{ g(\overline F|z) : z
\in \Rb^m\}$ and $\Q_{n-1}\bar g$ is the
$(n-1)$-quasiconvexification of $\bar g$,

\bigskip

(2) if $\ell=0$, then
\begin{eqnarray*}
\varphi^{(0)}(z) & = & \inf \Bigg\{\int_{\Rb^n\setminus
C_{1,\infty} } g (D \zeta) \, dx :\;\zeta \in W^{1,p}_{\rm
loc}(\Rb^n \setminus C_{1,\infty};\Rb^m), \, D \zeta \in L^p(\Rb^n
\setminus C_{1,\infty};\Rb^{m\times n}),\\&& \hspace{1cm}
\,\zeta-z \in L^{p}(0,+\infty;L^{p^*}(\Rb^{n-1};\Rb^m)) \text{ and
} \zeta \in L^{p}(-\infty,0;L^{p^*}(\Rb^{n-1};\Rb^m))\Bigg\}\,,
\end{eqnarray*}

\bigskip

(3) if $\ell \in (0,+\infty)$, then
\begin{eqnarray*}
\varphi^{(\ell)}(z) & = & \inf \Bigg\{\int_{\Rb^{n-1} \times
(-1,1)}g \big(D_\a \zeta | \ell D_n \zeta \big)\, dx : \zeta \in
W^{1,p}_{\rm loc}((\Rb^{n-1} \times (-1,1)) \setminus
C_{1,\infty};\Rb^m),\\
&&\hspace{1cm} D \zeta \in L^p(\Rb^{n-1} \times (-1,1);\Rb^m),
\;\qquad \zeta - z \in L^p((0,1);L^{p^*}(\Rb^{n-1};\Rb^m))\\
&&\hspace{7.5 cm} \zeta \in L^p((-1,0);
L^{p^*}(\Rb^{n-1};\Rb^m))\Bigg\}\,,
\end{eqnarray*}
where $C_{1,\infty}:=\{(x_\a,0)\in\Rb^n \; :\, 1\leq|x_\a| \}$.

\bigskip
Before giving a brief heuristic description of each regime, we
want to point out that whatever the value of $\ell$ is, the
interfacial energy density $\varphi^{(\ell)}$ corresponds to a
cohesive interface where the surface energy increases continuously
from zero with the jump in the deformation across the interface.

\bigskip
(1) The case $\ell=+\infty$ corresponds to $\d_j \ll r_j \ll
\e_j$, thus we expect $r_j$ to depend only on $\e_j$. In this case
we have a separation of scales effect. We first consider $r_j$ and
$\e_j$ as `fixed' and let $\d_j$ tend to zero as if we were
dealing with two pure dimensional reduction problems stated
separately on the upper and lower part (with respect to the sieve
plane) of $\O_j$. Then this first limit procedure yields two
functionals being both a copy of the functional in \cite{LDR}.
Since the two corresponding limit deformations $u^+$ and $u^-$
must match inside each connecting zone, the above two terms are
not completely decoupled. We are then in a situation quite similar
to that of \cite{AnsB,erratum}, except that here both periodically
`perforated' $(n-1)$-dimensional bodies are linked each other
through the `perforations'\ie through the holes of the sieve and
not through the sieve itself. Thus it is coherent to find a
critical size of order $\e^{(n-1)/(n-1-p)}$. Moreover this strong
separation between the phenomena of dimension reduction and
`perforation' leads to anisotropy as it can be seen, for instance,
also by an inspection of the proof of Lemma \ref{gliminfclose}
which shows that the extra interfacial energy term appears thanks
to suitable dilatations having a different scaling in the in-plane
and transverse variables. Finally we note that the formula for
$\varphi^{(\infty)}$ is given in terms of a `Le Dret-Raoult type'
functional involving the limit of the right capacitary scaling
(that is, involving the function $g$).

\smallskip

(2) The case $\ell=0$ corresponds to $r_j \ll \d_j \ll \e_j$. In
this case we expect that the critical size $r_j$ depends on both
$\d_j$ and $\e_j$. Indeed, as already pointed out, $r_j$ is of
order $\d_j^{1/(n-p)} \e_j^{(n-1)/(n-p)}$. Note that for
$\d_j=\e_j$ we recover $\e^{n/(n-p)}$ that is the critical size
obtained in \cite{Ans}; moreover $\varphi^{(0)}$ turns out to
coincide with the function $\varphi$ in \cite{Ans} (see Remark
\ref{Nadia}). Contrary to the previous case, now the isotropy is
preserved in fact here the dimensional reduction and `perforation'
processes are not completely decoupled: the reduction parameter
$\d_j$ is forced between both parameters $r_j$ and $\e_j$. This
can be seen also by noticing that now the scaling leading to the
interfacial energy is the same in every direction (see for
instance the proof of the $\Gamma$-limsup inequality). Moreover
now in $\varphi^{(0)}$ the reduction procedure is not explicit but
only witnessed by the boundary conditions expressed only on the
lateral part of the boundary of the considered domain.

\smallskip

(3) The case $\ell \in (0,+\infty)$ corresponds to $r_j \sim \d_j
\ll \e_j$. In this case the separation of scales effect does not
take place and the two previous scalings turn out to be equivalent
($R^{(0)}=\ell \, R^{(\infty)}$). Moreover we find that the
interfacial energy is continuous with respect to $\ell$ in the
extreme regimes\ie $R^{(\ell)} \varphi^{(\ell)}(z) \to
R^{(\infty)}\varphi^{(\infty)}(z)$ as $\ell \to +\infty$ and
$R^{(\ell)} \varphi^{(\ell)}(z) \to R^{(0)}\varphi^{(0)}(z)$ as
$\ell \to 0$. Finally, as in the previous case, the lateral
boundary conditions are the only mean describing the dimensional
reduction phenomenon in the procedure leading to
$\varphi^{(\ell)}$.

\bigskip

This paper is organized as follows: after recalling some useful
notation in Section \ref{notation}, we state the main results,
Theorem \ref{ABZ} and Theorem \ref{ABZ2}, in Section
\ref{mainres}. Then, in Section \ref{prel} we list some auxiliary
results as rescaled Poincar\'e type inequalities and joining
lemmas. Section \ref{close} is devoted to give a preliminary
definition of the interfacial energy density as limit of minimum
problems. In Section \ref{proof} we prove the $\Gamma$-convergence
result (Theorem \ref{ABZ}). It is only in Section
\ref{nonabstract} that we compute the explicit expression of the
interfacial energy density for each regime (Theorem \ref{ABZ2}).


\section{Notation}\label{notation}

\noindent Given $x \in \Rb^n$, we set $x=(x_\a,x_n)$ where
$x_\a:=(x_1,\ldots,x_{n-1})$ is the in-plane variable and
$D_\a:=\left(\frac{\partial}{\partial
x_1},\ldots,\frac{\partial}{\partial x_{n-1}}\right)$ (resp.
$D_n$) the derivative with respect to $x_\a$ (resp. $x_n$).

The notation $\Rb^{m \times n}$ stands for the set of $m \times n$
real matrices. Given a matrix $F \in \Rb^{m \times n}$, we write
$F=(\overline F|F_n)$ where $\overline F=(F_1, \ldots , F_{n-1})$
and $F_i$ denotes the $i$-th column of $F$, $1 \leq i \leq n$ and
$\overline F\in \Rb^{m \times (n-1)}$.

The Lebesgue measure in $\Rb^n$ will be denoted by $\mathcal L^n$
and the Hausdorff $(n-1)$-dimensional measure by $\HH^{n-1}$. Let
$A$ be an open subset of $\Rb^d$ ($d=n-1$, $d=n$). If $s\in
[1,+\infty]$, we use standard notation for Lebesgue and Sobolev
spaces $L^s(A;\rr^m)$ and $W^{1,s}(A;\rr^m)$.

Let $\o$ be a bounded open subset of $\Rb^{n-1}$ and $I=(-1,1)$,
we define $\O:=\o \times I$. In the sequel, we will identify
$L^s(\o;\Rb^m)$ (resp. $W^{1,s}(\o;\Rb^m)$) with the space of
functions $v \in L^s(\O;\Rb^m)$ (resp. $W^{1,s}(\O;\Rb^m)$) such
that $D_n v = 0$ in the sense of distribution.

For every $(a,b)\subset \Rb$ with $a<b$ and $q_1,q_2\geq 1$,
$L^{q_1}(a,b;L^{q_2}(\Rb^{(n-1)};\Rb^m))$ is the space of
measurable $m$-vectorial functions $\zeta$ such that
$$
\int^a_b\left(\int_{\Rb^{n-1}}|\zeta(x_\a,x_n)|^{q_2}\,dx_\a\right)^{\frac{q_1}{q_2}}\,dx_n<+\infty.
$$
Let $a \in \Rb^{n-1}$ and $\rho>0$, we denote by $B^{n-1}_\rho(a)$
the open ball of $\Rb^{n-1}$ of center $a$ and radius $\rho$ and
by $Q^{n-1}_\rho(a)$ the open cube of $\Rb^{n-1}$ with center $a$
and length side $\rho$. We write $B^{n-1}_\rho$ instead of
$B^{n-1}_\rho(0)$ not to overburden notation. Let $x_i^\e= i\e$
with $i\in \Zb^{n-1}$, we set $Q_{i,\e}^{n-1}:= Q^{n-1}_\e
(x_i^\e)$.

We define $U^{+a} := U \times (0,a)$ and $U^{-a}:= U \times
(-a,0)$ with $U \subseteq \Rb^{n-1}$ and $a>0$, while if $a=1$,
then $U^{+} = U^{+1}$ and $U^-= U^{-1}$.

We set $C_{1,\infty}:=\left\{ (x_\a,0) \in \Rb^n : 1\leq |x_\a|
\right\}$ and $C_{1,N}:=\left\{ (x_\a,0) \in \Rb^n : 1\leq |x_\a|
< N \right\}$ for every $N>1$.

Let $p \geq 1$, we denote by ${\rm Cap}_p(B^{n-1}_1; A)$ the
$p$-capacity of $B^{n-1}_1$ with respect to $A \subseteq \Rb^d$:
$$
{\rm Cap}_p(B^{n-1}_1 ;A)=\inf \left\{ \int_A |D \psi |^p\, dx :
\; \psi \in W^{1,p}_0(A) \text{ and }\psi=1 \text{ on }B^{n-1}_1
\right\}.
$$
The letter $c$ will stand for a generic strictly-positive constant
which may vary from line to line and expression to expression
within the same formula.


\section{Statements of the main results}\label{mainres}

\noindent Since we are going to work with varying domains, we have
to precise the meaning of `converging sequences'.
\begin{defi}\label{conv} Let
$ \Omega_j= \o^{+\d_j} \cup \o^{-\d_j} \cup \big( \o_{r_j,\e_j}
\times \{0\} \big)$. Given a sequence $(u_j) \subset
W^{1,p}(\O_j;\Rb^m)$, we define ${\hat
u}_j(x_\a,x_n):=u_j(x_\a,\d_j\, x_n)$. We say that $(u_j)$
converges {\rm (}resp. converges weakly{\rm )} to $(u^+,u^-) \in
W^{1,p}(\o;\Rb^m) \times W^{1,p}(\o;\Rb^m)$ if we have
\begin{eqnarray*}
\hat u^+_j & := & \hat u_j|_{\o^+} \to u^+ \hbox{ in}\quad
L^p(\o^+;\Rb^m) \quad (\hbox{resp. weakly in
}W^{1,p}(\o^+;\Rb^m)),\\
\hat {u}^-_j & := & \hat u_j|_{\o^-} \to u^- \hbox{ in}\quad
L^p(\o^-;\Rb^m) \quad (\hbox{resp. weakly in
}W^{1,p}(\o^-;\Rb^m)).
\end{eqnarray*}
Similarly if we replace $\O_j$ by $\o^{\pm \d_j}$.

\smallskip
We say that the sequence $(|D u_j|^p / \d_j)$ is equi-integrable
on $\o^{\pm \d_j}$ if $\big( \big|\big(D_\a \hat u_j
|\frac{1}{\d_j} D_n \hat u_j \big) \big|^p \big)$ is
equi-integrable on $\o^\pm$.
\end{defi}

\begin{rmk}\label{convrmk}{\rm By
virtue of Definition \ref{conv}, a sequence $(u_j) \subset
W^{1,p}(\O_j;\Rb^m)$ converges to $(u^+,u^-)\in W^{1,p}(\o;\Rb^m)
\times W^{1,p}(\o;\Rb^m)$ if and only if
\begin{equation}\label{strongconv}
\lim_{j\to +\infty}\frac{1}{\d_j}\int_{\o^{\pm\d_j}}|u_j -
u^\pm|^p\, dx = 0,
\end{equation}
while (\ref{strongconv}) and
\begin{equation}\label{weakconv}
\sup_{j \in \Nb} \frac{1}{\d_j}\int_{\o^{\pm\d_j}}|D u_j|^p\,dx =
\sup_{j \in \Nb} \int_{\o^\pm}\left|\left(D_\a \hat
u_j\Big|\frac{1}{\d_j} D_n \hat u_j \right) \right|^p dx <+\infty
\end{equation}
imply the weak convergence.}\end{rmk}

Note that Remark \ref{convrmk} is still valid if we consider the
domain $ \o^{+\d_j} \cup\o^{-\d_j}$ in place of $\O_j$.

The main results of this paper are the following:

\begin{thm}[$\Gamma$-convergence]\label{ABZ}
Let $1<p<n-1$. Let $\o$ be a bounded open subset of $\Rb^{n-1}$
satisfying $\HH^{n-1}(\partial \o)=0$ and $W:\Rb^{m \times n} \to
[0,+\infty)$ be a Borel function such that $W(0)=0$ and satisfying
a growth condition of order $p$ : there exists a constant
$\beta>0$ such that
\begin{equation}\label{pgrowth}
|F|^p - 1 \leq W(F) \leq \beta (|F|^p + 1), \quad \text{for every
} F \in \Rb^{m \times n}.
\end{equation}
Let $(\e_j)$, $(\d_j)$ and $(r_j)$ be sequences of strictly
positive numbers converging to zero such that
$$
\lim_{j \to  +\infty}\frac{\d_j}{\e_j}=0
$$
and set
$$
\ell:= \lim_{j \to +\infty}\frac{r_j}{\d_j}\,.
$$
If
$$
\ell \in (0,+\infty] ,\quad \text{ and }\quad 0<
R^{(\ell)}:=\lim_{j \to
+\infty}\frac{r_j^{n-1-p}}{\e_j^{n-1}}<+\infty
$$
or
$$
\ell =0 ,\quad \text{ and } \quad 0< R^{(0)}:=\lim_{j \to
+\infty}\frac{r_j^{n-p}}{\d_j\, \e_j^{n-1}}<+\infty\,,
$$
then, up to an extraction, the sequence of functionals $\mathcal
F_j : L^p(\O_j;\Rb^m) \to [0,+\infty]$ defined by
$$
\mathcal F_j(u):=\left\{
\begin{array}{ll}
\ds \frac{1}{\d_j} \int_{\O_j} W(D
u)\, dx  & \text{if }u \in W^{1,p}(\O_j;\Rb^m),\\
&\\
\ds +\infty & \text{otherwise}
\end{array}
\right.
$$ $\G$-converges to
$$
\mathcal F^{(\ell)}(u^+,u^-)=\int_\o \Q_{n-1}\overline W(D_\a
u^+)\, dx_\a + \int_\o \Q_{n-1}\overline W(D_\a u^-)\, dx_\a +
R^{(\ell)} \int_\o \varphi^{(\ell)} (u^+ - u^-)\, dx_\a
$$
on $W^{1,p}(\o;\Rb^m) \times W^{1,p}(\o;\Rb^m)$ with respect to
the convergence introduced in Definition {\rm\ref{conv}}, where
 $\overline W(\overline F):=\inf \{ W(\overline
F|z) : z \in \Rb^m\}$, $\Q_{n-1}\overline W$ is the
$(n-1)$-quasiconvexification of $\overline W$ and
$\varphi^{(\ell)}:\Rb^m \to [0,+\infty)$ is a locally Lipschitz
continuous function for any $\ell\in [0,+\infty]$.
\end{thm}

\begin{rmk}\label{andrea}{\rm
Note that if $\ell\in (0,+\infty]$ the only meaningful scaling for
$r_j$ is that of order $\e_j^{(n-1)/(n-1-p)}$\ie for both
$R^{(\ell)}=0$ and $R^{(\ell)}=+\infty$ we loose the asymptotic
memory of the sieve. In fact, if $R^{(\ell)}=0$, we obtain two
uncoupled problems in the limit, while if $R^{(\ell)}=+\infty$,
limit deformations $(u^+, u^-)$ with finite energy are continuous
across the mid-section ($u^+=u^-$ in $\omega$) as in Le Dret-Raoult
\cite{LDR}. Similarly, for $\ell=0$.}
\end{rmk}

\smallskip

\begin{rmk}\label{equiv}{\rm
If $\ell \in (0,+\infty)$ then
$$
0<R^{(\ell)}=\lim_{j \to +\infty} \frac{r_j^{n-1-p}}{\e_j^{n-1}} <
+\infty \qquad \hbox{if\, and\, only\, if} \qquad 0<R^{(0)}=
\lim_{j \to +\infty} \frac{r_j^{n-p}}{\d_j\,\e_j^{n-1}} <
+\infty\,;
$$
hence, in this case the two meaningful scalings are equivalent.}
\end{rmk}

\noindent The following result provides a characterization of the
interfacial energy density $\varphi^{(\ell)}$ for each $\ell \in
[0,+\infty]$.

\begin{thm}[Representation formulas]\label{ABZ2}
Let $p^*=(n-1)p/(n-1-p)$ be the Sobolev exponent in dimension
$(n-1)$. Then, upon extracting a subsequence, there exists the
limit
$$g(F):=\lim_{j\to +\infty} r_j^p \Q_n W (r_j^{-1} F),$$
for all $F \in \Rb^{m \times n}$, where $\Q_n W$ denotes the
$n$-quasiconvexification of $W$, so that:\\
 if $\ell \in (0,+\infty)$,
\begin{eqnarray*}
\varphi^{(\ell)}(z) & := & \inf \Bigg\{\int_{(\Rb^{n-1} \times I)
\setminus C_{1,\infty}}g \big(D_\a \zeta | \ell D_n \zeta \big)\,
dx : \zeta \in W^{1,p}_{\rm loc}((\Rb^{n-1} \times I) \setminus
C_{1,\infty};\Rb^m),\\
&&\hspace{0.3cm} D \zeta \in L^p ((\Rb^{n-1} \times I) \setminus
C_{1,\infty};\Rb^{m \times n} ),
\;\quad \zeta - z \in L^p (0,1;L^{p^*}(\Rb^{n-1};\Rb^m))\\
&&\hspace{6.5 cm} \zeta \in L^p (-1,0;
L^{p^*}(\Rb^{n-1};\Rb^m))\Bigg\};
\end{eqnarray*}
if $\ell=+\infty$
\begin{eqnarray*}
\varphi^{(\infty)}(z) & := & \inf \Bigg\{\int_{\Rb^{n-1}}\Bigl(
\Q_{n-1} \,\overline g (D_\a \zeta^+ ) + \Q_{n-1} \, \overline g
(D_\a \zeta^- ) \Bigr)\, dx_\a : \; \zeta^\pm \in
W^{1,p}_{\rm loc}(\Rb^{n-1};\Rb^m),\\
&&\hspace{3.9 cm}\zeta^+=\zeta^- \text{ in }B_1^{n-1}, \quad D_\a
\zeta^\pm
\in L^p(\Rb^{n-1} ;\Rb^{m \times (n-1)}),\\
&&\hspace{6.7cm} (\zeta^+ - z) \,, \zeta^- \in
L^{p^*}(\Rb^{n-1};\Rb^m) \Bigg\}\,,
\end{eqnarray*} where $\overline g(\overline F):=\inf \{ g(\overline
F|z) : z \in \Rb^m\}$ and $\Q_{n-1}\overline g$ is the
$(n-1)$-quasiconvexification of $\overline g$;\\\smallskip if
$\ell=0$
\begin{eqnarray*}
\varphi^{(0)}(z) & = & \inf \Bigg\{\int_{\Rb^n \setminus
C_{1,\infty}} g (D \zeta) \, dx :\; \zeta \in W^{1,p}_{\rm
loc}(\Rb^n \setminus C_{1,\infty};\Rb^m),\; D \zeta \in L^p(\Rb^n
\setminus
C_{1,\infty};\Rb^{m \times n}),\\
&& \hspace{2.0cm}  \zeta-z \in
L^p(0,+\infty;L^{p^*}(\Rb^{n-1};\Rb^m)),\; \zeta \in
L^p(-\infty,0;L^{p^*}(\Rb^{n-1};\Rb^m)) \Bigg\}\,,
\end{eqnarray*}
for all $z\in \Rb^m$.
\end{thm}

\bigskip

\begin{rmk}{\rm
Without loss of generality we may assume that $W$ is quasiconvex
(upon first relaxing the energy); hence, by (\ref{pgrowth}), $W$
satisfies the following $p$-Lipschitz condition (see e.g.
\cite{D}):
\begin{equation}\label{plip}
|W(F_1) -W(F_2)| \leq c\, (1+ |F_1|^{p-1}+ |F_2|^{p-1})|F_1-F_2|,
\quad \text{ for all }F_1,\; F_2 \in \Rb^{m \times n}\,.
\end{equation}
}\end{rmk}

\section{Preliminary results}\label{prel}

\subsection{Some rescaled Poincar\'e Inequalities}

\noindent Since we deal with varying domains, depending on
different parameters, it is useful to note how the constant in
Poincar\'e type inequalities rescale with respect to such
parameters.

\begin{lemma}\label{poincare}
Let $A$ be an open bounded and connected subset of $\Rb^{n-1}$
with Lipschitz boundary and let $A_\rho:=\rho A$ for $\rho>0$.
\begin{itemize}
\item[(i)]There exists a constant $c>0$ (depending only on
$(A,n,p)$) such that for every $\rho,\d>0$
$$\int_{A^{\pm \d}_\rho}|u - \overline u_{A_\rho^{\pm \d}}|^p\, dx
\leq c  \int_{A_\rho^{\pm \d}}\left( \rho^p |D_\a u|^p + \d^p|D_n
u|^p \right)\, dx,
$$
for every $u \in W^{1,p}(A^{\pm\d}_\rho;\Rb^m)$ where $\overline
u_{A_\rho^{\pm \d}}= \medi_{A^{\pm \d}_\rho} u \dx$.

\item[(ii)]If $B$ is an open and connected subset of $A$ with Lipschitz
boundary and $B_\rho := \rho B$ then there exists a constant $c>0$
(depending only on $(A,B,n,p)$) such that for every $\rho,\d>0$
$$\int_{A^{\pm \d}_\rho}|u - \overline u_{B^{\pm \d}_\rho}|^p\, dx
\leq c  \int_{A^{\pm \d}_\rho}\left( \rho^p |D_\a u|^p + \d^p|D_n
u|^p \right)\, dx,
$$
for every $u \in W^{1,p}(A^{\pm \d}_\rho;\Rb^m)$ where $\overline
u_{B^{\pm \d}_\rho}= \medi_{B^{\pm \d}_\rho} u \dx$.
\end{itemize}
\end{lemma}

\noindent {\it Proof. } Let us define $v(x_\a,x_n):=u(\rho x_\a,
\d x_n)$ then $v \in W^{1,p}(A^\pm;\Rb^m)$. By a change of
variable, we get that $\overline u_{A_\rho^{\pm \d}}=\overline
v_{A^\pm}$. Moreover, by the Poincar\'e Inequality, there exists a
constant $c=c(A,n,p)>0$ such that
\begin{eqnarray*}
\int_{A_\rho^{\pm \d}}|u - \overline u_{A_\rho^{\pm\d}}|^p\, dx &
= & \d \rho^{n-1} \int_{A^\pm}|v -
\overline v_{A^\pm}|^p\, dy\\
 & \leq & c\d \rho^{n-1} \int_{A^\pm}|D v|^p\, dy\\
 & = & c \int_{A_\rho^{\pm \d}}\left(\rho^p |D_\a u|^p +
 \d^p |D_n u|^p \right)\, dx
\end{eqnarray*}
and it completes the proof of (i). Now, if $B_\rho \subset
A_\rho$, we get that
\begin{eqnarray*}
&&\int_{A_\rho^{\pm\d}}|u - \overline u_{B_\rho^{\pm\d}}|^p\, dx\\
& \leq & c\Bigl(\int_{A_\rho^{\pm\d}}|u - \overline
u_{A_\rho^{\pm\d}}|^p\, dx +\d \rho^{n-1} \HH^{n-1}(A)|\overline
u_{A_\rho^{\pm\d}}-\overline u_{B_\rho^{\pm\d}}|^p\Bigr)\\
& \leq & c\int_{A_\rho^{\pm\d}} |u - \overline
u_{A_\rho^{\pm\d}}|^p\, dx + c  \frac{\HH^{n-1}(A)}{\HH^{n-1}(B)}
\left( \int_{B_\rho^{\pm\d}}|u - \overline u_{A_\rho^{\pm\d}}|^p
\, dx + \int_{B_\rho^{\pm\d}} |u - \overline
u_{B_\rho^{\pm\d}}|^p\, dx\right)\\
& \leq & c \int_{A_\rho^{\pm\d}} \left( \rho^p |D_\a u|^p + \d^p
|D_n u|^p \right)\, dx.
\end{eqnarray*}
\hfill$\Box$\\

\subsection{A joining lemma on varying domains}
\noindent If not otherwise specified, in all what follows the
convergence of a sequence of functions has to be intended in the
sense of Definition \ref{conv}.

\smallskip

The following lemma, is the key tool in the proof of Theorem
\ref{ABZ}. It is a technical result which allows to modify
sequences of functions `near' the sets
$B^{(n-1)}_{r_j}(x_i^{\e_j})$. It is very close in spirit to Lemma
3.4 in \cite{Ans} although now the geometry of the problem yields
a different construction involving suitable cylindrical (instead
of spherical) annuli to surround the connecting zones.

\begin{lemma}\label{important}
Let $(\e_j)$, $(\d_j)$ be sequences of strictly positive numbers
converging to $0$ and such that $\d_j\ll \e_j$. Let $(u_j) \subset
W^{1,p}(\o^{+\d_j}\cup \o^{-\d_j};\Rb^m)$ be a sequence converging
to $(u^+,u^-) \in W^{1,p}(\o;\Rb^m) \times W^{1,p}(\o;\Rb^m)$
satisfying $\sup_j\mathcal{F}_j(u_j)<+\infty$; let $k\in \Nb$. Set
$\rho_j=\g\e_j$ with $\g<1/2$ and
$$ Z_j:=\{i \in \Zb^{n-1}: \; {\rm dist}(x_i^{\e_j},\Rb^{n-1}
\setminus \o)
> \e_j\}\,.
$$
For every $i \in Z_j$, there exists $k_i \in \{0,\ldots,k-1\}$
such that having set
$$
C_j^i := \left\{x_\a \in \o :\;  2^{-k_i-1}\rho_j <|x_\a -
x_i^{\e_j}| < 2^{-k_i}\rho_j\right\},$$
\begin{equation}\label{aver} u_j^{i \pm} :=
\med_{(C_j^i)^{\pm\d_j}} u_j\, dx
\end{equation}
and
$$
\rho_j^i := \frac{3}{4}2^{-k_i} \rho_j,$$ there exists a sequence
$(w_j)\subset W^{1,p}(\o^{+\d_j}\cup \o^{-\d_j};\Rb^m)$ weakly
converging to $(u^+,u^-)$ such that
\begin{equation}\label{wj1} w_j=u_j \; \text{ in }\;  \Bigl(\o
\setminus \bigcup_{i \in Z_j} C_j^i \Bigr)^{\pm\d_j},
\end{equation}
\begin{equation}\label{wj2}
w_j=u_j^{i\pm} \; \text{ on } \; \big(\partial
B^{n-1}_{\rho_j^i}(x_i^{\e_j})\big)^{\pm\d_j}
\end{equation} and satisfying
\begin{equation}\label{wj3}
\limsup_{j\to +\infty}\frac{1}{\d_j}\int_{\o^{\pm\d_j}} \big|  W(D
w_j) - W(D u_j) \big|\, dx \leq \frac{c}{k}\,.
\end{equation}
\end{lemma}

\noindent {\it Proof. }For every $j\in \Nb$, $i\in Z_j$, $k \in
\Nb$ and $h \in \{0,\ldots, k-1\}$, we define
$$
C_j^{i,h}:=\left\{ x_\a \in \o :\; 2^{-h-1}\rho_j <
|x_\a-x_i^{\e_j}| < 2^{-h}\rho_j \right\},
$$
$$(u_j^{i,h})^\pm := \med_{(C_j^{i,h})^{\pm\d_j}}u_j\, dx$$ and
\begin{equation}\label{def-rhojih}
\rho_j^{i,h}:=\frac{3}{4}2^{-h}\rho_j.
\end{equation}
Let $\phi \equiv \phi_j^{i,h} \in \C_c^\infty(C_j^{i,h};[0,1])$ be
a cut-off function such that $\phi=1$ on $\partial
B^{n-1}_{\rho_j^{i,h}}(x_i^{\e_j})$ and $|D_\a \phi| \leq
c/\rho_j^{i,h}$. In $(C_j^{i,h})^{\pm\d_j}$, we set
$$
w_j^{i,h}(x):=\phi(x_\a) (u_j^{i,h})^\pm + (1-\phi(x_\a)) u_j,
$$
then
\begin{eqnarray*}
\int_{(C_j^{i,h})^{\pm\d_j}}|D w_j^{i,h}|^p \dx & \leq & c
\int_{(C_j^{i,h})^{\pm\d_j}} \left(|D_\a
\phi|^p|u_j-(u_j^{i,h})^\pm|^p + |D u_j|^p \right) dx\\
& \leq & c\int_{(C_j^{i,h})^{\pm\d_j}}
\left(\frac{|u_j-(u_j^{i,h})^\pm|^p}{(\rho_j^{i,h})^p}
 + |D u_j|^p \right) dx.
\end{eqnarray*}
Applying Lemma \ref{poincare} (i), with $\rho=\rho_j^{i,h} $ and
$A_\rho=C_j^{i,h}$, we have that
\begin{eqnarray}\label{1615}
&&\int_{(C_j^{i,h})^{\pm\d_j}} |D w_j^{i,h}|^p \dx \nonumber\\
&\leq& c \int_{(C_j^{i,h})^{\pm\d_j}}\left( |D_\a u_j|^p +
\Bigl({\d_j\over\rho_j^{i,h}}\Bigr)^p |D_n u_j|^p \right) dx + c
\int_{(C_j^{i,h})^{\pm\d_j}}
|D u_j|^p \, dx\nonumber\\
&\leq& m_j(k,\g)\, c \int_{(C_j^{i,h})^{\pm\d_j}} |D u_j|^p \, dx,
\end{eqnarray}
where by (\ref{def-rhojih})
$$m_j(k,\g):= \max\Bigl\{1, \Bigl({2^{k+1}\over 3\g}\Bigr)^p
\Bigl({\d_j\over \e_j}\Bigr)^p\Bigr\}
$$
and since $\d_j \ll \e_j$, $m_j(k,\g)\to 1$ as $j \to +\infty$. As
$$\sum_{h=0}^{k-1}
\int_{(C_j^{i,h})^{\pm\d_j}}(1+|D u_j|^p)\, dx \leq
\int_{B^{n-1}_{\rho_j}(x_i^{\e_j})^{\pm\d_j}}(1+ |D u_j|)^p\ dx,$$
there exists $k_i \in \{0,\ldots,k-1\}$ such that, having set
$C_j^i:=C_j^{i,k_i}$, we get
\begin{equation}\label{1627}
\int_{(C_j^i)^{\pm\d_j}}(1+|D u_j|^p)\, dx \leq \frac{1}{k}
\int_{B^{n-1}_{\rho_j}(x_i^{\e_j})^{\pm\d_j}}(1+ |D u_j|^p)\, dx.
\end{equation}
\begin{figure}
\centerline{\psfig{file=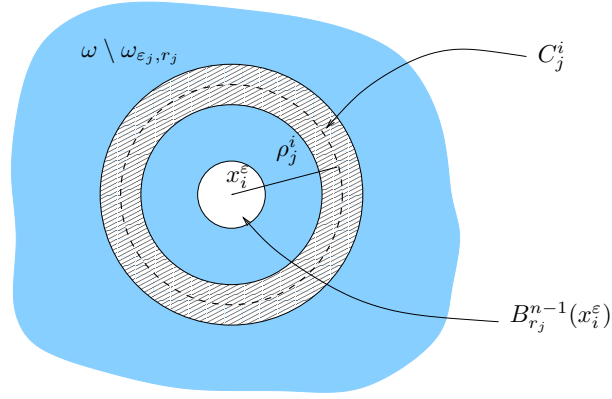,width=4in}} \vspace*{8pt}
\caption{The $(n-1)$-dimensional annuli $C^i_j$\,.}
\end{figure}
Hence, if we define the sequence
$$
w_j:=\left\{
\begin{array}{ll}
w_j^{i,k_i} & \text{in }
(C_j^i)^{\pm\d_j}\, \text{ for }\, i \in Z_j\\
&\\ u_j & \text{otherwise}\,,
\end{array}\right.
$$
by the $p$-growth condition (\ref{pgrowth}), (\ref{1615}),
(\ref{1627}) and Remark \ref{convrmk} we have
\begin{eqnarray*}
\frac{1}{\d_j}\int_{\o^{\pm\d_j}} \big|  W(D w_j) - W(D u_j)
\big|\, dx &=&  \sum_{i \in Z_j}
\frac{1}{\d_j}\int_{(C_j^i)^{\pm\d_j}} \big| W(D w_j^{i,k_i}) -
W(D u_j) \big|\,
dx\\
& \leq &{c\over k}\, m_j(k,\g) \sum_{i \in Z_j} \frac{1}{\d_j}
\int_{B^{n-1}_{\rho_j}(x_i^{\e_j})^{\pm\d_j}}(1+ |D u_j|^p)\,
dx\\
& \leq & \frac{c}{k}\, m_j(k,\g)  \left(1+\sup_{j \in \Nb}
\frac{1}{\d_j} \int_{\o^{\pm\d_j}} |Du_j|^p\dx \right)\\
&\le& \frac{c}{k}\, m_j(k,\g) \,,
\end{eqnarray*}
which concludes the proof of (\ref{wj3}). Note that, by
construction, $(w_j)$ satisfies (\ref{wj1}) and (\ref{wj2}) and it
converges weakly to $(u^+,u^-)$. In fact,
\begin{eqnarray*}
\frac{1}{\d_j}\int_{\o^{\pm\d_j}} |w_j -u^\pm|^p\, dx & = &
\frac{1}{\d_j}\sum_{i \in Z_j}\int_{ (C_j^i)^{\pm\d_j}} |\phi
u_j^{i\pm}
+(1-\phi)u_j - u^\pm|^p\, dx\\
&&+\frac{1}{\d_j}\int_{\o^{\pm\d_j} \setminus \bigcup_{i \in Z_j}
(C_j^i)^{\pm\d_j}} |u_j-u^\pm|^p\, dx\\
& \leq & \frac{c}{\d_j}\int_{\o^{\pm\d_j}} |u_j -u^\pm|^p\, dx
+\frac{c}{\d_j}\sum_{i \in Z_j}\int_{ (C_j^i)^{\pm\d_j}}| u_j -
u_j^{i\pm}|^p\, dx,
\end{eqnarray*}
while by Lemma \ref{poincare} (i) applied with $\rho=\rho_j^i$ and
since $\d_j \ll \e_j$, $\rho_j^i \leq \e_j$, we get
\begin{equation}\label{wconv1}
\frac{1}{\d_j}\int_{\o^{\pm\d_j}} |w_j -u^\pm|^p\, dx \leq
\frac{c}{\d_j}\int_{\o^{\pm\d_j}} |u_j -u^\pm|^p\, dx + c\e_j^p
\frac{1}{\d_j}\int_{\o^{\pm\d_j}} |D u_j|^p\, dx \,.
\end{equation}
Moreover by (\ref{1615}) we have
\begin{equation}\label{wconv2}
\frac{1}{\d_j}\int_{\o^{\pm\d_j}} |D w_j|^p\, dx \leq
\frac{c}{\d_j}\int_{\o^{\pm\d_j}} |D u_j|^p\, dx.
\end{equation}
Hence (\ref{wconv1}), (\ref{wconv2}), the convergence of $(u_j)$
towards $(u^+,u^-)$, $\sup_j\frac{1}{\d_j}\int_{\o^{\pm
\d_j}}|Du_j|^p\dx<+\infty$ together with Remark \ref{convrmk}
imply the weak convergence of $(w_j)$ towards
$(u^+,u^-)$. \hfill$\Box$\\

\begin{rmk}{\rm Note that to prove Lemma \ref{important}
we essentially use that $\rho_j< \e_j/2$ (but not necessarily
equal to $\g\e_j$) and $\lim_{j\to +\infty}(\d_j/\rho_j)=0$.
Hence, Lemma \ref{important} is still true if we replace the
assumptions $\d_j\ll \e_j$ and $\rho_j= \g\e_j$ by $\rho_j<
\e_j/2$ and $\lim_{j\to +\infty}(\d_j/\rho_j)=0$.

Since we will apply Lemma \ref{important} when $\rho_j= \g\e_j$
($\g<1/2$) and $\d_j\ll \e_j$, we prefer to prove it directly
under these assumptions. }\end{rmk}

\bigskip

If the sequence $(|D u_j|^p/\d_j)$ is equi-integrable on
$\o^{\pm\d_j}$ (see Definition \ref{conv}), then  we do not have
to choose for every $i \in Z_j$ a suitable annulus $C_j^i$ but we
may consider the same radius independently of $i$ as the following
lemma shows.

\begin{lemma}\label{important+equiint}
Let $(u_j)$, $(\e_j)$, $(\d_j)$, $(\rho_j)$ and $Z_j$ be as in
Lemma {\rm\ref{important}} and suppose that $(|D u_j|^p/\d_j)$ is
equi-integrable on $\o^{\pm\d_j}$. Set
$$
C_j^i := \left\{x_\a \in \o : \frac{2}{3}\rho_j <|x_\a -
x_i^{\e_j}| < \frac{4}{3} \rho_j\right\} \quad\text{ and }\quad
u_j^{i \pm} := \med_{(C_j^i)^{\pm\d_j}} u_j\, dx
$$
for every $i\in Z_j$. Then, there exists a sequence $(w_j) \subset
W^{1,p}(\o^{+\d_j}\cup \o^{-\d_j};\Rb^m)$ weakly converging to
$(u^+,u^-)$ such that
\begin{equation}\label{wj1bis} w_j=u_j \; \text{ in }\; \Bigl(\o
\setminus \bigcup_{i \in Z_j} C_j^{i}\Bigr)^{\pm\d_j},
\end{equation}
\begin{equation}\label{wj2bis}
w_j=u_j^{i\pm} \; \text{ on } \; \big( \partial
B^{n-1}_{\rho_j}(x_i^{\e_j}) \big) ^{\pm\d_j}
\end{equation} and
\begin{equation}\label{wj3bis}
\limsup_{j\to +\infty}\frac{1}{\d_j}\int_{\o^{\pm\d_j}} \big|  W(D
w_j) - W(D u_j) \big|\, dx \le o(1)\quad {\rm as}\quad \g \to
0^+\,.
\end{equation}
Moreover, the sequence $(|D w_j|^p/\d_j)$ is equi-integrable on
$\o^{\pm\d_j}$.
\end{lemma}

\noindent {\it Proof. }Let $\phi \equiv \phi_j^i \in
\C_c^\infty(C_j^i;[0,1])$ be a cut-off function such that $\phi=1$
on $\partial B^{n-1}_{\rho_j}(x_i^{\e_j})$ and $|D_\a \phi| \leq
c/\rho_j$. In $(C_j^i)^{\pm\d_j}$, we define
$$
w_j^{i}:=\phi(x_\a) u_j^{i\pm} + (1-\phi(x_\a)) u_j.
$$
Then, reasoning as in the proof of Lemma \ref{important}, we have
that
$$
\int_{(C_j^i)^{\pm\d_j}}W(D w_j^i)\, dx  \leq  c
\int_{(C_j^i)^{\pm\d_j}} (1+|D u_j|^p) \, dx.
$$
Hence, if we define
$$
w_j:=\left\{
\begin{array}{ll}
w_j^i & \text{in }
(C_j^i)^{\pm\d_j}\, \text{ for }\,i \in Z_j,\\
&\\ u_j & \text{otherwise},
\end{array}\right.$$
$w_j$ satisfies (\ref{wj1bis}) and (\ref{wj2bis}). Moreover,
\begin{eqnarray*}
\frac{1}{\d_j} \int_{\o^{\pm\d_j}} \big| W(D w_j) - W(D u_j)
\big|\, dx & \leq &  \sum_{i \in Z_j} \frac{1}{\d_j}
\int_{(C_j^i)^{\pm\d_j}} \big|
W(D w_j^i) - W(D u_j) \big|\,dx\\
& \leq & c \sum_{i \in Z_j} \frac{1}{\d_j}
\int_{(B^{n-1}_{4\rho_j/3}(x_i^{\e_j}) \cap \o)^{\pm\d_j}}(1+ |D
u_j|^p)\, dx.
\end{eqnarray*}
Since $\#(Z_j) \leq c/\e_j^{n-1}$, we get that
$$\mathcal H^{n-1}\Bigg( \bigcup_{i \in Z_j}
(B^{n-1}_{4\rho_j/3}(x_i^{\e_j}) \cap \o) \Bigg) \leq c \g^{n-1}$$
and by the equi-integrability of $(|D u_j|^p/\d_j)$ we obtain
(\ref{wj3bis}). Finally, the weak convergence of $(w_j)$ can be
proved as in Lemma \ref{important} while the equi-integrability of
$(|D w_j|^p/\d_j)$ is just a consequence of the definition of
$(w_j)$. \hfill$\Box$


\section{A preliminary analysis of the energy contribution `close' to
the connecting zones}\label{close} \noindent For later references,
in the following section we study the asymptotic behavior of a
sequence of functions which will turn out to represent the energy
contribution `close' to the connecting zones. The results listed
in this section will be applied in Section 6 to prove the
$\Gamma$-convergence of $(\mathcal{F}_j)$ as well as in Section 7
to compute the explicit formula for $\varphi^{(\ell)}$.
\bigskip

Before starting, let us recall that we consider the domain $
\Omega_j= \o^{+\d_j} \cup \o^{-\d_j} \cup \big( \o_{r_j,\e_j}
\times \{0\} \big)$ where $\o_{r_j,\e_j}:= \bigcup_{i
\in\Zb^{n-1}}B^{n-1}_{r_j}(x_i^{\e_j}) \cap \o$. Our
$\G$-convergence analysis deals with the case where the thickness
$\d_j$ of $\O_j$  is much smaller than the period of distribution
of the connecting zones $\e_j$\ie
$$
\lim_{j\to +\infty}\frac{\d_j}{\e_j}=0\,.
$$
Moreover, we can exclude that $r_j\ge \e_j/2$ otherwise the zones
may overlap. More precisely, we assume that $r_j\ll \e_j$\ie
\begin{equation}\label{e/r}
\lim_{j\to +\infty}\frac{r_j}{\e_j} = 0\,.
\end{equation}
This choice will be justify a posteriori since (\ref{e/r}) will be
the only admissible assumption to get a non trivial
$\Gamma$-convergence result (see Remark \ref{andrea}).

Finally, it remains to fix the behavior of $r_j$ with respect to
$\d_j$. Let us define
$$\ell:=\lim_{j \to +\infty}\frac{r_j}{\d_j}.$$
This yields to consider all the possible scenarii, namely to
distinguish between the cases: $\ell$ finite, infinite or zero.
\bigskip

For any fixed $\ell\in [0,+\infty]$, we consider the sequence of
functions $(\varphi^{(\ell)}_{\g,j})$ defined in (\ref{phigj}) and
(\ref{phigj12}). Propositions \ref{midi} and \ref{midi12}
establish the existence of the function $\varphi^{(\ell)}$ as the
(locally uniform) limit of $(\varphi^{(\ell)}_{\g,j})$ as $j\to
+\infty$ and $\g\to 0^+$ while Proposition \ref{rs} will allow us
to prove that $\varphi^{(\ell)}$ is actually the interfacial
energy density in $\mathcal{F}^{(\ell)}$ (see e.g. Proposition
\ref{gliminfclose}).

\subsection{The case $\ell \in (0,+\infty]$}

\noindent Setting $N_j= \e_j/r_j$, we define the space $$X^\g_j
(z) := \Big\{\zeta \in W^{1,p}((B^{n-1}_{\g N_j} \times I)
\setminus C_{1,\g N_j};\Rb^m):\;  \zeta=z \text{ on } (\partial
B^{n-1}_{\g N_j})^+ ,\,  \zeta=0 \text{ on } (\partial B^{n-1}_{\g
N_j})^-\Big\}\,,$$ where $I=(-1,1)$ and we consider the following
minimum problem
\begin{equation}\label{phigj}
\varphi^{(\ell)}_{\g,j}(z):=\inf \left\{ \int_{(B_{\g N_j}^{n-1}
\times I) \setminus C_{1,\g N_j}} r_j^p \, W\left(r_j^{-1} D_\a
\zeta | \d_j^{-1} D_n \zeta \right) dx : \quad \zeta \in
X^\g_j(z)\right\}\,.
\end{equation}
\begin{figure}
\centerline{\psfig{file=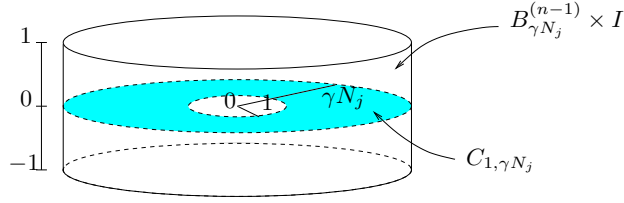,width=4in}} \vspace*{8pt}
\caption{The domain $(B^{(n-1)}_{\g N_j}\times I) \setminus C_{1,\g
N_j}$\,.}
\end{figure}
In the next proposition we study the behavior of
$(\varphi^{(\ell)}_{\g,j})$ as $j \to +\infty$ and $\g \to 0^+$.

\begin{proposition}\label{midi}
Let $\ell\in (0,+\infty]$. If
\begin{equation}\label{rinfty}0<R^{(\ell)}:=\lim_{j \to +\infty}
\frac{r_j^{n-1-p}}{\e_j^{n-1}} < +\infty
\end{equation}
then,

\smallskip

(i) there exists a constant $c>0$ (independent of $j$ and $\g$)
such that
$$
0 \leq \varphi^{(\ell)}_{\g,j}(z) \leq  c \left( |z|^p + \g^{n-1}
\right)
$$
for all $z \in \Rb^m$, $j\in \Nb$ and $\g>0$;

\smallskip

(ii) there exists a constant $c>0$ (independent of $j$ and $\g$)
such that
\begin{equation}\label{Stima1}
|\varphi^{(\ell)}_{\g,j}(z)-\varphi^{(\ell)}_{\g,j}(w)| \leq
c\,|z-w|\, \big(\g^{(n-1)(p-1)/p} + r_j^{p-1} + |z|^{p-1} +
|w|^{p-1} \big)
\end{equation}
for every $z,w \in \Rb^m$, $j\in \Nb$ and $\g>0$;

\smallskip

(iii) for every fixed $\gamma>0$, up to subsequences,
$\varphi^{(\ell)}_{\g,j}$ converges locally uniformly on $\Rb^m$
to $\varphi^{(\ell)}_{\g}$ as $j\to +\infty$ and
\begin{equation}\label{Stima12}
|\varphi^{(\ell)}_{\g}(z)-\varphi^{(\ell)}_{\g}(w)| \leq
c\,|z-w|\, \big(\g^{(n-1)(p-1)/p} + |z|^{p-1} + |w|^{p-1} \big)
\end{equation}
for every $z,w \in \Rb^m$ ;

\smallskip

(iv) up to subsequences, $\varphi^{(\ell)}_{\g}$ converges locally
uniformly on $\Rb^m$, as $\g \to 0^+$, to a continuous function
$\varphi^{(\ell)}:\Rb^m \to [0,+\infty)$ satisfying
\begin{equation}\label{AA}
0 \leq \varphi^{(\ell)} (z) \leq c |z|^p\,, \qquad
 \qquad |\varphi^{(\ell)}(z) - \varphi^{(\ell)}(w)|
 \leq c\,|z-w| \big( |z|^{p-1}+|w|^{p-1}\big)
\end{equation}
for every $z,w\in \rr^m$.
\end{proposition}

\noindent {\it Proof. } Fix $\g>0$, then $\g N_j> 2$ for $j$ large
enough.

(i) According to the $p$-growth condition (\ref{pgrowth}),
\begin{equation}\label{cap}
0 \leq \varphi^{(\ell)}_{\g,j}(z) \leq \beta \left(\C_{\g,j}(z) +
\HH^{n-1}(B_1^{n-1}) \g^{n-1}\frac{\e_j^{n-1}}{r_j^{n-1-p}}
\right),
\end{equation} where
$$
\C_{\g,j}(z):= \inf \left\{ \int_{(B_{\g N_j}^{n-1} \times I)
\setminus C_{1,\g N_j}} \left|\left(D_\a \zeta \Big|
\frac{r_j}{\d_j} D_n \zeta \right)\right|^p dx : \quad \zeta \in
X^\g_j(z)\right\}\,.
$$
Since $\C_{\g,j}(z)$ is invariant by rotations, reasoning as in
\cite{Ans} Section 4.1, we can consider the  minimization problem
with respect to a particular class of scalar test functions as
follows
\begin{eqnarray}\label{CgjStima}
\nonumber{\C_{\g,j}(z)\over |z|^p} & = &  \inf \Bigg\{
\int_{(B_{\g N_j}^{n-1} \times I) \setminus C_{1,\g N_j}}
\left|\left(D_\a \psi \Big| \frac{r_j}{\d_j} D_n \psi
\right)\right|^p \dx \,:\, \psi \in W^{1,p}((B^{n-1}_{\g
N_j} \times I) \setminus C_{1,\g N_j}),\\
&& \nonumber\hspace{5cm}\psi=1 \text{ on } \big(\partial
B^{n-1}_{\g N_j} \big)^+ \text{ and } \psi=0 \text{ on }
\big(\partial
B^{n-1}_{\g N_j}\big)^- \Bigg\}\\
 & \leq & \nonumber\inf \Bigg\{ \int_{B_{\g N_j}^{n-1}} \big( |D_\a
\psi^+ |^p + |D_\a \psi^- |^p \big) \dx : \quad (\psi^+
-1)\,,\,\psi^-
\in W^{1,p}_0(B^{n-1}_{\g N_j})\\
&& \hspace{7cm}\text{ and }\psi^+=\psi^-  \text{ in } B^{n-1}_{1}
\Bigg\}.
\end{eqnarray}
Let $\psi^\pm_1$ be the unique minimizer of the strictly convex
minimization problem (\ref{CgjStima}). It turns out that
$\psi^\pm_2:=1-\psi_1^\mp$ is also a minimizer. Thus by
uniqueness, $\psi_1^\pm=\psi_2^\pm$ and in particular,
$\psi^\pm_1=1/2$ in $B_1^{n-1}$. Hence,
\begin{eqnarray}\label{pascap}
\C_{\g,j}(z) & \leq & |z|^p \inf \Bigg\{ \int_{B_{\g N_j}^{n-1}}
\big( |D_\a \psi^+ |^p + |D_\a \psi^- |^p \big)\dxa \,: \, (\psi^+
-1)\,,\, \psi^- \in
W^{1,p}_0 (B^{n-1}_{\g N_j}),\nonumber\\
&& \hspace{7cm}\text{ and } \psi^+=\psi^-=\frac{1}{2}  \text{ in }
B^{n-1}_{1} \Bigg\}\nonumber\\
 & = & 2 |z|^p \inf \Bigg\{
\int_{B_{\g N_j}^{n-1}} |D_\a \psi |^p \dxa : \quad \psi \in
W^{1,p}_0(B^{n-1}_{\g N_j}) \text{ and } \psi=\frac{1}{2}  \text{
in } B^{n-1}_{1} \Bigg\}\nonumber\\
 & = &  \frac{|z|^p}{2^{p-1}} \inf \Bigg\{
\int_{B_{\g N_j}^{n-1}} |D_\a \psi |^p \dxa : \quad \psi \in
W^{1,p}_0(B^{n-1}_{\g N_j}) \text{ and } \psi=1  \text{ in }
B^{n-1}_{1} \Bigg\}\nonumber\\
& = & \frac{|z|^p}{2^{p-1}} {\rm Cap}_p \big(B^{n-1}_{1};B_{\g
N_j}^{n-1} \big).
\end{eqnarray}
Since
$$\lim_{j \to +\infty} {\rm Cap}_p \big(B^{n-1}_{1};B_{\g
N_j}^{n-1} \big) = {\rm Cap}_p \big(B^{n-1}_{1};\Rb^{n-1}\big) <
+\infty\,;
$$
hence, by (\ref{rinfty}), (\ref{cap}) and (\ref{pascap}) we
conclude the proof of (i).

\bigskip

(ii) For every $\eta>0$, there exists $\zeta_{\g,j} \in X^\g_j(z)$
such that
\begin{equation}\label{fabia}
\int_{(B_{\g N_j}^{n-1} \times I) \setminus C_{1,\g N_j}} r_j^p\,
W\left(r_j^{-1} D_\a \zeta_{\g,j} | \d_j^{-1} D_n \zeta_{\g,j}
\right) dx \leq \varphi^{(\ell)}_{\g,j}(z) +\eta.
\end{equation}
We want to modify $\zeta_{\g,j}$ in order to get an admissible
test function for $\varphi^{(\ell)}_{\g,j}(w)$. More precisely, we
just have to modify $\zeta_{\g,j}$ on a neighborhood of $(\partial
B_{\g N_j}^{n-1})^+$ to change the boundary condition $z$ into
$w$. To this aim we introduce a cut-off function $\theta \in
\C^\infty_c(\Rb^{n-1};[0,1])$, independent of $x_n$, such that
$$
\theta(x_\a)=\left\{
\begin{array}{rcl}
1 & \text{if} & x_\a \in B^{n-1}_1,\\
\\
 0 & \text{if} & x_\a \not\in B^{n-1}_2
\end{array}
\quad \text{ and }\quad |D_\a \theta|\leq c\,. \right.
$$
Hence, we define $\tilde\zeta_{\g,j} \in X^\g_j(w)$ as follows
$$
\tilde \zeta_{\g,j}=\left\{
\begin{array}{lcl}
\zeta_{\g,j}+ (1-\theta(x_\a)) (w-z)  & \text{in} &
(B^{n-1}_{\g N_j})^+\\
\\
\zeta_{\g,j}& \text{in} & (B^{n-1}_{\g N_j})^- \cup \big(
B_1^{n-1}\times \{0\} \big)\,.
\end{array}
\right.
$$
By (\ref{fabia}), since $\zeta_{\g,j}=\tilde \zeta_{\g,j}$ in
$(B^{n-1}_{\g N_j})^- $, we have that
\begin{eqnarray*}
&&\varphi^{(\ell)}_{\g,j}(w)-\varphi^{(\ell)}_{\g,j}(z)\\
 & \leq & r_j^p
\int_{(B_{\g N_j}^{n-1} \times I) \setminus C_{1,\g N_j}} \Bigl(
W\big(r_j^{-1} D_\a \tilde \zeta_{\g,j} | \d_j^{-1} D_n \tilde
\zeta_{\g,j} \big) - W\big(r_j^{-1} D_\a \zeta_{\g,j} | \d_j^{-1}
D_n
\zeta_{\g,j} \big)\Bigr) \dx\, +\eta\\
& = & r_j^p \int_{(B_{\g N_j}^{n-1})^+}\Bigl( W\big(r_j^{-1} D_\a
\tilde \zeta_{\g,j} | \d_j^{-1} D_n \tilde \zeta_{\g,j} \big)-
W\big(r_j^{-1} D_\a \zeta_{\g,j} | \d_j^{-1} D_n \zeta_{\g,j}
\big) \Bigl)\dx \,+\eta\,.
\end{eqnarray*}
By (\ref{plip}) and H\"older's Inequality, we obtain that
\begin{eqnarray*}
&&\varphi^{(\ell)}_{\g,j}(w)-\varphi^{(\ell)}_{\g,j}(z)\,-\eta\\
&\leq&  c\, \int_{(B_{\g N_j}^{n-1})^+}\Bigg( r_j^{p-1} + \left|
\left( D_\a \zeta_{\g,j} \Big| \frac{r_j}{\d_j} D_n \zeta_{\g,j}
\right) \right|^{p-1} + \left| \left( D_\a \tilde \zeta_{\g,j}
\Big| \frac{r_j}{\d_j} D_n \tilde
\zeta_{\g,j} \right) \right|^{p-1}\Bigg)\\
&& \hspace{2cm} \times \left| \left( D_\a \tilde \zeta_{\g,j}
-D_\a \zeta_{\g,j} \Big| \frac{r_j}{\d_j} (D_n \tilde \zeta_{\g,j}
- D_n \zeta_{\g,j})
\right) \right| \dx  \\
&\leq & c \int_{(B_{\g N_j}^{n-1})^+}\Bigg( r_j^{p-1} + 2 \left|
\left( D_\a \zeta_{\g,j} \Big| \frac{r_j}{\d_j} D_n \zeta_{\g,j}
\right) \right|^{p-1} + |D_\a \theta|^{p-1}\,
|w-z|^{p-1}\Bigg) |D_\a \theta|\, |w-z| \dx \\
&\leq& c\, |z-w|^p \int_{B_{\g N_j}^{n-1}} |D_\a \theta|^p\,
 dx_\a + c\, r_j^{p-1} |z-w| \int_{B_{\g N_j}^{n-1}} |D_\a \theta|\,
 dx_\a\\
 && + 2c\, |z-w| \; \|D_\a \theta \|_{L^p(B_{\g N_j}^{n-1};\Rb^{n-1})}
\left\| \left( D_\a \zeta_{\g,j} \Big| \frac{r_j}{\d_j} D_n
\zeta_{\g,j} \right) \right\|^{p-1}_{L^p\big((B_{\g
N_j}^{n-1})^+;\Rb^{m \times n}\big)}\,.
\end{eqnarray*}
Since $\g N_j > 2$ and ${\rm Supp}(\theta) \subset B_2^{n-1}$, we
obtain that
\begin{eqnarray}\label{pasta}
&&\hspace{-1cm}\varphi^{(\ell)}_{\g,j}(w)-\varphi^{(\ell)}_{\g,j}(z)\nonumber\\
&&\hspace{-0.5cm}\leq c |z-w| \Bigg( |z-w|^{p-1} + r_j^{p-1}+
\left\| \left( D_\a \zeta_{\g,j} \Big| \frac{r_j}{\d_j} D_n
\zeta_{\g,j} \right) \right\|^{p-1}_{L^p\big((B_{\g
N_j}^{n-1})^+;\Rb^{m \times n}\big)} \Bigg)\, +\eta.
\end{eqnarray}
By the $p$-growth condition (\ref{pgrowth}), (\ref{fabia}) and
(i), we have that
\begin{eqnarray}\label{pasta2}
&&\int_{(B_{\g N_j}^{n-1})^+} \left| \left( D_\a \zeta_{\g,j}
\Big| \frac{r_j}{\d_j} D_n \zeta_{\g,j} \right)
\right|^p \dx\nonumber\\
&\leq&\int_{(B_{\g N_j}^{n-1})^+} r_j^p \, W\left(r_j^{-1} D_\a
\zeta_{\g,j} | \d_j^{-1} D_n \zeta_{\g,j} \right) dx
+ r_j^p \, \HH^{n-1}\big( B_{\g N_j}^{n-1} \big)\nonumber\\
&\leq& \varphi^{(\ell)}_{\g,j}(z)
+\eta + c\g^{n-1} \frac{\e_j^{n-1}}{r_j^{n-1-p}}\nonumber\\
& \leq&  c(|z|^p + \g^{n-1} ) + \eta + c\g^{n-1}
\frac{\e_j^{n-1}}{r_j^{n-1-p}}.
\end{eqnarray}
Hence, by (\ref{pasta}), (\ref{pasta2}) and (\ref{rinfty}) we have
that
$$
\varphi^{(\ell)}_{\g,j}(w)-\varphi^{(\ell)}_{\g,j}(z) \leq c\,
|z-w| \Big( |z|^{p-1}+|w|^{p-1} + r_j^{p-1}+ \g^{(n-1)(p-1)/p} +
\eta^{(p-1)/p} \Big) +\eta
$$
and (\ref{Stima1}) follows by the arbitrariness of $\eta$.

\bigskip

By (ii) and Ascoli-Arzela's Theorem we have that, up to
subsequences, $\varphi^{(\ell)}_{\g,j}$ converges uniformly on
compact sets of $\rr^m$ to $\varphi^{(\ell)}_{\g}$ as $j\to
+\infty$. Moreover, passing to the limit in (\ref{Stima1}) as
$j\to +\infty$ we get
$$
|\varphi^{(\ell)}_{\g}(w)-\varphi^{(\ell)}_{\g}(z)| \leq c\, |z-w|
\Big( |z|^{p-1}+|w|^{p-1} + \g^{(n-1)(p-1)/p}\Big)\,.
$$
Hence, we can apply again Ascoli-Arzela's Theorem to conclude
that, up to subsequences, $\varphi^{(\ell)}_{\g}$ converges
uniformly on compact sets of $\rr^m$ to $\varphi^{(\ell)}$ as
$\g\to 0^+$. In particular, $\varphi^{(\ell)} :\Rb^m \to
[0,+\infty)$ is a continuous function and
$$
0 \leq \varphi^{(\ell)} (z) \leq c |z|^p\,, \qquad
 \qquad |\varphi^{(\ell)}(z) - \varphi^{(\ell)}(w)|
 \leq c\,\big( |z|^{p-1}+|w|^{p-1}\big)|z-w|
$$
for every $z,w\in \rr^m$. \hfill$\Box$


\begin{figure}
\centerline{\psfig{file=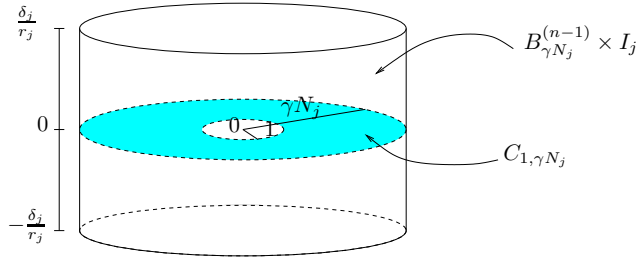,width=4in}} \vspace*{8pt}
\caption{The domain $(B^{(n-1)}_{\g N_j}\times I_j)\setminus C_{1,\g
N_j}$\,.}
\end{figure}

\subsection{The case $\ell =0$}

\noindent In this case we expect that the energy contribution due
to the presence of the sieve is obtained studying the behavior, as
$j\to +\infty$ and $\g\to 0^+$, of the sequence
($\varphi^{(0)}_{\g,j}$) defined as follows
\begin{eqnarray}\label{phigj12}
\varphi^{(0)}_{\g,j}(z) & := & \frac{\d_j}{r_j} \inf \left\{
\int_{(B_{\g N_j}^{n-1} \times I) \setminus C_{1,\g N_j}} r_j^p \,
W\left(r_j^{-1} D_\a \zeta | \d_j^{-1} D_n
\zeta \right)dx : \quad \zeta \in X^\g_j(z)\right\}\nonumber\\
 & = & \inf \left\{ \int_{(B_{\g N_j}^{n-1} \times
I_j ) \setminus C_{1,\g N_j}} r_j^p \, W(r_j^{-1} D \zeta )\, dx :
\quad \zeta \in Y^\g_j(z)\right\}
\end{eqnarray}
where $I_j:=(-\d_j/r_j , \d_j/r_j)$ and
\begin{eqnarray*}
Y^\g_j(z) = \Big\{\zeta \in W^{1,p}((B_{\g N_j}^{n-1} \times I_j)
\setminus C_{1,\g N_j} ;\Rb^m ) & : & \zeta=z \text{ on }(\partial
B_{\g N_j}^{n-1})^{+(\d_j/r_j)},\\
&& \zeta=0 \text{ on } (\partial B_{\g N_j}^{n-1})^{-(\d_j/r_j)}
\Big\}.
\end{eqnarray*}
Note that in this case we are interested in the limit behavior of
a sequence that is obtained from the one corresponding to
$\ell\in(0,+\infty]$ multiplying it by $\d_j/r_j$ (see
(\ref{phigj12}) and recall (\ref{phigj})). Let us try to motivate
this choice.

Let $\ell\in(0,+\infty)$, then starting from $(\ref{phigj})$ by a
change of variable it is immediate to check that
\begin{equation}\label{uguale}
\varphi_{\g,j}^{(\ell)}(z)=\frac{r_j}{\d_j}\inf \left\{
\int_{(B_{\g N_j}^{n-1} \times I_j ) \setminus C_{1,\g N_j}} r_j^p
\, W(r_j^{-1} D \zeta )\, dx : \quad \zeta \in Y^\g_j(z)\right\}.
\end{equation}
Now assuming that
$\lim_{j\to+\infty}r_j^{n-p}/(\d_j\,\e_j^{n-1})<+\infty$ (or
equivalently that
$\lim_{j\to+\infty}r_j^{n-1-p}/\e_j^{n-1}<+\infty$; see Remark
\ref{equiv}) we know that the sequence $(\varphi_{\g,j}^{(\ell)})$
converges to $\ell\, \tilde \varphi^{(\ell)}$, for some $\tilde
\varphi^{(\ell)}$, locally uniformly in $\Rb^m$, as $j\to+\infty$
and $\g\to 0^+$ (Proposition \ref{midi}). Then if
$\ell\in(0,+\infty)$, studying the limit behavior of
(\ref{phigj12}) is perfectly equivalent to study the limit
behavior of (\ref{phigj}). While if
$\ell=\lim_{j\to+\infty}r_j/\d_j=0$, (\ref{uguale}) suggests that,
to recover nontrivial information in the limit, we have to study
the asymptotic behavior of the sequence obtained from
(\ref{uguale}) dividing it by $r_j/\d_j$, that is to study the
asymptotic behavior of the sequence given by (\ref{phigj12}).
\bigskip

Following the line of the proof of Proposition \ref{midi}, we want
to establish an analogous result for the sequence
$(\varphi^{(0)}_{\g,j})$.

\begin{proposition}\label{midi12}
Let $\ell=0$. If
\begin{equation}\label{rzero}
0<R^{(0)}=\lim_{j \to +\infty} \frac{r_j^{n-p}}{\e_j^{n-1} \d_j} <
+\infty
\end{equation}
then,

\smallskip

(i) there exists a constant $c>0$ (independent of $j$ and $\g$)
such that
$$
0 \leq \varphi^{(0)}_{\g,j}(z) \leq  c \left( |z|^p + \g^{n-1}
\right)
$$
for all $z \in \Rb^m$, $j\in \Nb$ and $\g>0$;

\smallskip

(ii) there exists a constant $c>0$ (independent of $j$ and $\g$)
such that
\begin{equation}\label{Stima2}
|\varphi^{(0)}_{\g,j}(z)-\varphi^{(0)}_{\g,j}(w)| \leq c\, |z-w|
\,\big(\g^{(n-1)(p-1)/p} + r_j^{n-1} + |z|^{p-1} + |w|^{p-1} \big)
\end{equation}
for every $z,w \in \Rb^m$, $j\in \Nb$ and $\g>0$;

\smallskip

(iii) for every fixed $\gamma>0$, up to subsequences,
$\varphi^{(0)}_{\g,j}$ converges locally uniformly in $\rr^m$ to
$\varphi^{(0)}_{\g}$ as $j\to +\infty$, and
\begin{equation}\label{Stima23}
|\varphi^{(0)}_{\g}(z)-\varphi^{(0)}_{\g}(w)| \leq c\, |z-w| \Big(
\g^{(n-1)(p-1)/p} + |z|^{p-1}+|w|^{p-1}\Big)
\end{equation}
for every $z,w \in \Rb^m$;

\smallskip

(iv) up to subsequences, $\varphi^{(0)}_{\g}$ converges locally
uniformly in $\rr^m$, as $\g \to 0^+$, to a continuous function
$\varphi^{(0)}:\Rb^m \to [0,+\infty)$ satisfying
\begin{equation}\label{AA12}
0 \leq \varphi^{(0)} (z) \leq c |z|^p , \qquad
 \qquad |\varphi^{(0)}(z) - \varphi^{(0)}(w)|
 \leq c\,|z-w|\, \big( |z|^{p-1}+|w|^{p-1}\big)
\end{equation}
for every $z,w\in \rr^m$.

\end{proposition}
\noindent {\it Proof. } Fix $\g>0$, then $\g N_j> 2$ and
$\d_j/r_j>2$ for $j$ large enough.

(i) According to the $p$-growth condition (\ref{pgrowth}),
\begin{equation}\label{cap12}
0 \leq \varphi^{(0)}_{\g,j}(z) \leq \beta \left( \C_{\g,j}(z) + 2
\HH^{n-1}(B_1^{n-1})\, \g^{n-1}\,\frac{
\d_j\,\e_j^{n-1}}{r_j^{n-p}} \right),
\end{equation}
where
$$
\C_{\g,j}(z)= \inf \left\{ \int_{(B_{\g N_j}^{n-1} \times I_j)
\setminus C_{1,\g N_j}} |D \zeta |^p \, dx : \quad \zeta \in
Y^\g_j(z) \right\}.
$$
Arguing similarly than in the proof of Proposition \ref{midi}, we
can rewrite
\begin{eqnarray}\label{CgjScalar}
\nonumber{\C_{\g,j}(z)\over |z|^p} & = & \inf \Bigg\{ \int_{(B_{\g
N_j}^{n-1} \times I_j) \setminus C_{1,\g N_j}} |D \psi|^p\, dx
\,:\; \psi \in W^{1,p}((B^{n-1}_{\g N_j}
\times I_j) \setminus C_{1,\g N_j}),\nonumber\\
&& \hspace{1.6cm} \psi=1 \text{ on } (\partial B^{n-1}_{\g
N_j})^{+(\d_j/r_j)} \, , \quad \psi=0 \text{ on } (\partial
B^{n-1}_{\g N_j})^{-(\d_j/r_j)} \Bigg\}.
\end{eqnarray}
Let $\psi_1$ be the unique minimizer of the strictly convex
minimization problem (\ref{CgjScalar}). It turns out that
$\psi_2(x_\a,x_n):=1-\psi_1(x_\a,-x_n)$ is also a minimizer. Thus
by uniqueness, $\psi_1=\psi_2$ and in particular,
$\psi_1=\psi_2=1/2$ on $B_1^{n-1} \times \{0\}$. Thus
\begin{eqnarray}\label{pascap12}
\C_{\g,j}(z) & = & 2 |z|^p \inf \Bigg\{ \int_{(B_{\g
N_j}^{n-1})^{+(\d_j/r_j)}} |D \psi |^p \, dx : \quad \psi \in
W^{1,p}((B_{\g N_j}^{n-1})^{+(\d_j/r_j)}),\nonumber\\
&&\hspace{2cm}\psi = 0 \text{ on } (\partial B_{\g
N_j}^{n-1})^{+(\d_j/r_j)}  \text{ and } \psi=\frac{1}{2} \text{ on
}
B^{n-1}_1 \times \{0\} \Bigg\}\nonumber\\
 & = &  \frac{|z|^p}{2^{p-1}} \inf \Bigg\{
\int_{(B_{\g N_j}^{n-1})^{+(\d_j/r_j)}} |D \psi |^p \, dx :
\quad \psi \in W^{1,p}((B^{n-1}_{\g N_j})^{+(\d_j/r_j)}),\nonumber\\
&&\hspace{2cm}\psi = 0 \text{ on } (\partial B_{\g
N_j}^{n-1})^{+(\d_j/r_j)}  \text{ and } \psi=1 \text{ on }
B^{n-1}_{1} \times \{0\} \Bigg\}\nonumber\\
& \le & \frac{|z|^p}{2^p} {\rm Cap}_p \big(B^{n-1}_{1};B_{\g
N_j}^{n-1}\times I_j \big)\,.
\end{eqnarray}
Since
$$\lim_{j \to +\infty} {\rm Cap}_p
\big(B^{n-1}_{1};B_{\g N_j}^{n-1}\times I_j \big) = {\rm Cap}_p
\big(B^{n-1}_{1};\Rb^n\big)<+\infty\,;
$$
hence, by (\ref{rzero}), (\ref{cap12}) and (\ref{pascap12}) we
conclude the proof of (i).

\bigskip

(ii) We can proceed as in the proof of Proposition \ref{midi} (ii)
using a different cut-off function also depending on $x_n$.
Namely, let $\theta \in \C^\infty_c(\Rb^n;[0,1])$ be such that
$$\theta(x_\a,x_n)=\left\{
\begin{array}{rcl}
1 & \text{if} & (x_\a,x_n) \in B^{n-1}_1 \times (-1,1),\\
\\
 0 & \text{if} & (x_\a,x_n) \not\in B^{n-1}_2 \times (-2,2)
\end{array}
\quad \text{ and }\quad |D \theta|\leq c. \right.
$$
Hence, if $\zeta_{\g,j} \in Y_j^\g(z)$ is a sequence which `almost
attains' the infimum value $\varphi_{\g,j}^{(0)}$, we define
$\tilde\zeta_{\g,j} \in Y^\g_j(w)$ as follows
$$
\tilde \zeta_{\g,j}=\left\{
\begin{array}{lcl}
\zeta_{\g,j}+ (1-\theta(x)) (w-z)  & \text{in} &
(B^{n-1}_{\g N_j})^{+(\d_j/r_j)},\\
\\
\zeta_{\g,j}& \text{in} & \big( (B^{n-1}_{\g N_j})^{-(\d_j/r_j)}
\big) \cup \big( B_1^{n-1} \times \{0\} \big).
\end{array}
\right.
$$
By (\ref{rzero}) we conclude the proof of (ii) reasoning as in the
proof of Proposition \ref{midi} (ii).

The proof of (iii) and (iv) follows the line of the proof of (iii)
and (iv) in Proposition \ref{midi}. \hfill$\Box$

\smallskip

Now we are able to describe the energy contribution close to the
connecting zones as $j\to +\infty$ and $\g\to 0^+$.

\begin{proposition}[Discrete approximation of the interfacial
energy]\label{rs} Let $(u_j)\subset W^{1,p}(\O_j;\Rb^m)\cap
L^{\infty}(\O_j;\Rb^m)$ be a sequence converging to $(u^+,u^-)\in
W^{1,p}(\o;\Rb^m)\times
 W^{1,p}(\o;\Rb^m)$ such that\, $\sup_j\mathcal{F}_j(u_j)<+\infty$ and
satisfying $\sup_{j \in \Nb} \|u_j\|_{L^\infty(\O_j;\Rb^m)}
<+\infty$. Let $(u_j^{i\pm})$ be as in {\rm (\ref{aver})}. If
$$\ell\in(0,+\infty]\quad \text{ and }\quad
0<R^{(\ell)}=\lim_{j\to+\infty}\frac{r_j^{n-1-p}}{\e_j^{n-1}}<+\infty$$
or $$\ell=0\quad \text{ and }\quad
0<R^{(0)}=\lim_{j\to+\infty}\frac{r_j^{n-p}}{\d_j\e_j^{n-1}}<+\infty
$$
then
\begin{equation}\label{conv+-}
\lim_{\g\to 0^+}\limsup_{j \to +\infty}\int_\o \Bigl|\sum_{i \in
Z_j}\varphi^{(\ell)}_{\g,j}(u_j^{i+}-
u_j^{i-})\chi_{Q_{i,\e_j}^{n-1}} -
\varphi^{(\ell)}(u^+-u^-)\Bigr|\,dx_\a=0\,,
\end{equation}
for every $\ell\in[0,+\infty]$.
\end{proposition}

\noindent {\it Proof. } Since $\sup_{j \in \Nb}
\|u_j\|_{L^\infty(\O_j;\Rb^m)}<+\infty$ by Propositions \ref{midi}
or \ref{midi12} we have that
\begin{eqnarray*}
&&\limsup_{j\to +\infty} \int_\o \Bigl|\sum_{i \in
Z_j}\varphi^{(\ell)}_{\g,j}(u_j^{i+}-
u_j^{i-})\chi_{Q_{i,\e_j}^{n-1}} -
\varphi^{(\ell)}(u^+-u^-)\Bigr|\,dx_\a\\
&\le& \limsup_{j\to +\infty} \int_\o \sum_{i \in
Z_j}\Bigl|\varphi^{(\ell)}_{\g,j}(u_j^{i+}- u_j^{i-}) -
\varphi^{(\ell)}(u_j^{i+}-u_j^{i-})\Bigr|\chi_{Q_{i,\e_j}^{n-1}}\,dx_\a\\
&&+ \limsup_{j\to +\infty} \int_\o \Bigl|\sum_{i \in
Z_j}\varphi^{(\ell)}(u_j^{i+}- u_j^{i-})\chi_{Q_{i,\e_j}^{n-1}} -
\varphi^{(\ell)}(u^+-u^-)\Bigr|\,dx_\a\\
&\le& o(1) + \limsup_{j\to +\infty} \int_\o \Bigl|\sum_{i \in
Z_j}\varphi^{(\ell)}(u_j^{i+}- u_j^{i-})\chi_{Q_{i,\e_j}^{n-1}} -
\varphi^{(\ell)}(u^+-u^-)\Bigr|\,dx_\a, \end{eqnarray*} as $\g \to
0^+$. By (\ref{AA}) or (\ref{AA12}) and H\"older's Inequality we
have that
\begin{eqnarray*}
&&\limsup_{j\to +\infty}\int_\o \Bigl|\sum_{i \in
Z_j}\varphi^{(\ell)}(u_j^{i+}- u_j^{i-}) \chi_{Q^{n-1}_{i,\e_j}}
- \varphi^{(\ell)}(u^+-u^-) \Bigr|\, dx_\a\\
&=& \limsup_{j\to +\infty}\sum_{i \in Z_j} \int_{Q^{n-1}_{i,\e_j}}
|\varphi^{(\ell)}(u_j^{i+}- u_j^{i-})
- \varphi^{(\ell)}(u^+-u^-)|\, dx_\a\\
&\leq& c\, \limsup_{j\to +\infty}\Bigl(\sum_{i \in
Z_j}\int_{Q^{n-1}_{i,\e_j}}\big|u_j^{i+}-u^+|^p + |u_j^{i-}
-u^{-}|^p\dxa\Bigr)^{1/p}\,.
\end{eqnarray*}
Hence, it remains to prove that
\begin{equation}\label{lim0}
\limsup_{j\to +\infty}\sum_{i \in Z_j}\int_{Q^{n-1}_{i,\e_j}}
|u^{\pm}- u_j^{i\pm}|^p \, dx_\a =0\,.
\end{equation}
By Lemma \ref{poincare} (ii) applied with $\rho=\e_j$,
$B_\rho=C_j^i$ and $A_\rho=Q^{n-1}_{i,\e_j}$ and since $\d_j \ll
\e_j$, we have
\begin{eqnarray}\label{corep}
\int_{Q^{n-1}_{i,\e_j}} |u^\pm - u_j^{i\pm}|^p\, dx_\a &\leq&
{c\over \d_j}\Bigl(\int_{(Q^{n-1}_{i,\e_j})^{\pm\d_j}}
|u_j-u^{\pm}|^p\dx + \int_{(Q^{n-1}_{i,\e_j})^{\pm\d_j}}
|u_j-u_j^{i\pm}|^p\dx\Bigr)\nonumber
\\
&\le& {c\over \d_j}\int_{(Q^{n-1}_{i,\e_j})^{\pm\d_j}}
|u_j-u^{\pm}|^p\dx
+\frac{c\,\e_j^p}{\d_j}\int_{(Q^{n-1}_{i,\e_j})^{\pm\d_j}} |D
u_j|^p \dx\,,
\end{eqnarray}
for all $i \in Z_j$; hence, summing up on $i \in Z_j$, we find
$$
\sum_{i \in Z_j} \int_{Q^{n-1}_{i,\e_j}} |u_j - u_j^{i\pm}|^p\,
dx_\a \leq {c\over \d_j}\int_{\o^{\pm\d_j}} |u_j-u^{\pm}|^p\dx +
\frac{c\,\e_j^p}{\d_j}\int_{\o^{\pm\d_j}} |D u_j|^p \dx\,,
$$
then passing to the limit as $j\to +\infty$ by the convergence of
$(u_j)$ towards $(u^+,u^-)$ and $\sup_j\mathcal{F}_j(u_j)<+\infty$
we get (\ref{lim0}) and then (\ref{conv+-}). \hfill$\Box$


\section{$\Gamma$-convergence result}\label{proof}

\subsection{The liminf inequality}

\noindent Let $(u_j) \subset W^{1,p}(\O_j;\Rb^m)\cap
L^\infty(\O_j;\Rb^m)$ be a sequence converging to $(u^+,u^-)\in
W^{1,p}(\o,\rr^m)\times W^{1,p}(\o,\rr^m)$  such that $\sup_{j \in
\Nb} \|u_j\|_{L^\infty(\O_j;\Rb^m)}<+\infty$ and
$$
\liminf_{j \to +\infty}\mathcal F_j(u_j)<+\infty\,.
$$
By Lemma \ref{important}, for every fixed $k \in \Nb$, there
exists a sequence $(w_j)\subset W^{1,p}(\O_j;\Rb^m)\cap
L^\infty(\O_j;\Rb^m)$ weakly converging to $(u^+,u^-)$ satisfying
(\ref{wj1}), (\ref{wj2}) and such that
\begin{eqnarray}\label{lowerbd}
&&\liminf_{j \to +\infty}\frac{1}{\d_j} \left( \int_{\o^{+\d_j}}
W(D u_j)\, dx
+ \int_{\o^{-\d_j}} W(D u_j)\, dx \right)\nonumber\\
&\geq& \liminf_{j \to +\infty}\frac{1}{\d_j} \left(
\int_{\o^{+\d_j}} W(D w_j)\, dx +  \int_{\o^{-\d_j}} W(D
w_j)\, dx \right) - \frac{c}{k}\nonumber\\
&\geq& \liminf_{j \to +\infty}\frac{1}{\d_j} \left( \int_{(\o
\setminus E_j)^{+\d_j}} W(D w_j)\, dx + \int_{(\o
\setminus E_j)^{-\d_j}} W(D w_j)\, dx \right)\nonumber\\
&& + \liminf_{j \to +\infty}\frac{1}{\d_j} \left(
\int_{E_j^{+\d_j}} W(D w_j)\, dx +  \int_{E_j^{-\d_j}} W(D w_j)\,
dx \right) - \frac{c}{k},
\end{eqnarray}
where $E_j:=\bigcup_{i \in Z_j}B^{n-1}_{\rho_j^i}(x_i^{\e_j})$.

We first consider the energy contribution `far' from the
connecting zones. In this case, we suitably modify the sequence
$(w_j)$ in order to get a constant inside each half cylinder
$B^{(n-1)}_{\rho_j^i}(x_i^{\e_j})^{\pm \d_j}$. Then, we apply the
classical result of dimensional reduction proved in \cite{LDR} to
$\o^{+\d_j}$ and $\o^{-\d_j}$, separately.

\begin{proposition}\label{gliminffar}
We have
\begin{eqnarray*}
&&\liminf_{j \to +\infty}\frac{1}{\d_j} \left( \int_{(\o \setminus
E_j)^{+\d_j}} W(D w_j)\, dx + \int_{(\o \setminus E_j)^{-\d_j}}
W(D
w_j)\, dx \right)\\
&& \geq \int_\o \left(\Q_{n-1}\overline W (D_\a u^+) +
\Q_{n-1}\overline W( D_\a u^-) \right) \, dx_\a.
\end{eqnarray*}
\end{proposition}

\noindent {\it Proof. } We define
\begin{equation}\label{defvj}
v_j:=\left\{
\begin{array}{lcl}
w_j & \text{ in } & (\o \setminus E_j)^{\pm\d_j},\\
u_j^{i\pm} & \text{ in } &
B^{n-1}_{\rho_j^i}(x_i^{\e_j})^{\pm\d_j} \text{ if } i \in Z_j.
\end{array}
\right.
\end{equation}
Then $(v_j) \subset W^{1,p}(\O_j;\Rb^m)$ converges weakly to
$(u^+,u^-)$. In fact,
\begin{equation}\label{vjpm1}
\sup_{j \in \Nb} \frac{1}{\d_j}\int_{\o^{\pm\d_j}}|D v_j|^p\, dx
\leq \sup_{j \in \Nb} \frac{1}{\d_j}\int_{\o^{\pm\d_j}}|D u_j|^p\,
dx <+\infty.
\end{equation}
Moreover, since $\rho_j^i < \rho_j <\e_j/2$, then
$B^{n-1}_{\rho_j^i}(x_i^{\e_j}) \subset Q^{n-1}_{i,\e_j}$; hence,
$$\int_{\o^{\pm\d_j}}|v_j - u^\pm|^p\, dx  \leq  \int_{(\o \setminus
E_j)^{\pm\d_j}}|w_j - u^\pm|^p\, dx\nonumber + \sum_{i \in
Z_j}\int_{(Q^{n-1}_{i,\e_j})^{\pm\d_j}}|u^\pm - u_j^{i\pm}|^p\, dx
$$
and, by (\ref{corep}), we obtain that
\begin{eqnarray}\label{vjpm2}
\frac{1}{\d_j}\int_{\o^{\pm\d_j}}|v_j - u^\pm|^p\, dx & \leq &
\frac{1}{\d_j} \int_{\o^{\pm\d_j}}|w_j - u^\pm|^p\, dx +
\frac{c}{\d_j}\int_{\o^{\pm\d_j}}|u_j -
u^\pm|^p\, dx\nonumber\\
&&+  c \, \e_j^p \,\sup_{j \in \Nb}
\frac{1}{\d_j}\int_{\o^{\pm\d_j}} |D u_j|^p\, dx.
\end{eqnarray}
Passing to the limit as $j\to +\infty$ in (\ref{vjpm2}), by
(\ref{vjpm1}) and Remark \ref{convrmk} we get that $(v_j)$
converges weakly to $(u^+,u^-)$.

Since $W(0)=0$, by (\ref{defvj}) and \cite{LDR} Theorem 2, we have
\begin{eqnarray*}
&&\liminf_{j \to +\infty}\frac{1}{\d_j}\left( \int_{(\o \setminus
E_j)^{+\d_j}} W(D w_j)\, dx + \int_{(\o
\setminus E_j)^{-\d_j}} W(D w_j)\, dx \right)\\
&=&\liminf_{j \to +\infty}\frac{1}{\d_j}\left( \int_{(\o \setminus
E_j)^{+\d_j}} W(D v_j)\, dx + \int_{(\o
\setminus E_j)^{-\d_j}} W(D v_j)\, dx \right)\\
&=&\liminf_{j \to +\infty}\frac{1}{\d_j}\left( \int_{\o^{+\d_j}}
W(D v_j)\, dx + \int_{\o^{-\d_j}} W(D v_j)\, dx
\right)\\
&\geq& \int_\o \Q_{n-1}\overline W (D_\a u^+)\, dx_\a +\int_\o
\Q_{n-1} \overline W(D_\a u^-)\, dx_\a\,.
\end{eqnarray*}
\hfill$\Box$\\

Now let us deal with the contribution `near' the connecting zones.
We always work under the assumption
$$\ell\in(0,+\infty]\quad \text{ and }\quad
0<R^{(\ell)}=\lim_{j\to+\infty}\frac{r_j^{(n-1-p)}}{\e_j^{n-1}}<+\infty,$$
or $$\ell=0\quad \text{ and }\quad
0<R^{(0)}=\lim_{j\to+\infty}\frac{r_j^{(n-p)}}{\d_j\e_j^{n-1}}<+\infty.
$$
In the following proposition we suitably modify $(w_j)$ in each
surrounding cylinder in order to get an admissible test function
for the minimum problem (\ref{phigj}) or (\ref{phigj12}).

\begin{proposition}\label{gliminfclose}
Let $\ell \in [0,+\infty]$. Then
$$
\liminf_{j \to +\infty}\frac{1}{\d_j} \left( \int_{E_j^{+\d_j}}
W(D w_j)\, dx + \int_{E_j^{-\d_j}} W(D w_j)\, dx \right) \geq
R^{(\ell)} \int_\o \varphi^{(\ell)} (u^+-u^-) \, dx_\a + o(1)\,,
$$
as $\g\to 0^+$.
\end{proposition}

\noindent {\it Proof. }Let $\ell \in (0,+\infty]$, the case
$\ell=0$ can be treated similarly. Let $i \in Z_j$ and
$N_j=\frac{\e_j}{r_j}$. Since $\rho_j^i <\g \e_j$, we can define
$$
\zeta_j^i:= \left\{
\begin{array}{ll}
w_j(x_i^{\e_j} + r_j\, y_\a, \d_j\, y_n) - u_j^{i-} & \text{in }
\big(B_{\rho_j^i/r_j}^{n-1} \times I \big) \setminus
C_{1,\rho_j^i/r_j}\\
&\\
(u_j^{i+} - u_j^{i-}) & \text{in } \big(B_{\g N_j}^{n-1}
\setminus B_{\rho_j^i/r_j}^{n-1}\big)^+ \\
&\\
0 & \text{in } \big(B_{\g N_j}^{n-1} \setminus
B_{\rho_j^i/r_j}^{n-1} \big)^- \,,
\end{array}
\right.
$$
where $N_j={\e_j/r_j}$. Then $\zeta_j^i \in W^{1,p}((B_{\g
N_j}^{n-1} \times I) \setminus C_{1,\g N_j};\Rb^m)$, $\zeta_j^i =
(u_j^{i+} - u_j^{i-})$ on $\big(\partial B_{\g N_j}^{n-1}\big)^+ $
and $\zeta_j^i=0$ on $\big(\partial B_{\g N_j}^{n-1}\big)^-$.
Since $W(0)=0$, changing variable, by (\ref{phigj}) we get
\begin{eqnarray}\label{1136}
&&\frac{1}{\d_j} \left(
\int_{B_{\rho_j^i}^{n-1}(x_i^{\e_j})^{+\d_j}} W(D w_j)\, dx +
\int_{B_{\rho_j^i}^{n-1}(x_i^{\e_j})^{-\d_j}} W(D
w_j)\, dx \right)\nonumber\\
&= & r_j^{n-1}\left(\int_{\big( B_{\rho_j^i/r_j}^{n-1} \big)^+}
W\Bigl(r_j^{-1} D_\a \zeta_j^i |\d_j^{-1} D_n \zeta_j^i \Bigr)\,
dy + \int_{\big( B_{\rho_j^i/r_j}^{n-1}\big)^-} W \Bigl(r_j^{-1}
D_\a \zeta_j^i |\d_j^{-1} D_n
\zeta_j^i \Bigr)\, dy \right)\nonumber\\
&=& r_j^{n-1}\int_{( B_{\g N_j}^{n-1} \times I) \setminus C_{1,\g
N_j}} W\Bigl(r_j^{-1} D_\a \zeta_j^i |\d_j^{-1} D_n \zeta_j^i
\Bigr)\dy
\nonumber\\
&\geq& r_j^{n-1-p}\varphi^{(\ell)}_{\g,j}(u_j^{i+} - u_j^{i-})\,.
\end{eqnarray}
Summing up in (\ref{1136}), for $i \in Z_j$, we get  that
\begin{eqnarray}\label{1137}
&&\nonumber \frac{1}{\d_j} \left( \int_{E_j^{+\d_j}} W(D w_j)\, dx
+ \int_{E_j^{-\d_j}} W(D
w_j)\, dx \right)\\
&=&\nonumber \sum_{i \in Z_j}\frac{1}{\d_j} \left(
\int_{B_{\rho_j^i}^{n-1}(x_i^{\e_j})^{+\d_j}} W(D w_j)\, dx +
\int_{B_{\rho_j^i}^{n-1}(x_i^{\e_j})^{-\d_j}} W(D w_j)\, dx
\right)\\
&\geq& r_j^{n-1-p}\,\sum_{i \in
Z_j}\varphi^{(\ell)}_{\g,j}(u_j^{i+} - u_j^{i-})=
\frac{r_j^{n-1-p}}{\e_j^{n-1}} \sum_{i \in Z_j}\e_j^{n-1}
\varphi^{(\ell)}_{\g,j}(u_j^{i+} - u_j^{i-})\,.
\end{eqnarray}
Passing to the limit as $j \to +\infty$ we get, by (\ref{rinfty})
and Proposition \ref{rs}, that
\begin{eqnarray*}
&&\liminf_{j \to +\infty}\frac{1}{\d_j} \left( \int_{E_j^{+\d_j}}
W(D w_j)\, dx + \int_{E_j^{-\d_j}} W(D
w_j)\, dx \right)\\
&\geq &  R^{(\ell)} \int_\o \varphi^{(\ell)}(u^+ - u^-)\, dx_\a
\\
&& + R^{(\ell)} \, \liminf_{j \to +\infty}\int_\o \Bigl(\sum_{i
\in Z_j}\varphi^{(\ell)}_{\g,j}(u_j^{i+}-
u_j^{i-})\chi_{Q_{i,\e_j}^{n-1}} -
\varphi^{(\ell)}(u^+-u^-)\Bigr)\,dx_\a\\
&=& R^{(\ell)} \int_\o \varphi^{(\ell)}(u^+ - u^-)\, dx_\a +
o(1)\,,
\end{eqnarray*}
as $\g\to 0^+$, which completes the proof.
\hfill$\Box$\\

We now prove the liminf inequality for any arbitrary converging
sequence.
\begin{lemma}\label{gammaliminf}
Let $\ell \in [0,+\infty]$. For every sequence $(u_j)$ converging
to $(u^+,u^-)$ we have
\begin{eqnarray*}
\liminf_{j\to +\infty} \mathcal F_j(u_j) &\geq& \int_\o \Q_{n-1}
\overline W(D_\a u^+)\, dx_\a + \int_\o \Q_{n-1} \overline W(D_\a
u^-)\,
dx_\a \\
&&+ R^{(\ell)} \int_\o \varphi^{(\ell)}(u^+ - u^-)\, dx_\a\,.
\end{eqnarray*}
\end{lemma}

\noindent {\it Proof. } Let $(u_j)\to (u^+,u^-)$ be such that
$\liminf_{j\to +\infty} \mathcal F_j(u_j)<+\infty$.  Reasoning as
in \cite{Ans} Proposition 5.2, by \cite{BDV} Lemma 3.5, upon
passing to a subsequence, for every $M>0$ and $\eta>0$, we have
the existence of $R_M>M$ and of a Lipschitz function $\Phi_M \in
\C^1_c(\Rb^m;\Rb^m)$ with ${\rm Lip}(\Phi_M)=1$ such that
$$\Phi_M(z)=\left\{
\begin{array}{lll}
z & \text{if} & |z|< R_M,\\
&&\\
0 & \text{if} & |z| >2 R_M
\end{array}
\right.
$$ and
\begin{equation}\label{janv}
\liminf_{j \to +\infty} \mathcal F_j(u_j) \geq \liminf_{j \to
+\infty} \mathcal F_j(\Phi_M (u_j)) -\eta\,.
\end{equation}
Note that $(\Phi_M (u_j))\subset W^{1,p}(\O_j;\Rb^m)\cap
L^\infty(\O_j;\Rb^m)$, $\sup_{j \in \Nb}\| \Phi_M(u_j)
\|_{L^\infty(\O_j;\Rb^m)}<R_M$ and it converges to
$(\Phi_M(u^+),\Phi_M(u^-))$ as $j \to +\infty$. Hence, if we apply
(\ref{lowerbd}), Propositions \ref{gliminffar} and
\ref{gliminfclose} to $(\Phi_M (u_j))$ in place of $(u_j)$,
letting $\g \to 0$ and $k \to +\infty$, we get that
\begin{eqnarray}\label{fev}
\liminf_{j \to +\infty} \mathcal F_j(\Phi_M(u_j)) & \geq & \int_\o
\Q_{n-1} \overline W(D_\a \Phi_M(u^+))\, dx_\a + \int_\o \Q_{n-1}
\overline
W(D_\a \Phi_M(u^-))\, dx_\a\nonumber\\
&& + R^{(\ell)} \int_\o \varphi^{(\ell)}(\Phi_M(u^+) -
\Phi_M(u^-))\, dx_\a.
\end{eqnarray}
Moreover $\Phi_M(u^\pm) \rightharpoonup u^\pm$ weakly in
$W^{1,p}(\o;\Rb^m)$ as $M\to +\infty$; hence, by (\ref{janv}),
(\ref{fev}), the lower semicontinuity of $\int_\o \Q_{n-1}
\overline W(D_\a u)\dxa$ with respect to the weak
$W^{1,p}(\o;\Rb^m)$-convergence, and (\ref{AA}) we have that
\begin{eqnarray}\label{mars}
&&\liminf_{j \to +\infty} \mathcal F_j(u_j)\nonumber\\
&\geq& \int_\o \Q_{n-1} \overline W(D_\a u^+)\, dx_\a + \int_\o
\Q_{n-1} \overline W(D_\a u^-)\, dx_\a+ R^{(\ell)} \int_\o
\varphi^{(\ell)}(u^+ - u^-)\, dx_\a-\eta\,,
\end{eqnarray}
and by the arbitrariness of $\eta$, the thesis. \hfill$\Box$


\subsection{The limsup inequality}
\noindent For every $ (u^+,u^-)\in W^{1,p}(\o,\rr^m)\times
W^{1,p}(\o,\rr^m)$ the limsup inequality is obtained by suitably
modifying the recovery sequences $(u_j^{\pm})$ for the
$\Gamma$-limits of $$ \frac{1}{\d_j}\int_{\o^{+\d_j}} W( Du )\, dx
\quad \text{ and } \quad \frac{1}{\d_j}\int_{\o^{-\d_j}} W( Du )\,
dx.$$

\begin{lemma}
Let $\ell \in [0,+\infty]$ and let $\o$ be an open bounded subset
of $\rr^{n-1}$ such that $\HH^{n-1}(\partial \o)=0$. Then, for all
$(u^+,u^-)\in W^{1,p}(\o,\rr^m)\times W^{1,p}(\o,\rr^m)$ and for
all $\eta>0$ there exists a sequence $(\bar u_j)\subset
W^{1,p}(\O_j;\Rb^m)$ converging to $(u^+,u^-)$ such that
\begin{eqnarray*}
\limsup_{j\to +\infty} \mathcal F_j(\bar u_j) &\le& \int_\o
\Q_{n-1} \overline W(D_\a u^+)\, dx_\a + \int_\o \Q_{n-1}
\overline W(D_\a
u^-)\, dx_\a \\
&&+ R^{(\ell)} \int_\o \varphi^{(\ell)}(u^+ - u^-)\, dx_\a + \eta
R^{(\ell)} \HH^{n-1}(\o)\,.
\end{eqnarray*}
\end{lemma}

\noindent{\it Proof.} The proof of the limsup is divided into
three steps. We first construct a sequence $(\bar u_j) \subset
W^{1,p}(\O_j;\Rb^m)$ that we expect to be a recovery sequence. In
the second step we prove that $(\bar u_j)$ converges to
$(u^+,u^-)$. Finally, we prove that it satisfies the limsup
inequality. We first deal with the case $\ell \in (0,+\infty]$.\\

{\bf Step 1: Definition of a recovery sequence.} Let $u^\pm \in
W^{1,p}(\o;\Rb^m) \cap L^\infty(\o;\Rb^m)$. According to
\cite{LDR} Theorem 2 and \cite{Bo&Fo} Theorem 1.1, there exist two
sequences $(u^{\pm}_j)\subset W^{1,p}(\o^{\pm\d_j};\Rb^m)$ such
that $u_j^{\pm} \to u^{\pm}$, the sequences of gradients $(|D
u_j^{\pm}|^p/\d_j)$ are equi-integrable on $\o^{\pm\d_j}$,
respectively, and
\begin{equation}\label{canne}
\lim_{j \to +\infty}{1\over \d_j}\int_{\o^{\pm\d_j}} W( Du^\pm_j)
\, dx = \int_\o \Q_{n-1} \overline W(D_\a u^\pm)\, dx_\a\,.
\end{equation}
Moreover, using a truncation argument (as in \cite{AnsB} Lemma
6.1, Step 2) we may assume without loss of generality that
$$
\sup_{j \in \Nb}\|u^{\pm}_j\|_{L^\infty(\o^{\pm\d_j};\Rb^m)} <
+\infty\,.
$$
Let $u_j:=u_j^+ \chi_{\o^{+\d_j}} + u_j^- \chi_{\o^{-\d_j}}\in
W^{1,p}(\o^{+\d_j} \cup \o^{-\d_j};\Rb^m)$ and let $(w_j)$ be the
sequence obtained from $(u_j)$ as in Lemma
\ref{important+equiint}, then $\sup_{j \in
\Nb}\|w_j\|_{L^\infty(\o^{\pm\d_j};\Rb^m)}< +\infty$.

\bigskip

We first define $(\bar u_j)$ `far' from the connecting zones\ie
\begin{equation}\label{defubar}
\bar u_j:=w_j\; \text{  in }\; \Bigl( \o \setminus \bigcup_{i \in
\Zb^{n-1}}B_{\rho_j}^{n-1}(x_i^{\e_j})\Bigr)^{\pm\d_j}\,.
\end{equation}
Then we pass to define $(\bar u_j)$ on each
$B_{\rho_j}^{n-1}(x_i^{\e_j})^{\pm\d_j}$ making a distinction
between the indices $i\in Z_j$ and $i\in \Zb^{n-1}\setminus Z_j$.

If $i \in Z_j$, by (\ref{phigj}), for every $\eta>0$ there exists
$\zeta_{\g,j}^i \in X^\g_j(u_j^{i+}- u_j^{i-})$ such that
\begin{equation}\label{presque}
\int_{(B_{\g N_j}^{n-1} \times I) \setminus C_{1,\g N_j}}r_j^p \,
W\left(r_j^{-1} D_\a \zeta_{\g,j}^i | \d_j^{-1} D_n \zeta_{\g,j}^i
\right) dx \leq \varphi^{(\ell)}_{\g,j}(u_j^{i+}- u_j^{i-}) +
\eta.
\end{equation}
Then, we define
\begin{equation}\label{defubarZj}
\bar u_j:=
\zeta_{\g,j}^i\left(\frac{x_\a-x_i^{\e_j}}{r_j},\frac{x_n}{\d_j}\right)+u_j^{i-}\;
\text{ in }\;  B^{n-1}_{\rho_j}(x_i^{\e_j})^{\pm\d_j}\,, \quad
i\in Z_j\,.
\end{equation}
In particular, $\bar u_j=u_j^{i\pm}=w_j$ on $\big(\partial
B^{n-1}_{\rho_j}(x_i^{\e_j})\big)^{\pm\d_j}$.

Let us now deal with the contact zones not well contained in
$\o$\ie with the indices $i \not\in Z_j$. For fixed $\g>0$ and $j$
large enough we have that $\g N_j>2$. Let $\psi \in
W^{1,p}(B^{n-1}_2;[0,1])$ be such that $\psi =1$ on $\partial
B_2^{n-1}$ and $\psi=0$ in $B_1^{n-1}$ and define
$$
\psi_{\g,j}(x):=\left\{
\begin{array}{lll}
0 & \text{ in } & (B_{\g N_j}^{n-1})^- \\
\psi(x_\a) & \text{ in } & (B_2^{n-1})^+ \\
1 & \text{ in } & (B_{\g N_j}^{n-1} \setminus B_2^{n-1})^+ \,.
\end{array}
\right.
$$
Then $\psi_{\g,j} \in W^{1,p}((B^{n-1}_{\g N_j} \times I)
\setminus C_{1,\g N_j};[0,1])$, $\psi_{\g,j} = 1$ on
$\big(\partial B_{\g N_j}^{n-1}\big)^+$ and $\psi_{\g,j} = 0$ on
$\big( \partial B_{\g N_j}^{n-1}\big)^-$. Let $w_j^{\pm}=w_j\,
\chi_{\o^{\pm\d_j}}$, we extend both of them to the whole
$\o\times (-\d_j,\d_j)$ by reflection\ie we define
$\tilde{w}_j^{\pm}(x_\a, x_n)= w_j^{\pm}(x_\a, -x_n)$ for $x\in
\o^{\mp\d_j}$ and $\tilde{w}_j^{\pm}(x)= w_j^{\pm}(x)$ for $x\in
\o^{\pm\d_j}$. Hence, we define
\begin{equation}\label{defubarZn}
\bar u_j:=\psi_{\g,j}
\left(\frac{x_\a-x_i^{\e_j}}{r_j},\frac{x_n}{\d_j}
\right)\tilde{w}_j^+ + \Bigg( 1- \psi_{\g,j}
\left(\frac{x_\a-x_i^{\e_j}}{r_j},\frac{x_n}{\d_j} \right)
\Bigg)\tilde{w}_j^-
\end{equation}
in $\big( B^{n-1}_{\rho_j}(x_i^{\e_j}) \times (-\d_j,\d_j)\big)
\cap \O_j$ and for $i\in \Zb^{n-1}\setminus Z_j$. In particular,
we have that $\bar u_j=w_j$ on $\big( \partial
B^{n-1}_{\rho_j}(x_i^{\e_j}) \times (-\d_j,\d_j)\big) \cap \O_j$;
thus $(\bar u_j) \subset W^{1,p}(\O_j;\Rb^m)$.

\bigskip

{\bf Step 2: The sequence $(\bar u_j)$ weakly converges to
$(u^+,u^-)$}. Let us check (\ref{strongconv}) and
(\ref{weakconv}). We will only treat the upper cylinder $\o^{+
\d_j}$, the lower part being analogous. First
\begin{eqnarray}\label{canne0}
&&\frac{1}{\d_j}\int_{\o^{+\d_j}} |\bar u_j - u^+|^p \, dx
\nonumber\\
&=& \frac{1}{\d_j}\int_{ \left( \o \setminus \bigcup_{i \in
\Zb^{n-1}} B_{\rho_j}^{n-1}(x_i^{\e_j})\right)^{+\d_j}} |w^+_j -
u^+|^p \,
dx\nonumber\\
&&+ \frac{1}{\d_j}\sum_{i \in Z_j} \int_{
B_{\rho_j}^{n-1}(x_i^{\e_j})^{+\d_j}} \left|\zeta_{\g,j}^i
\left(\frac{x_\a-x_i^{\e_j}}{r_j},\frac{x_n}{\d_j}\right)+u_j^{i-}
-u^+\right|^p \, dx\nonumber\\
&&+ \frac{1}{\d_j}\sum_{i \in \Zb^{n-1}\setminus Z_j}
\int_{\big(\o \cap B_{\rho_j}^{n-1}(x_i^{\e_j})\big)^{+\d_j}}
\left|\psi_{\g,j}
\left(\frac{x_\a-x_i^{\e_j}}{r_j},\frac{x_n}{\d_j}\right)(w_j^+
- \tilde{w}^-_j) + \tilde{w}^-_j - u^+ \right|^p \, dx\nonumber\\
&\leq & \frac{1}{\d_j}\int_{\o^{+\d_j}} |w_j - u^+|^p \, dx +
c\sum_{i \in Z_j} \int_{
B_{\rho_j}^{n-1}(x_i^{\e_j})} |u^+ - u_j^{i+}|^p \, dx_\a\nonumber\\
&&+ \frac{c}{\d_j}\sum_{i \in Z_j} \int_{
B_{\rho_j}^{n-1}(x_i^{\e_j})^{+\d_j}} \left|\zeta_{\g,j}^i
\left(\frac{x_\a-x_i^{\e_j}}{r_j},\frac{x_n}{\d_j}\right)-
(u_j^{i+} -
u_j^{i-}) \right|^p dx\nonumber\\
&&+{c\over \d_j}\, \int_{\left(\o \cap \bigcup_{i\in \Zb^{n-1}
\setminus Z_j}B_{\rho_j}^{n-1}(x_i^{\e_j})\right)^{+\d_j}}
\left(|w^+_j|^p + |\tilde{w}_j^-|^p + |u^+|^p\right)\dx\,.
\end{eqnarray}
Since $\lim_{j\to +\infty}\HH^{n-1}\Bigl(\o \cap \bigcup_{i\in
\Zb^{n-1}\setminus Z_j}B_{\rho_j}^{n-1}(x_i^{\e_j})\Bigr)=0$ and
$\sup_{j \in \Nb}\|w^{\pm}_j\|_{L^\infty(\o^{\pm\d_j};\Rb^m)}<
+\infty$, we have that
\begin{equation}\label{canne3}
\lim_{j\to +\infty}{c\over \d_j}\, \int_{\left(\o \cap
\bigcup_{i\in \Zb^{n-1} \setminus
Z_j}B_{\rho_j}^{n-1}(x_i^{\e_j})\right)^{+\d_j}} \left(|w^+_j|^p +
|\tilde{w}_j^-|^p + |u^+|^p\right)\dx=0\,.
\end{equation}
Moreover, reasoning as in the proof of Proposition \ref{rs} (see
inequality (\ref{corep})), we have that
\begin{equation}\label{canne4}
\lim_{j \to +\infty} \sum_{i \in Z_j} \int_{
B_{\rho_j}^{n-1}(x_i^{\e_j})} |u^+ - u_j^{i+}|^p \, dx_\a=0\,,
\end{equation}
and, by the convergence $w_j\to (u^+,u^-)$, it remains only to
prove that
\begin{equation}\label{canne2}
\lim_{j \to +\infty} \frac{1}{\d_j}\sum_{i \in Z_j} \int_{
B_{\rho_j}^{n-1}(x_i^{\e_j})^{+\d_j}} \left|\zeta_{\g,j}^i
\left(\frac{x_\a-x_i^{\e_j}}{r_j},\frac{x_n}{\d_j}\right)-
(u_j^{i+} - u_j^{i-}) \right|^p dx=0\,.
\end{equation}
In fact, changing variable, we get that
\begin{eqnarray*}
&&\frac{1}{\d_j}\sum_{i \in Z_j} \int_{
B_{\rho_j}^{n-1}(x_i^{\e_j})^{+\d_j}} \left|\zeta_{\g,j}^i
\left(\frac{x_\a-x_i^{\e_j}}{r_j},\frac{x_n}{\d_j}\right)-
(u_j^{i+}- u_j^{i-}) \right|^p dx\\
&=& r_j^{n-1} \sum_{i \in Z_j} \int_{(B_{\g N_j} ^{n-1})^+}
\left|\zeta_{\g,j}^i(x)- (u_j^{i+} - u_j^{i-}) \right|^p \dx\,,
\end{eqnarray*}
and by, Poincar\'e's Inequality
$$
\int_{B_{\g N_j}^{n-1}} \left|\zeta_{\g,j}^i (x_\a,x_n)- (u_j^{i+}
- u_j^{i-}) \right|^p dx_\a \leq c\, (\g N_j)^p \int_{B_{\g
N_j}^{n-1}}  |D_\a \zeta_{\g,j}^i (x_\a,x_n)|^p \dx_\a
$$
for a.e. $x_n \in (0,1)$. Hence, by  the $p$-growth condition
(\ref{pgrowth}) and (\ref{presque}) if we integrate with respect
to $x_n$ and sum up in $i \in Z_j$, we get that
\begin{eqnarray}\label{canne1}
&&\frac{1}{\d_j}\sum_{i \in Z_j} \int_{
B_{\rho_j}^{n-1}(x_i^{\e_j})^{+\d_j}} \left|\zeta_{\g,j}^i
\left(\frac{x_\a-x_i^{\e_j}}{r_j},\frac{x_n}{\d_j}\right)-
(u_j^{i+} - u_j^{i-}) \right|^p dx\nonumber\\
&\leq & c\, r_j^{n-1} \g^p N_j^p \sum_{i \in Z_j} \int_{(B_{\g
N_j}^{n-1})^+}  |D_\a \zeta_{\g,j}^i |^p \, dx\nonumber\\
&\leq & c\, r_j^{n-1} \g^p N_j^p \sum_{i \in Z_j} \int_{(B_{\g
N_j}^{n-1})^+}  \left|\left( D_\a \zeta_{\g,j}^i \Big|
\frac{r_j}{\d_j}
D_n \zeta_{\g,j}^i \right) \right|^p \, dx\nonumber\\
&\leq & c\, r_j^{n-1} \g^p N_j^p \sum_{i \in Z_j} \left(
\varphi^{(\ell)}_{\g,j}(u_j^{i+} - u_j^{i-}) +\eta + r_j^p\,
\HH^{n-1}( B_{\g N_j}^{n-1}) \right)\nonumber\\
&\leq & c\,\g^p\, \e_j^p\,\frac{ r_j^{n-1-p}}{\e_j^{n-1}} \left(
\sum_{i \in Z_j} \e_j^{n-1} \varphi^{(\ell)}_{\g,j}(u_j^{i+} -
u_j^{i-}) +\left(\eta + c\,\g^{n-1}\,
\frac{\e_j^{n-1}}{r_j^{n-1-p}}\right) \HH^{n-1}(\o)\right).
\end{eqnarray}
By Proposition \ref{rs} and (\ref{rinfty}), passing to the limit
as $j\to +\infty$ in (\ref{canne1}), we get (\ref{canne2}).

It remains to prove that (\ref{weakconv}) holds. In fact,
\begin{eqnarray}\label{gradubar}
\nonumber
&& \frac{1}{\d_j}\int_{\o^{+ \d_j}}|D {\bar u}_j|^p\,dx \\
&=& \nonumber\frac{1}{\d_j}\int_{\left( \o \setminus \bigcup_{i
\in \Zb^{n-1}} B_{\rho_j}^{n-1}(x_i^{\e_j})\right)^{+ \d_j}}|D
w^\pm_j|^p\,dx\\
&&\nonumber +  \frac{1}{\d_j}\int_{\bigcup_{i\in Z_j}
B_{\rho_j}^{n-1}(x_i^{\e_j})^{+ \d_j}}
\left|\left(r_j^{-1}D_\a\zeta_{\g,j}^i
\Bigl(\frac{x_\a-x_i^{\e_j}}{r_j},\frac{x_n}{\d_j}\Bigr)\Big|
\d_j^{-1} D_n\zeta_{\g,j}^i
\Bigl(\frac{x_\a-x_i^{\e_j}}{r_j},\frac{x_n}{\d_j}
\Bigr)\right)\right|^p\dx\\
&&+ \frac{1}{\d_j}\int_{\left(\bigcup_{i\in\Zb^{n-1}\setminus
Z_j}B_{\rho_j}^{n-1}(x_i^{\e_j})\cap \o\right)^{ + \d_j}}|D {\bar
u}_j|^p\,dx \,.
\end{eqnarray}
It can be easily shown that
\begin{eqnarray}\label{gradubarZj}
&&\nonumber\frac{1}{\d_j}\int_{\bigcup_{i\in Z_j}
B_{\rho_j}^{n-1}(x_i^{\e_j})^{+ \d_j}}
\left|\left(r_j^{-1}D_\a\zeta_{\g,j}^i
\Bigl(\frac{x_\a-x_i^{\e_j}}{r_j},\frac{x_n}{\d_j}\Bigr)\Big|
\d_j^{-1} D_n\zeta_{\g,j}^i
\Bigl(\frac{x_\a-x_i^{\e_j}}{r_j},\frac{x_n}{\d_j}
\Bigr)\right)\right|^p\dx\\
&\le& { r_j^{n-1-p}\over \e_j^{n-1}} \Bigl(\sum_{i\in Z_j}
\e_j^{n-1}\,\varphi^{(\ell)}_{\g,j}(u_j^{i+} - u_j^{i-}) \Bigr)+
\HH^{n-1}(\o) \Bigl(\eta  { r_j^{n-1-p}\over \e_j^{n-1}}  +
\g^{n-1}\Bigr)\,;
\end{eqnarray}
while,
\begin{eqnarray}\label{gradubarZn}
&&\nonumber\frac{1}{\d_j}\int_{\left(
\bigcup_{i\in\Zb^{n-1}\setminus
Z_j}B_{\rho_j}^{n-1}(x_i^{\e_j})\cap \o \right)^{+ \d_j}}|D {\bar
u}_j|^p\,dx \\
& \le & \nonumber c\, \sum_{i\in \Zb^{n-1}\setminus Z_j}\left(
\frac{1}{r_j^p \d_j} \int_{(B_{\rho_j}^{n-1}(x_i^{\e_j})\cap
\o)^{+\d_j}} \Bigl|D_\a \psi_{\g,j}
\Bigl(\frac{x_\a-x_i^{\e_j}}{r_j},\frac{x_n}{\d_j} \Bigr)\Bigr|^p
\, \big(|w^+_j|^p + |\tilde{w}^-_j|^p\big)\ dx\right.
\\
&&\left.\nonumber\qquad\qquad + {1\over \d_j}
\int_{(B_{\rho_j}^{n-1}(x_i^{\e_j})\cap \o)^{+\d_j}} (|Dw^+_j|^p
+|D \tilde{w}^-_j|^p)\, dx \right)\\
&\le&\nonumber c\,\sum_{i\in \Zb^{n-1}\setminus Z_j}\left(
r_j^{n-1-p}\int_{B_2^{n-1}}|D_\a \psi|^p\, dx_\a +{1\over \d_j}
\int_{(B_{\rho_j}^{n-1}(x_i^{\e_j})\cap \o)^{+ \d_j}}
|Dw^+_j|^p\,dx\right.\\
&& \nonumber \left. \qquad \qquad +{1\over \d_j}
\int_{(B_{\rho_j}^{n-1}(x_i^{\e_j})\cap
\o)^{- \d_j}} |Dw^-_j|^p\,dx\right)\\
 \nonumber &\le& c\,\sum_{i\in \Zb^{n-1}\setminus Z_j}\left(
{r_j^{n-1-p}\over \e_j^{n-1}} \HH^{n-1}(Q_{i,\e_j}^{n-1}) +
{1\over \d_j} \int_{(B_{\rho_j}^{n-1}(x_i^{\e_j})\cap
\o)^{\pm\d_j}}
|Dw^\pm_j|^p\,dx\right.\\
&& \left. \qquad \qquad +{1\over \d_j}
\int_{(B_{\rho_j}^{n-1}(x_i^{\e_j})\cap \o)^{- \d_j}}
|Dw^-_j|^p\,dx\right)\,.
\end{eqnarray}
Note that the previous sum can be computed over all $i \in
\Zb^{n-1} \setminus Z_j$ such that $Q^{n-1}_{i,\e_j} \cap \o \neq
\emptyset$. Let
$$\o'_j:= \bigcup_{i \in \Zb^{n-1}
\setminus Z_j ,\, Q^{n-1}_{i,\e_j} \cap \o \neq \emptyset}
Q^{n-1}_{i,\e_j},$$ then
\begin{equation}\label{babadjian}
\sum_{i \in \Zb^{n-1} \setminus Z_j ,\, Q^{n-1}_{i,\e_j} \cap \o
\neq \emptyset} \mathcal H^{n-1}(Q^{n-1}_{i,\e_j}) = \mathcal
H^{n-1}(\o'_j) \to \mathcal H^{n-1}(\partial \o)=0.
\end{equation}
Moreover, by Lemma \ref{important+equiint} we have that
$\sup_{j}{1\over \d_j} \int_{ \o^{\pm\d_j}}
|Dw^\pm_j|^p\,dx<+\infty$; hence, by Proposition \ref{rs},
(\ref{rinfty}), (\ref{gradubar}), (\ref{gradubarZj}) and
(\ref{gradubarZn}) we get (\ref{weakconv}).

\bigskip

{\bf Step 3: The sequence $(\bar u_j)$ is a recovery sequence.} We
now prove the limsup inequality.
\begin{eqnarray}\label{step1}
&&\limsup_{j \to +\infty} \int_{\o^{\pm \d_j}} W(D \bar u_j)\, dx
\nonumber\\
&=&  \limsup_{j \to +\infty}  \frac{1}{\d_j} \Bigg( \int_{\left(\o
\setminus \bigcup_{i \in
\Zb^{n-1}}B_{\rho_j}^{n-1}(x_i^{\e_j})\right)^{\pm\d_j}} W( D \bar
u_j)\, dx + \int_{\bigcup_{i \in
Z_j}B_{\rho_j}^{n-1}(x_i^{\e_j})^{\pm\d_j}} W( D \bar u_j)\, dx
\nonumber\\
&&\qquad\qquad\qquad + \int_{\left(\o \cap \bigcup_{i \in
\Zb^{n-1} \setminus
Z_j}B_{\rho_j}^{n-1}(x_i^{\e_j})\right)^{\pm\d_j}} W( D \bar
u_j)\, dx\Bigg)\,.
\end{eqnarray}
We deal with the first term in (\ref{step1}). The definition of
$\bar u_j$ (\ref{defubar}), Lemma \ref{important+equiint} and
(\ref{canne}), yield
\begin{eqnarray}\label{step2}
&&\nonumber\limsup_{j \to +\infty} \frac{1}{\d_j} \int_{\left(\o
\setminus \bigcup_{i \in
\Zb^{n-1}}B_{\rho_j}^{n-1}(x_i^{\e_j})\right)^{\pm\d_j}} W( D \bar
u_j)\, dx\nonumber\\
&=& \nonumber\limsup_{j \to +\infty} \frac{1}{\d_j}\int_{\left(\o
\setminus \bigcup_{i \in
\Zb^{n-1}}B_{\rho_j}^{n-1}(x_i^{\e_j})\right)^{\pm\d_j}} W( D
w_j)\dx\nonumber\\
&\leq &\nonumber \limsup_{j \to +\infty}  \frac{1}{\d_j}
\int_{\o^{\pm\d_j}} W( D u^{\pm}_j)\dx + o(1)\nonumber\\
&=& \int_\o \Q_{n-1} \overline W(D_\a u^\pm)\, dx_\a + o(1)\,,
\end{eqnarray}
as $\g\to 0^+$. For every $i \in Z_j$, by (\ref{defubarZj}) and
(\ref{presque}) we get that
\begin{eqnarray*}
&&\frac{1}{\d_j} \Bigg(
\int_{B^{n-1}_{\rho_j}(x_i^{\e_j})^{+\d_j}}W(D \bar u_j)\, dx +
\int_{B^{n-1}_{\rho_j}(x_i^{\e_j})^{-\d_j}}W(D \bar
u_j)\, dx \Bigg)\\
&=&  r_j^{n-1}\int_{(B_{\g N_j}^{n-1} \times I) \setminus C_{1,\g
N_j}} W\left(r_j^{-1} D_\a \zeta_{\g,j}^i | \d_j^{-1} D_n
\zeta_{\g,j}^i \right)
dx\\
&\leq &  r_j^{n-1-p} \left( \varphi^{(\ell)}_{\g,j}(u_j^{i+}-
u_j^{i-}) + \eta \right)\,;
\end{eqnarray*}
hence, by (\ref{rinfty}) and Proposition \ref{rs} we get
\begin{eqnarray}\label{step3}
&&\nonumber\limsup_{j \to +\infty} \frac{1}{\d_j}\Bigg(
\int_{\bigcup_{i \in Z_j} B^{n-1}_{\rho_j}(x_i^{\e_j})^{+\d_j}}W(D
\bar u_j)\, dx + \int_{\bigcup_{i \in Z_j}
B^{n-1}_{\rho_j}(x_i^{\e_j})^{-\d_j}}W(D \bar u_j)\, dx
\Bigg)\\
&\leq &\nonumber R^{(\ell)} \int_\o \varphi^{(\ell)} (u^+ - u^-)\,
dx_\a
+ R^{(\ell)}\, \HH^{n-1}(\o) \,\eta \\
&&\nonumber + \limsup_{j \to +\infty}\int_\o \Bigl|\sum_{i \in
Z_j}\varphi^{(\ell)}_{\g,j}(u_j^{i+}-
u_j^{i-})\chi_{Q_{i,\e_j}^{n-1}} -
\varphi^{(\ell)}(u^+-u^-)\Bigr|\,dx_\a\\
&=& R^{(\ell)} \int_\o \varphi^{(\ell)} (u^+ - u^-)\, dx_\a +
R^{(\ell)}\, \HH^{n-1}(\o) \,\eta + o(1)\,,
\end{eqnarray}
as $\g\to 0^+$. Finally, for $i \not\in Z_j$, by the $p$-growth
condition (\ref{pgrowth}) and (\ref{gradubarZn}), we obtain
\begin{eqnarray*}
&&\frac{1}{\d_j} \Bigg(\int_{\left(\bigcup_{i\in
\Zb^{n-1}\setminus Z_j}B_{\rho_j}^{n-1}(x_i^{\e_j})\cap \o
\right)^{\pm\d_j}} W(D \bar
u_j)\, dx \Bigg)\\
&\leq &\sum_{i\in \Zb^{n-1}\setminus Z_j}  \frac{\beta}{\d_j}
\Bigg( \int_{(B_{\rho_j}^{n-1}(x_i^{\e_j})\cap \o)^{\pm\d_j}} (1
+|D \bar u_j|^p)\, dx \Bigg)\\
&\le&c\, \HH^{n-1}\Bigl( \bigcup_{i\in \Zb^{n-1}\setminus Z_j}
B_{\rho_j}^{n-1}(x_i^{\e_j})\cap \o\Bigr)\\
&&+ c\,\sum_{i\in \Zb^{n-1}\setminus Z_j}\left({ r_j^{n-1-p}\over
\e_j^{n-1}} \HH^{n-1}(Q_{i,\e_j}^{n-1})  + {1\over \d_j}
\int_{(B_{\rho_j}^{n-1}(x_i^{\e_j})\cap \o)^{+ \d_j}}
|Dw^+_j|^p\,dx \right.\\
&& \left. \hspace{6cm} + {1\over \d_j}
\int_{(B_{\rho_j}^{n-1}(x_i^{\e_j})\cap \o)^{- \d_j}}
|Dw^-_j|^p\,dx\right)\,.
\end{eqnarray*}
Since
$$
\lim_{j\to +\infty} \HH^{n-1}\Bigl(\bigcup_{i\in
\Zb^{n-1}\setminus Z_j} B_{\rho_j}^{n-1}(x_i^{\e_j})\cap
\o\Bigr)=0\,,
$$
by (\ref{rinfty}), the equi-integrability of $(|D
w_j^{\pm}|^p/\d_j)$ on $\o^{\pm\d_j}$ and (\ref{babadjian}), we
deduce
\begin{eqnarray}\label{step4}
\limsup_{j \to +\infty}\frac{1}{\d_j}\int_{\left(\o \cap
\bigcup_{i \in \Zb^{n-1} \setminus
Z_j}B_{\rho_j}^{n-1}(x_i^{\e_j})\right)^{\pm\d_j}} W(D \bar u_j)\,
dx =0\,.
\end{eqnarray}
Gathering (\ref{step1})-(\ref{step4}) and passing to the limit as
$\g\to 0^+$ we get the limsup inequality
for every $u^\pm \in W^{1,p}(\o;\Rb^m) \cap L^\infty(\o;\Rb^m)$.\\

We remove the boundedness assumption simply noting that any
arbitrary $W^{1,p}(\o;\Rb^m)$ function can approximated by a
sequence of functions belonging to $W^{1,p}(\o;\Rb^m) \cap
L^\infty(\o;\Rb^m)$, with respect to the strong
$W^{1,p}(\o;\Rb^m)$-convergence. Then, by the lower semicontinuity
of the $\Gamma$-limsup and the continuity of
$$
(v^+,v^-) \mapsto \int_\o \Q_{n-1} \overline W(D_\a v^+)\, dx_\a +
\int_\o \Q_{n-1} \overline W(D_\a v^-)\, dx_\a + R^{(\ell)}
\int_\o \varphi^{(\ell)} (v^+ - v^-)\, dx_\a
$$
with respect to the strong $W^{1,p}(\o;\Rb^m)$-convergence we get
the thesis for $\ell\in (0,+\infty]$.

\bigskip

If $\ell=0$, we can follow the line of the previous case with
slight changes. Let us start by dealing with Step 1. First, we
have to notice that for the definition of $(\bar u_j)$ in
$B^{n-1}_{\rho_j}(x_i^{\e_j})^{\pm\d_j}$, for $i\in Z_j$, we have
to consider, for any $\eta>0$, a function $\zeta_{\g,j} \in
Y^\g_j(z)$ such that
$$
\int_{(B_{\g N_j}^{n-1} \times I_j) \setminus C_{1,\g N_j}} r_j^p
\, W\left(r_j^{-1} D \zeta_{\g,j} \right) dx \leq
\varphi^{(0)}_{\g,j}(z) +\eta\,;
$$
hence,
$$
\bar u_j(x_\a,x_n):=
\zeta_{\g,j}^i\left(\frac{x_\a-x_i^{\e_j}}{r_j},\frac{x_n}{r_j}\right)+u_j^{i-}\;\text{
in }\; B^{n-1}_{\rho_j}(x_i^{\e_j})^{\pm\d_j}\,, \quad \text{ for
}\, i\in Z_j\,.
$$
While for the definition of $(\bar u_j)$ in
$B^{n-1}_{\rho_j}(x_i^{\e_j})^{\pm\d_j}$, for $i\in
\Zb^{n-1}\setminus Z_j$, we have to introduce a suitable function
$\psi_{\g,j}$ different from the one used in (\ref{defubarZn}). In
fact, for a fixed $\g>0$ and $j$ large enough we can always assume
that $\g N_j >2$ and $\d_j/r_j >2$. Let $\psi \in
W^{1,p}(B_2^{n-1} \times (0,2);[0,1])$ such that $\psi = 0$ on
$B_1^{n-1} \times \{0\}$ and $\psi=1$ on $\partial B_2^{n-1}
\times (0,2)$. We then define
$$
\psi_{\g,j}(x):=\left\{
\begin{array}{lll}
0 & \text{ in } & (B_{\g N_j}^{n-1})^{-(\d_j/r_j)},\\
\psi(x) & \text{ in } & (B_2^{n-1})^{+2},\\
1 & \text{ in } &  (B_{\g N_j}^{n-1})^{+(\d_j/r_j)} \setminus
(B_2^{n-1})^{+2} .
\end{array}
\right.
$$
The functions $\psi_{\g,j}$ belong to $W^{1,p}((B^{n-1}_{\g N_j}
\times I_j) \setminus C_{1,\g N_j};[0,1])$ and satisfy
$\psi_{\g,j} = 1$ on $(\partial B_{\g N_j}^{n-1})^{+(\d_j/r_j)}$
and $\psi_{\g,j} = 0$ in $(B_{\g N_j}^{n-1})^{-(\d_j/r_j)}$.
Hence, we define
$$
\bar u_j:=\psi_{\g,j}
\left(\frac{x_\a-x_i^{\e_j}}{r_j},\frac{x_n}{r_j}
\right)\tilde{w}_j^+ + \Bigg( 1- \psi_{\g,j}
\left(\frac{x_\a-x_i^{\e_j}}{r_j},\frac{x_n}{r_j} \right)
\Bigg)\tilde{w}_j^-
$$
in $ \big( B^{n-1}_{\rho_j}(x_i^{\e_j}) \times (-\d_j,\d_j)\big)
\cap \O_j$ and  for $i\in \Zb^{n-1}\setminus Z_j$. In particular,
we have that $\bar u_j=w_j$ on $\big( \partial
B^{n-1}_{\rho_j}(x_i^{\e_j}) \times (-\d_j,\d_j)\big) \cap \O_j$.

Taking into account the definition of $(\bar u_j)$ we can proceed
as in Steps 2 and 3 also for $\ell=0$.
\hfill$\Box$\\


\section{Representation formula for the interfacial energy
density}\label{nonabstract}

\noindent This section is devoted to describe explicitly the
interfacial energy density $\varphi^{(\ell)}$ for $\ell\in
[0,+\infty]$. As in \cite{Ans}, we expect to find a capacitary
type formula for each regime $\ell\in(0,+\infty)$, $\ell=+\infty$
and $\ell=0$.

We recall that $\varphi^{(\ell)}$ is the pointwise limit of the
sequence $(\varphi^{(\ell)}_{\g,j})$, as $j\to +\infty$ and $\g\to
0^+$ where for $\ell\in (0,+\infty]$
$$
\varphi^{(\ell)}_{\g,j}(z)=\inf \left\{ \int_{(B_{\g N_j}^{n-1}
\times I) \setminus C_{1,\g N_j}} r_j^p \, W\left(r_j^{-1}
\Bigl(D_\a \zeta \Big| {r_j\over\d_j} D_n \zeta \Bigr)\right) dx :
\quad \zeta \in X^\g_j(z)\right\},
$$
while for $\ell=0$,
$$
\varphi^{(0)}_{\g,j}(z)= \inf \left\{ \int_{(B_{\g N_j}^{n-1}
\times I_j ) \setminus C_{1,\g N_j} } r_j^p \, W(r_j^{-1} D
\zeta)\, dx : \quad \zeta \in Y^\g_j(z)\right\}
$$
(see Section \ref{close}). The main difficulty occurring in the
description of $\varphi^{(\ell)}$ is due to the fact that the
above minimum problems are stated on (increasingly) varying
domains. This do not permit, for example, to deal with a direct
$\Gamma$-convergence approach in order to apply the classical
result on the convergence of associated minimum problems. Thus the
proof of the representation formula will be performed in three
main steps: we first prove an auxiliary $\G$-convergence result
for a suitable sequence of energies stated on a fixed domain, then
we describe the functional space occurring in the limit capacitary
formula, finally, we prove that $\varphi^{(\ell)}$ is described by
a representation formula of capacitary-type.
\bigskip

We introduce some convenient notation for the sequel. Let
$g_j:\Rb^{m \times n} \to [0,+\infty)$ be the sequence of
functions given by
$$g_j(F):=r_j^p \, W(r_j^{-1}F)$$
for every $F\in \Rb^{m\times n}$. By (\ref{pgrowth}) and
(\ref{plip}) it follows that
\begin{equation}\label{pgrowthWj}
|F|^p - r_j^p \leq g_j(F) \leq \beta ( r_j^p + |F|^p), \quad
\text{for all }F \in \Rb^{m \times n}
\end{equation}
and the following $p$-Lipschitz condition holds:
$$|g_j(F_1)-g_j(F_2)| \leq c
(r_j^{p-1}+|F_1|^{p-1}+|F_2|^{p-1})|F_1-F_2|, \quad \text{ for all
}\, F_1,\, F_2 \in \Rb^{m \times n}.
$$
Then, according to Ascoli-Arzela's Theorem, up to subsequences,
$g_j$ converges locally uniformly in $\Rb^{m \times n}$ to a
function $g$ satisfying:
\begin{equation}\label{pgrowthg} |F|^p \leq g(F) \leq \beta |F|^p,
\quad \text{for all }\,F \in \Rb^{m \times n}\end{equation}and
\begin{equation}\label{plipg}
|g(F_1)-g(F_2)| \leq c (|F_1|^{p-1}+|F_2|^{p-1})|F_1-F_2|, \quad
\text{ for all }\, F_1,\, F_2\in \Rb^{m \times n}.
\end{equation}

\subsection{The case $\ell \in (0,+\infty)$}

\noindent We define
\begin{eqnarray*}
X_N(z) & := & \Big\{ \zeta \in W^{1,p}((B_N^{n-1} \times I)
\setminus C_{1,N};\Rb^m) :\quad \zeta =z \,\text{ on }\, (
\partial B_N^{n-1} )^+ \\&& \hspace{6cm}\text{ and }
\zeta =0 \,\text{ on }\, ( \partial B_N^{n-1})^- \Big\}
\end{eqnarray*}
for $N>1$ and $I=(-1,1)$. We recall the following
$\Gamma$-convergence result.
\begin{proposition}\label{lfinite} Let
 $$
 \ell=\lim_{j\to +\infty}{r_j\over\d_j} \in (0,+\infty)\,,
 $$
then the sequence of functionals $G^{(\ell)}_j:L^p((B_N^{n-1}
\times I) \setminus C_{1,N};\Rb^m) \to [0,+\infty]$, defined by
$$
G_j^{(\ell)}(\zeta):=\left\{\begin{array}{ll}\ds \int_{(B_N^{n-1}
\times I) \setminus C_{1,N}} g_j \left(D_\a
\zeta\Big|\frac{r_j}{\d_j} D_n \zeta \right)\, dx & \text{ if
}\zeta \in X_N(z)\\ &\\+\infty & \text{
otherwise\,,}\end{array}\right.
$$
$\G$-converges, with respect to the $L^p$-convergence, to
$$
G^{(\ell)}(\zeta):=\left\{\begin{array}{ll}\ds \int_{(B_N^{n-1}
\times I) \setminus C_{1,N}}
g(D_\a \zeta | \ell D_n \zeta )\, dx & \text{ if }\zeta \in X_N(z)\\
&\\ +\infty & \text{ otherwise\,.}\end{array}\right.
$$
\end{proposition}

\noindent {\it Proof. } Since $\ell=\lim_{j\to +\infty} (r_j/\d_j)
\in (0,+\infty)$, by the locally uniform convergence of $g_j$ to
$g$ we have that the sequence of quasiconvex functions $F\mapsto
g_j(\overline F|(r_j/\d_j) F_n)$ pointwise converges to $F\mapsto
g(\overline F|\ell F_n)$. Hence the conclusion comes from
\cite{BD} Propositions 12.8 and 11.7. \hfill$\Box$

\begin{rmk}\label{bochner}{\rm
We denote by $p^*$ the Sobolev exponent in dimension $(n-1)$ i.e.
$$p^*:=\frac{(n-1)p}{n-1-p}.$$ We recall that if $(a,b) \subset \Rb$,
the space $L^p(a,b;L^{p^*}(\Rb^{n-1};\Rb^m))$ is a reflexive and
separable Banach space (see e.g. \cite{Adams} or \cite{Y}). Hence,
by the Banach-Alaoglu-Bourbaki Theorem, any bounded sequence
admits a weakly converging subsequence. }
\end{rmk}

\begin{proposition}[Limit space]\label{bla1}
Let
\begin{equation}\label{assl}
\ell=\lim_{j\to +\infty} {r_j\over \d_j}\in (0, +\infty)\,,\qquad
0<R^{(\ell)}=\lim_{j \to +\infty} \frac{r_j^{n-1-p}}{\e_j^{n-1}} <
+\infty
\end{equation}
and let $(\zeta_{\g,j}) \in X^\g_j(z)$ such that, for every fixed
$\g>0$,
\begin{equation}\label{w1ploc}
\sup_{j \in \Nb}\int_{(B_{\g N_j}^{n-1} \times I) \setminus
C_{1,\g N_j}} g_j\Bigl(D_\a\zeta_{\g,j}\Big| {r_j\over \d_j}D_n
\zeta_{\g,j}\Bigr)\, dx \le c\,.
\end{equation}
Then, there exists a sequence $\tilde\zeta_j \in
W^{1,p}_{\loc}((\Rb^{n-1}\times I)\setminus C_{1,\infty}; \Rb^m)$
such that
$$
\tilde\zeta_j = \zeta_{\g,j} \quad \text{ in } \quad (B^{n-1}_{\g
N_j}\times I)\setminus C_{1, \g N_j}
$$
and such that, up to subsequences, it converges weakly to $\zeta$
in $W^{1,p}_{\rm loc}((\Rb^{n-1}\times I)\setminus
C_{1,\infty};\Rb^m)$. Moreover, the function $\zeta$ satisfies the
following properties
\begin{equation}\label{adm}
\left\{
\begin{array}{l}
D \zeta \in L^p((\Rb^{n-1}\times I)\setminus C_{1,\infty};\Rb^{m
\times
n}),\\
\\
\zeta-z \in L^{p}(0, 1;L^{p^*}(\Rb^{n-1};\Rb^m)),\\
\\
\zeta \in L^{p}(-1,0;L^{p^*}(\Rb^{n-1};\Rb^m))\,.
\end{array}
\right.
\end{equation}
\end{proposition}

\noindent {\it Proof. } By (\ref{pgrowthWj}), (\ref{assl}) and
(\ref{w1ploc}) we deduce that, for every fixed $\g>0$,
\begin{equation}\label{boundgradl}
\sup_{j \in \Nb}\int_{(B_{\g N_j}^{n-1} \times I) \setminus
C_{1,\g N_j}} \Bigl|\Bigl(D_\a\zeta_{\g,j}\Big| {r_j\over \d_j}D_n
\zeta_{\g,j}\Bigr)\Bigr|^p\, dx \le c\,.
\end{equation}
We define
$$
\tilde \zeta_j:= \left\{ {\begin{array}{lll} z & \text{in} &
\big(\mathbb{R}^{n - 1} \setminus B_{\g N_j }^{n - 1} \big)^+,\\
\\
\zeta_{\g,j}  & \text{in} & (B_{\g N_j }^{n - 1}  \times  I)
 \setminus C_{1,\g N_j},\\
 \\
0 & \text{in} & \big(\mathbb{R}^{n - 1} \setminus B_{\gamma N_j
}^{n - 1}  \big)^-\,;
\\
\end{array} } \right.
$$
hence,
$$
\tilde \zeta_j(\cdot,x_n)-z \in W^{1,p}(\Rb^{n-1};\Rb^m)\quad
\text{for \, a.e.\; } x_n\in(0,1)
$$
and
$$
\tilde \zeta_j(\cdot,x_n) \in W^{1,p}(\Rb^{n-1};\Rb^m)\quad
\text{for \, a.e.\; } x_n\in(-1,0)\,.
$$
Moreover by (\ref{boundgradl}) we get that
\begin{equation}\label{stima}
\int_{(\Rb^{n-1} \times I) \setminus C_{1,\infty}}
\Bigl|\Bigl(D_\a\tilde\zeta_{j}\Big| {r_j\over \d_j}D_n
\tilde\zeta_{j}\Bigr)\Bigr|^p \, dx = \int_{(B_{\g N_j}^{n-1}
\times
 I ) \setminus C_{1,\g N_j}} \Bigl|\Bigl(D_\a\zeta_{\g,j}\Big|
{r_j\over \d_j}D_n \zeta_{\g,j}\Bigr)\Bigr|^p \, dx \leq c\,.
\end{equation}
Since $p<n-1$, according to the Sobolev Inequality (see e.g.
\cite{Adams}), there exists a constant $c=c(n,p)>0$ (independent
of $x_n$) such that
\begin{equation}\label{p*1}
\left(\int_{\Rb^{n-1}}|\tilde \zeta_j(x_\a,x_n)-z|^{p^*}\,
dx_\a\right)^{p/p^*} \leq c \,\int_{\Rb^{n-1}} |D_\a \tilde
\zeta_j(x_\a,x_n)|^p\, dx_\a
\end{equation}
for a.e. $x_n \in (0,1)$, and
\begin{equation}\label{p*2}
\left(\int_{\Rb^{n-1}}|\tilde \zeta_j(x_\a,x_n)|^{p^*}\,
dx_\a\right)^{p/p^*} \leq c\, \int_{\Rb^{n-1}} |D_\a \tilde
\zeta_j(x_\a,x_n)|^p\, dx_\a
\end{equation}
for a.e. $x_n \in (-1,0)$. If we integrate (\ref{p*1}) and
(\ref{p*2}) with respect to $x_n$, by (\ref{stima}) and Remark
\ref{bochner}, we get that there exist $\zeta_1\in
L^p(0,1;L^{p^*}(\Rb^{n-1};\Rb^m))$ and $\zeta_2\in
L^p(-1,0;L^{p^*}(\Rb^{n-1};\Rb^m))$ such that, up to subsequences,
$$
\left\{
\begin{array}{lll}
\tilde \zeta_j-z\rightharpoonup \zeta_1 & \text{in} &
L^p(0,1;L^{p^*}(\Rb^{n-1};\Rb^m)),
\\
\tilde \zeta_j\rightharpoonup \zeta_2 & \text{in} &
L^p(-1,0;L^{p^*}(\Rb^{n-1};\Rb^m)),\\
\\
D \tilde \zeta_j \rightharpoonup D \zeta_1 & \text{in} &
L^p((\Rb^{n-1})^+;\Rb^{m \times n}),\\
D \tilde \zeta_j \rightharpoonup D \zeta_2 & \text{in} &
L^p((\Rb^{n-1})^-;\Rb^{m \times n}).
\end{array}
\right.$$ In particular, we have that
$$
\left\{
\begin{array}{lll}
\tilde\zeta_j \rightharpoonup \zeta_1+z  & \text{in} &
W^{1,p}_{\rm loc}((\Rb^{n-1})^+;\Rb^m),\\
\\
\tilde \zeta_j\rightharpoonup\zeta_2 & \text{in} & W^{1,p}_{\rm
loc}((\Rb^{n-1})^-;\Rb^m)\,.
\end{array}
\right.$$ Then, since $\zeta_1+z=\zeta_2$ on $B^{n-1}_1$ in the
sense of traces, we can define
$$
\zeta:= \left\{ {\begin{array}{lll}
{\zeta_1 + z  } & \text{in} & (\Rb^{n - 1})^+  \\
{\zeta_2  } & \text{in} & (\Rb^{n - 1})^-
\cup \big( B^{n-1}_1 \times \{0\} \big), \\
\end{array} } \right.
$$
and it satisfies (\ref{adm}).\hfill$\Box$

\bigskip

Now we are able to describe the interfacial energy density
$\varphi^{(\ell)}$ as the following nonlinear capacitary formula.

\begin{proposition}[Representation formula]\label{reprell}
We have
\begin{eqnarray*}
\varphi^{(\ell)}(z) &= & \inf \Bigg\{\int_{(\Rb^{n-1} \times I)
\setminus C_{1,\infty}}g \big(D_\a \zeta | \ell D_n \zeta \big)\,
dx :\zeta \in W^{1,p}_{\rm loc}((\Rb^{n-1} \times I) \setminus
C_{1,\infty};\Rb^m),\\
&&\hspace{0.6cm} D \zeta \in L^p((\Rb^{n-1} \times I)\setminus
C_{1,\infty};\Rb^{m \times n}),\; \zeta - z \in
L^{p}(0,1;L^{p^*}(\Rb^{n-1};\Rb^m))\\
&&\hspace{5.8cm} \text{ and } \zeta \in
L^{p}(-1,0;L^{p^*}(\Rb^{n-1};\Rb^m))\Bigg\}
\end{eqnarray*}
for every $z\in \Rb^m$.
\end{proposition}

\noindent {\it Proof. } We define
\begin{eqnarray*}
\psi^{(\ell)}(z) &:= & \inf \Bigg\{\int_{(\Rb^{n-1} \times I)
\setminus C_{1,\infty}}g \big(D_\a \zeta | \ell D_n \zeta \big)\,
dx :\zeta \in W^{1,p}_{\rm loc}((\Rb^{n-1} \times I) \setminus
C_{1,\infty};\Rb^m),\\
&&\hspace{0.6cm} D \zeta \in L^p((\Rb^{n-1} \times I)\setminus
C_{1,\infty};\Rb^{m \times n}),\; \zeta - z \in
L^{p}(0,1;L^{p^*}(\Rb^{n-1};\Rb^m))\\
&&\hspace{5.8cm} \text{ and } \zeta \in
L^{p}(-1,0;L^{p^*}(\Rb^{n-1};\Rb^m))\Bigg\},
\end{eqnarray*}
we want to prove that $\varphi^{(\ell)}(z)= \psi^{(\ell)}(z)$ for
every $z\in \Rb^m$. For every fixed $\eta>0$, by definition of
$\varphi_{\g,j}^{(\ell)}(z)$ (see (\ref{phigj})), there exists
$\zeta_{\g,j} \in X^\g_j(z)$ such that
$$
\int_{(B_{\g N_j}^{n-1} \times I) \setminus C_{1,\g N_j}}g_j\left(
D_\a \zeta_{\g,j} \Big| \frac{r_j}{\d_j} D_n \zeta_{\g,j} \right)
dx \leq \varphi_{\g,j}^{(\ell)}(z) + \eta.
$$
By Proposition \ref{midi}(i) we have that (\ref{w1ploc}) is
fulfilled, then by Propositions \ref{bla1} and \ref{lfinite} we
get
\begin{eqnarray*}
\lim_{j \to +\infty} \varphi_{\g,j}^{(\ell)}(z) +\eta & \geq &
\liminf_{j \to +\infty} \int_{(B_{\g N_j}^{n-1} \times I)
\setminus C_{1,\g N_j}} g_j\left(D_\a \tilde\zeta_j
\Big| \frac{r_j}{\d_j} D_n\tilde\zeta_j \right) dx\nonumber\\
& \geq & \liminf_{j \to +\infty} \int_{(B_N^{n-1} \times
I)\setminus C_{1,N}}g_j\left( D_\a \tilde\zeta_j \Big|
\frac{r_j}{\d_j} D_n \tilde\zeta_j
\right) dx\\
& \geq & \int_{(B_N^{n-1} \times I)\setminus C_{1,N}}g( D_\a \zeta
| \ell D_n \zeta)\, dx
\end{eqnarray*}
with $\zeta \in W^{1,p}_{\rm loc}((\Rb^{n-1} \times I) \setminus
C_{1,\infty};\Rb^m)$ satisfying (\ref{adm}). Note that for every
fixed $\g>0$ and $j$ large enough we can always assume that $\g
N_j>N$ for some fixed $N>2$. Hence, passing to the limit as $N
\to+\infty$ and $\g \to 0^+$, we obtain
\begin{equation}\label{ineq1}
\varphi^{(\ell)}(z) + \eta \geq \int_{(\Rb^{n-1} \times I)
\setminus C_{1,\infty}}g (D_\a \zeta | \ell D_n \zeta )\, dx
\geq\psi^{(\ell)}(z)
\end{equation}
and by the arbitrariness of $\eta$ we get the first inequality.

We now prove the converse inequality. For every fixed $\eta>0$
there exists $\zeta \in W^{1,p}_{\rm loc}((\Rb^{n-1} \times I)
\setminus C_{1,\infty};\Rb^m)$ satisfying (\ref{adm}) such that
\begin{equation}\label{almostmin}
\int_{(\Rb^{n-1} \times I) \setminus C_{1,\infty}} g(D_\a \zeta |
\ell D_n\zeta)\, dx \leq \psi^{(\ell)}(z) + \eta\,.
\end{equation}
Let $N>2$, for every fixed $\g>0$ and $j$ large enough we have
that $\g N_j>N$. We consider a cut-off function $\theta_N \in
\C^\infty_c(B^{n-1}_N;[0,1])$ such that $\theta_N=1$ in
$B^{n-1}_{N/2}$, $|D_\a\theta_N| \leq c/N$ and we define
$$
\zeta_N:= \left\{ \begin{array}{lll} \theta _N (x_\a )\zeta + (1 -
\theta _N (x_\a))z &
\text{in} & (B^{n-1}_N)^{+},\\
&&\\
 \theta _N (x_\alpha )\zeta & \text{in} & (B^{n-1}_N)^{-}\cup
(B^{n-1}_1 \times \{0\})
\end{array} \right.
$$
so that $\zeta_N \in X_N(z)$. By Proposition \ref{lfinite}, there
exists a sequence ($\zeta_j^N) \subset X_N(z)$ strongly converging
to $\zeta_N$ in $L^p((B^{n-1}_N \times I) \setminus
C_{1,N};\Rb^m)$ such that
\begin{equation}\label{contrex}
\int_{(B^{n-1}_N \times I) \setminus C_{1,N}} g ( D_\a \zeta_N |
\ell D_n \zeta_N)\, dx = \lim_{j \to +\infty} \int_{(B^{n-1}_N
\times I) \setminus C_{1,N}} g_j \left( D_\a \zeta_j^N \Big|
\frac{r_j}{\d_j} D_n \zeta_j^N \right)\, dx
\end{equation}
Let us define $\zeta_{\g,j} \in X_j^\g (z)$ as
$$
\zeta_{\g,j}:= \left\{ \begin{array}{lll}
z & \text{in} & (B^{n-1}_{\g N_j} \setminus B^{n-1}_N)^+,\\
&&\\
\zeta_j^N &
\text{in} & (B^{n-1}_N \times I) \setminus C_{1,N},\\
&&\\
0 & \text{in} & (B^{n-1}_{\g N_j} \setminus B^{n-1}_N)^-.
\end{array} \right.
$$
Consequently, $\zeta_{\g,j}$ is an admissible test function for
(\ref{phigj}) and since $g_j(0)=0$ we get that
\begin{eqnarray*}
\varphi_{\g,j}^{(\ell)}(z) & \leq & \int_{(B^{n-1}_{\g N_j} \times
I) \setminus C_{1,\g N_j}} g_j\Bigl( D_\a \zeta_{\g,j} \Big|
\frac{r_j}{\d_j} D_n \zeta_{\g,j} \Bigr)\, dx\\
 & = & \int_{(B^{n-1}_N \times I) \setminus C_{1,N}} g_j\Bigl(D_\a
\zeta_N^j \Big| \frac{r_j}{\d_j} D_n \zeta_N^j \Bigr)\, dx.
\end{eqnarray*}
Passing to the limit as $j \to +\infty$, using (\ref{contrex}) and
the $p$-growth condition (\ref{pgrowthg}) satisfied by $g$, we
obtain
\begin{eqnarray}\label{est1}
\lim_{j \to +\infty} \varphi_{\g,j}^{(\ell)}(z) & \leq  &
\int_{(B^{n-1}_N \times I) \setminus C_{1,N}} g (D_\a \zeta_N |
\ell
D_n \zeta_N )\, dx \nonumber\\
&\le& \int_{(B^{n-1}_{N/2} \times I) \setminus C_{1,N/2}} g(D_\a
\zeta | \ell D_n \zeta )\, dx +c \int_{(B^{n-1}_N \setminus
B^{n-1}
_{N/2})^+} |D \zeta_N |^p\, dx \nonumber\\
&& \hspace{1cm} +c \int_{(B^{n-1}_N \setminus B^{n-1}_{N/2})^-} |D
\zeta_N |^p\dx \,.
\end{eqnarray} Let us examine the contribution of the gradient in
(\ref{est1}),
\begin{eqnarray}\label{est2}
&&\nonumber\int_{(B^{n-1}_N \setminus B^{n-1}_{N/2})^+}|D \zeta_N
|^p\dx
+ \int_{(B^{n-1}_N \setminus B^{n-1}_{N/2})^-} |D \zeta_N |^p\dx\\
&&\hspace{2cm}\le\nonumber c \int_{(B^{n-1}_N \setminus
B^{n-1}_{N/2})^+}(|D_\a \theta_N|^p |\zeta - z|^p + |D \zeta|^p
)\,dx\\
&&\hspace{3cm} \nonumber+ c \int_{(B^{n-1}_N \setminus
B^{n-1}_{N/2})^-}(|D_\a \theta_N|^p |\zeta |^p
+ |D \zeta|^p )\, dx\\
&&\hspace{2cm}\le \nonumber c\,\left( \int_{(\Rb^{n-1}\setminus
B^{n-1}_{N/2})^+} |D \zeta|^p \dx +
\int_{(\Rb^{n-1}\setminus B^{n-1}_{N/2})^-} |D \zeta|^p \dx \right)\\
&&\hspace{3cm}+ \frac{c}{N^p}\left(\int_{(B^{n-1}_N \setminus
B^{n-1}_{N/2})^+} |\zeta-z|^p\dx + \int_{(B^{n-1}_N \setminus
B^{n-1}_{N/2})^-} |\zeta|^p\dx\right)\,.
\end{eqnarray} Since $p^*>p$ we can apply H\"older Inequality
with $q=p^*/p$ obtaining
\begin{eqnarray}\label{est3}
&&\hspace{-1cm}\nonumber \frac{c}{N^p}\left(\int_{(B^{n-1}_N
\setminus B^{n-1}_{N/2})^+} |\zeta-z|^p\dx + \int_{(B^{n-1}_N
\setminus B^{n-1}_{N/2})^-}
|\zeta|^p\right)\\
&&\hspace{1cm}\le \nonumber c\, \int_0^1\Bigl(\int_{B^{n-1}_N
\setminus B^{n-1}_{N/2}}
|\zeta-z|^{p^*}\dx_\a\Bigr)^{p/p^*}\dx_n\\
&&\hspace{5cm} \nonumber  + c\, \int_{-1}^0\Bigl( \int_{B^{n-1}_N
\setminus
B^{n-1}_{N/2}} |\zeta|^{p^*}\dx_\a\Bigr)^{p/p^*}\dx_n\\
&& \hspace{1cm}\le\nonumber \nonumber c\,
\int_0^1\Bigl(\int_{\Rb^{n-1} \setminus B^{n-1}_{N/2}}
|\zeta-z|^{p^*}\dx_\a\Bigr)^{p/p^*}\dx_n\\
&& \hspace{5cm}  + c\, \int_{-1}^0\Bigl( \int_{\Rb^{n-1} \setminus
B^{n-1}_{N/2}} |\zeta|^{p^*}\dx_\a\Bigr)^{p/p^*}\dx_n.
\end{eqnarray} Hence by (\ref{adm}), (\ref{est2}) and (\ref{est3})
we have that, for every fixed $\g>0$,
$$
\lim_{N \to +\infty}\int_{(B^{n-1}_N \setminus B^{n-1}_{N/2})^\pm}
|D \zeta_N|^p\dx=0
$$
which thanks to (\ref{almostmin}) and (\ref{est1}) implies that
$$
\lim_{j\to +\infty} \varphi_{\g,j}^{(\ell)}(z) \le
\psi^{(\ell)}(z) + \eta.
$$
Then we get the converse inequality by letting $\g \to 0^+$ and by
the arbitrariness of $\eta$.\hfill$\Box$

\subsection{The case $\ell=+\infty$}

\noindent In this case the study leading to the representation
formula for $\varphi^{(\infty)}$ involves a dimensional reduction
problem stated on a varying domain. As before, we start proving
some $\Gamma$-convergence results (see Propositions \ref{separate}
and \ref{linfinite}) for suitable sequences of functionals stated
on fixed domains. This will allow as to apply some well-known
$\Gamma$-convergence and integral representation theorems due to
Le Dret-Raoult \cite{LDR} and Braides-Fonseca-Francfort \cite{BFF}
respectively.

\smallskip

Let $G_j^\pm :L^p((B_N^{n-1})^\pm;\Rb^m) \to [0,+\infty]$ be
defined by
$$
G_j^+(\zeta):=\left\{\begin{array}{ll}\ds \int_{(B_N^{n-1})^+}
g_j\left(D_\a \zeta\Big|\frac{r_j}{\d_j} D_n \zeta \right)\, dx &
\text{ if } \left\{\begin{array}{l}\zeta \in
W^{1,p}((B_N^{n-1})^+;\Rb^m)\\
\zeta=z \text{ on } (\partial B_N^{n-1})^+\end{array}\right.\\&\\
+\infty & \text{ otherwise}\end{array}\right.
$$
and
$$
G_j^-(\zeta):=\left\{\begin{array}{ll}\ds \int_{(B_N^{n-1})^-}
g_j\left(D_\a \zeta\Big|\frac{r_j}{\d_j} D_n \zeta \right)\, dx &
\text{ if }\left\{\begin{array}{l}\zeta \in
W^{1,p}((B_N^{n-1})^-;\Rb^m)\\
\zeta=0 \text{ on } (\partial
B_N^{n-1})^-\end{array}\right.\\&\\+\infty & \text{
otherwise}.\end{array}\right.
$$
\begin{proposition}\label{separate}
Let
$$
\ell=\lim_{j\to +\infty}{r_j\over \d_j}=+\infty\,,
$$
then, the sequences of functionals $(G_j^\pm$) $\G$-converge, with
respect to the $L^p$-convergence, to
$$G^+(\zeta):=\left\{\begin{array}{ll}\ds \int_{B_N^{n-1}} \Q_{n-1}\,
\overline g (D_\a \zeta)\, dx_\a&\text{ if } \zeta-z \in
W^{1,p}_0(B_N^{n-1};\Rb^m)\\&
\\ +\infty & \text{ otherwise}\end{array}\right.
$$
and
$$
G^-(\zeta):=\left\{\begin{array}{ll}\ds \int_{B_N^{n-1}}
\Q_{n-1}\, \overline g (D_\a \zeta)\, dx_\a&\text{ if }\zeta \in
W^{1,p}_0(B_N^{n-1};\Rb^m)\\
&\\+\infty & \text{ otherwise}\,,\end{array}\right.
$$
respectively, where $\overline g(\overline F) = \inf \{
g(\overline F|F_n): \; F_n \in \Rb^m \}$ for every $\overline F
\in \Rb^{m\times (n-1)}$.
\end{proposition}

\noindent {\it Proof. }We prove the $\Gamma$-convergence result
only for ($G_j^+$), the other one being analogous. According to
\cite{BFF} Theorem 2.5 and Lemma 2.6 there exists a continuous
function $\hat g :\Rb^{m \times (n-1)} \to [0,+\infty)$ such that,
up to subsequence, ($G_j^+$) $\G$-converges to
$$
G^+(\zeta):=\left\{\begin{array}{ll}\ds \int_{B_N^{n-1}} \hat g
(D_\a \zeta)\, dx_\a &\text{ if } \zeta-z \in
W^{1,p}_0(B_N^{n-1};\Rb^m)\\&\\ +\infty & \text{otherwise}\,.
\end{array}\right.
$$
Hence, it remains to show that $\hat g=\Q_{n-1}\, \overline g$. By
\cite{BFF} Lemma 2.6,  it is enough to consider
$W^{1,p}$-functions without boundary condition; hence, it will
suffice to deal with affine functions. Let $\zeta (x_\a):=
\overline F \cdot x_\a$, by \cite{BFF} Theorem 2.5, there exists a
sequence $(\zeta_j) \subset W^{1,p}((B_N^{n-1})^+;\Rb^m)$ (the
so-called recovery sequence) converging to $\zeta$ in
$L^p((B_N^{n-1})^+;\Rb^m)$, such that
\begin{equation}\label{1239}
\hat g(\overline F)\, c_N =G^+(\zeta)= \lim_{j \to +\infty}
\int_{(B_N^{n-1})^+} g_j\left(D_\a \zeta_j \Big| \frac{r_j}{\d_j}
D_n \zeta_j\right)\, dx
\end{equation}
where $c_N= \HH^{n-1}(B_N^{n-1})$. Moreover, by \cite{Bo&Fo}
Theorem 1.1, we can assume, without loss of generality, that the
sequence $\big(\big|\big(D_\a\zeta_j | \frac{r_j}{\d_j} D_n
\zeta_j \big)\big|^p \big)$ is equi-integrable. By (\ref{1239})
and (\ref{pgrowthWj}), we have that
$$
\sup_{j \in \Nb} \int_{(B_N^{n-1})^+} \Bigl|\Bigl(D_\a \zeta_j
\Big| \frac{r_j}{\d_j} D_n \zeta_j\Bigr)\Bigr|^p\dx \le c\,;
$$
hence, for every fixed $M>0$, if we define
$$
A_j^M:=\left\{x \in (B_N^{n-1})^+ : \quad \left| \left(D_\a
\zeta_j(x) \Big|\frac{r_j}{\d_j} D_n \zeta_j(x) \right)
\right|\leq M \right\}\,,
$$
we get that $\mathcal L^n((B_N^{n-1})^+ \setminus A_j^M) \leq
c/M^p$ for some constant $c>0$ independent of $j$ and $M$. Fix
$M>0$, by (\ref{1239}), we have
\begin{equation}\label{1325}
\hat g(\overline F)\, c_N \geq \limsup_{j \to +\infty}
\int_{A_j^M} g_j\left(D_\a\zeta_j \Big| \frac{r_j}{\d_j} D_n
\zeta_j \right)\, dx.
\end{equation}
Moreover, for all $x \in A_j^M$,
$$
\left| g_j\left(D_\a \zeta_j(x) \Big| \frac{r_j}{\d_j}
D_n\zeta_j(x) \right) - g\left(D_\a \zeta_j(x)
\Big|\frac{r_j}{\d_j} D_n \zeta_j(x) \right) \right| \leq
\sup_{|F|\leq M} |g_j(F)-g(F)|,
$$
and then,
\begin{eqnarray*}
&&\int_{A_j^M} \left| g_j\left(D_\a\zeta_j \Big| \frac{r_j}{\d_j}
D_n \zeta_j \right) - g\left(D_\a \zeta_j \Big| \frac{r_j}{\d_j}
D_n\zeta_j \right)
\right|dx\nonumber\\
&\leq& c_N\, \sup_{|F| \leq M}|g_j(F)-g(F)|.
\end{eqnarray*}
Hence, by the local uniform convergence of $g_j$ to $g$, we have
that
$$
\lim_{j\to +\infty}\int_{A_j^M}\left(g_j\left(D_\a\zeta_j \Big|
\frac{r_j}{\d_j} D_n \zeta_j \right) - g\left(D_\a \zeta_j \Big|
\frac{r_j}{\d_j} D_n\zeta_j \right)\right)\dx=0.
$$
By (\ref{1325}), we get
\begin{equation}\label{1330}
\hat g(\overline F) \, c_N \geq \limsup_{j \to +\infty}
\int_{A_j^M} g\left(D_\a\zeta_j \Big| \frac{r_j}{\d_j} D_n \zeta_j
\right)\, dx.
\end{equation}
Note that, since $\mathcal L^n((B_N^{n-1})^+ \setminus A_j^M) \to
0$ as $M \to +\infty$, by the $p$-growth condition
(\ref{pgrowthg}) and the equi-integrability assumption, we find
\begin{equation}\label{1331}
\limsup_{j \to +\infty} \int_{(B_N^{n-1})^+ \setminus A_j^M}
g\left(D_\a \zeta_j \Big|\frac{r_j}{\d_j} D_n \zeta_j \right)\,
dx= o(1)\,, \qquad \text{as}\quad M\to +\infty\,.
\end{equation}
Consequently, (\ref{1330}) and (\ref{1331}) imply that
\begin{equation}\label{1332}
\hat g(\overline F)\, c_N \geq \limsup_{j \to +\infty}
\int_{(B_N^{n-1})^+} g\left(D_\a \zeta_j \Big|\frac{r_j}{\d_j} D_n
\zeta_j \right)\, dx.
\end{equation}
Finally, from \cite{LDR} Theorem 2, we  know that
$$
\liminf_{j \to +\infty} \int_{(B_N^{n-1})^+} g\left(D_\a \zeta_j
\Big|\frac{r_j}{\d_j} D_n \zeta_j \right)\, dx \ge
\Q_{n-1}\,\overline g(\overline F)\, c_N\,;
$$
hence, by (\ref{1332}) we obtain that $\hat g(\overline F) \geq
\Q_{n-1}\,\overline g(\overline F)$.

We  now prove the converse inequality. By \cite{LDR} Theorem 2,
there exists a sequence $(\zeta_j)$ belonging to
$W^{1,p}((B_N^{n-1})^+;\Rb^m)$ and converging to $\zeta$ in
$L^p((B_N^{n-1})^+;\Rb^m)$ such that
\begin{equation}\label{cat}
\Q_{n-1} \overline g(\overline F)\, c_N =\lim_{j \to +\infty}
\int_{(B_N^{n-1})^+} g\left(D_\a \zeta_j \Big| \frac{r_j}{\d_j}
D_n \zeta_j\right)\, dx\,.
\end{equation}
Without loss of generality, we can still assume that the sequence
$\big(\big|\big(D_\a \zeta_j | \frac{r_j}{\d_j}D_n \zeta_j
\big)\big|^p\big)$ is equi-integrable. Thus arguing as above, from
(\ref{cat}) we deduce
\begin{equation}\label{1333}
\Q_{n-1}\,  \overline g(\overline F)\, c_N \geq \limsup _{j \to
+\infty} \int_{(B_N^{n-1})^+} g_j\left(D_\a \zeta_j \Big|
\frac{r_j}{\d_j} D_n\zeta_j \right)\, dx\,.
\end{equation}
Now, by \cite{BFF} Theorem 2.5, we have that
$$
\liminf _{j \to +\infty} \int_{(B_N^{n-1})^+} g_j\left(D_\a
\zeta_j \Big| \frac{r_j}{\d_j} D_n\zeta_j \right)\, dx \ge \hat
g(\overline F)\, c_N\,;
$$
hence,
 $\Q_{n-1}\,
\overline g(\overline F) \geq \hat g(\overline F)$, which
concludes the proof. \hfill$\Box$

\begin{rmk}\label{bcxn=0}
{\rm By \cite{LDR} Theorem 2, for every $\zeta \in
W^{1,p}(B_N^{n-1};\Rb^m)$ the recovery sequence is given by
$\zeta_j(x_\a,x_n):=\zeta(x_\a) + (\d_j/r_j)\, x_n\, b_j(x_\a)$
for a suitable sequence of functions $(b_j) \subset
\C^\infty_c(B_N^{n-1};\Rb^m)$. Note that by definition $(\zeta_j)$
keeps the boundary conditions of $\zeta$. Reasoning as in the
proof of Proposition \ref{separate} we can observed that
$(\zeta_j)$ is also a recovery sequence for $(G_j^+)$ (see e.g.
(\ref{1333})). The same remark holds for  $(G_j^-)$.}
\end{rmk}

\bigskip

\begin{proposition}\label{linfinite}
Let
$$
\ell= \lim_{j\to +\infty} {r_j\over\d_j}= +\infty\,,
$$
then the sequence of functionals $G^{(\infty)}_j:L^p((B_N^{n-1}
\times I) \setminus C_{1,N};\Rb^m) \to[0,+\infty]$ defined by
$$
G^{(\infty)}_j(\zeta):=\left\{\begin{array}{ll}\ds
\int_{(B_N^{n-1} \times I) \setminus C_{1,N}} g_j \left(D_\a
\zeta\Big|\frac{r_j}{\d_j} D_n \zeta \right)\, dx & \text{ if
}\zeta \in X_N(z)\\ & \\+\infty & \text{
otherwise}\end{array}\right.
$$
$\G$-converges, with respect to the $L^p$-convergence, to
$$
G^{(\infty)}(\zeta):=\left\{\begin{array}{ll}\ds \int_{(B_N^{n-1}
\times I) \setminus C_{1,N}} \Q_{n-1}\, \overline g(D_\a \zeta)\,
dx & \text{ if }\zeta \in X_N(z)\text{ and
}D_n \zeta=0\\\\
+\infty & \text{ otherwise}\,.\end{array}\right.
$$
\end{proposition}
\noindent {\it Proof. } The $\liminf$ inequality is a
straightforward consequence of Proposition \ref{separate}.

Dealing with the $\limsup$ inequality, let us consider $\zeta \in
X_N(z)$ with $D_n \zeta=0$. We denote by $\zeta^\pm \in
W^{1,p}(B_N^{n-1}(0);\Rb^m)$ the restriction of $\zeta$ to
$(B_N^{n-1})^+$ and $(B_N^{n-1})^-$, respectively. By Proposition
\ref{separate} and Remark \ref{bcxn=0}, there exist two sequences
$(\zeta_j^\pm)  \subset W^{1,p} ((B_N^{n-1})^\pm;\Rb^m)$ such that
\begin{equation}\label{2128}
\begin{array}{l}\zeta_j^+ \to \zeta^+ \text{ in
}L^p((B_N^{n-1})^+;\Rb^m)\,,
\quad \zeta_j^+=z \text{ on } ( \partial B_N^{n-1})^+\\\\
\zeta_j^-\to \zeta^- \text{ in }L^p((B_N^{n-1})^-;\Rb^m)\,, \quad
\zeta_j^-=0 \text{ on } (
\partial B_N^{n-1})^-
\end{array}
\end{equation}
and
\begin{eqnarray}\label{2129}
\nonumber \lim_{j \to +\infty} \int_{(B_N^{n-1})^+} g_j\left(D_\a
\zeta^+_j \Big| \frac{r_j}{\d_j} D_n\zeta^+_j \right)\, dx&=&
\int_{B_N^{n-1}} \Q_{n-1}\, \overline g(D_\a \zeta^+)\,
dx_\a\\
\lim_{j \to +\infty} \int_{(B_N^{n-1})^-} g_j\left(D_\a \zeta^-_j
\Big| \frac{r_j}{\d_j} D_n\zeta^-_j \right)\, dx&=&
\int_{B_N^{n-1}} \Q_{n-1}\, \overline g(D_\a \zeta^-)\, dx_\a\,.
\end{eqnarray}
Moreover, since $\zeta \in W^{1,p}((B^{n-1}_N \times I)\setminus
C_{1,N};\Rb^m)$, by Remark \ref{bcxn=0}, $(\zeta_j^+)$ and
$(\zeta_j^-)$ have the same trace on $B^{n-1}_1 \times \{0\}$;
hence, $\zeta_j^+=\zeta^-_j=\zeta$ on $B^{n-1}_1 \times \{0\}$.
Then we can define
$$
\bar \zeta_j:= \left\{
\begin{array}{lll}
\zeta_j^+  & \text{in} & (B_N^{n - 1})^+,\\
\zeta & \text{on} & B_1^{n-1} \times \{0\},\\
\zeta_j^- & \text{in} & (B_N^{n - 1})^-,  \\
\end{array}  \right.
$$
with $\bar \zeta_j \in W^{1,p}((B^{n-1}_N \times I)\setminus
C_{1,N};\Rb^m)$. In particular, by (\ref{2128}) we have that $\bar
\zeta_j \in  X_N(z)$ and $\bar\zeta_j \to \zeta$ in
$L^p((B_N^{n-1} \times I) \setminus C_{1,N};\Rb^m)$. Finally, by
(\ref{2129}) , we have
\begin{eqnarray*}
\lim_{j \to +\infty} G^{(\infty)}_j(\bar \zeta_j) & = & \lim_{j
\to+\infty}\int_{(B_N^{n-1} \times I) \setminus C_{1,N}}
g_j\left(D_\a\bar\zeta_j\Big| \frac{r_j}{\d_j} D_n \bar\zeta_j
\right)\, dx\\
& = &\int_{B_N^{n-1}} \Q_{n-1}\, \overline g(D_\a \zeta^+)\,
dx_\a+ \int_{B_N^{n-1}} \Q_{n-1}\, \overline g(D_\a \zeta^-)\,
dx_\a\\ & = & \int_{(B_N^{n-1} \times I) \setminus C_{1,N}}
\Q_{n-1}\, \overline g (D_\a \zeta)\, dx
\end{eqnarray*} which completes the proof of the
$\limsup$ inequality. \hfill$\Box$

\bigskip

\begin{proposition}[Limit space]\label{bla2}
Let
$$
\ell=\lim_{j\to +\infty} {r_j\over \d_j}= +\infty\,,\qquad
0<R^{(\infty)}=\lim_{j \to +\infty} \frac{r_j^{n-1-p}}{\e_j^{n-1}}
< +\infty
$$
and let $\zeta_{\g,j} \in  X^\g_j(z)$ such that, for every fixed
$\g>0$,
\begin{equation}\label{boundgradinfinity}
\sup_{j \in \Nb}\int_{(B_{\g N_j}^{n-1} \times I) \setminus
C_{1,\g N_j}} g_j\Bigl(D_\a\zeta_{\g,j}\Big| {r_j\over \d_j}D_n
\zeta_{\g,j}\Bigr)\, dx \le c\,.
\end{equation}
Then, there exists a sequence $\tilde\zeta_j \in
W^{1,p}_{\loc}((\Rb^{n-1}\times I)\setminus C_{1,\infty}; \Rb^m)$
such that
$$
\tilde\zeta_j= \zeta_{\g,j} \quad \text{ in }\quad (B^{n-1}_{\g
N_j} \times I)\setminus C_{1, \g N_j}
$$
and such that, up to subsequences, it converges weakly to
$\zeta^+$ in $W^{1,p}_{\rm loc}((\Rb^{n-1})^+;\Rb^m)$ and to
$\zeta^-$ in $W^{1,p}_{\rm loc}((\Rb^{n-1})^- ;\Rb^m)$. Moreover,
the functions $\zeta^\pm$ satisfy the following properties
$$\left\{
\begin{array}{l}
\zeta^\pm \in W^{1,p}_{\rm loc}(\Rb^{n-1};\Rb^m),\\
\\
\zeta^+=\zeta^- \quad\text{ in }\, B_1^{n-1},\\
\\
D_\a \zeta^\pm \in L^p(\Rb^{n-1};\Rb^{m\times (n-1)}),\\
\\
(\zeta^+ -z) \text{ and } \zeta^- \in L^{p^*}(\Rb^{n-1};\Rb^m).
\end{array}
\right.$$
\end{proposition}

\noindent {\it Proof. } We can reason as in Proposition \ref{bla1}
using the fact that, by (\ref{boundgradinfinity}),
$$
\int_{(\Rb^{n-1})^{\pm}}|D_n \tilde\zeta_j|^p\, dx \leq c\,
\Bigl(\frac{\d_j}{r_j}\Bigr)^p\,;
$$
hence, in the limit we have that $D_n \zeta=0$ a.e. in
$(\Rb^{n-1})^{\pm}$. \hfill$\Box$

\begin{proposition}[Representation formula]\label{reprinf}
We have
\begin{eqnarray*}
\varphi^{(\infty)}(z) &= & \inf \Bigg\{\int_{\Rb^{n-1}}\big(
\Q_{n-1}\, \overline g (D_\a \zeta^+ ) + \Q_{n-1}\, \overline g
(D_\a\zeta^- ) \big)\, dx_\a :\; \zeta^\pm \in W^{1,p}_{\rm
loc}(\Rb^{n-1};\Rb^m),\\
&&\hspace{3.8cm} \zeta^+=\zeta^- \text{ in }B_1^{n-1}, \quad
D_\a\zeta^\pm\in L^p(\Rb^{n-1} ;\Rb^{m \times
(n-1)}),\\
&&\hspace{6cm} (\zeta^+ - z) \text{ and } \zeta^- \in
L^{p^*}(\Rb^{n-1};\Rb^m) \Bigg\}
\end{eqnarray*}
for every $z\in \Rb^m$.
\end{proposition}
\noindent {\it Proof. } Reasoning as in the proof of Proposition
\ref{reprell}, by Propositions \ref{linfinite} and \ref{bla2} we
get the representation formula for $\varphi^{(\infty)}$.
\hfill$\Box$

\bigskip

\subsection{The case $\ell=0$}

\noindent We first recall the following $\Gamma$-convergence
result.
\begin{proposition}\label{lzero} The sequence of functionals
$G^{(0)}_j:L^p((B_N^{n-1} \times (-N,N)) \setminus C_{1,N};\Rb^m)
\to [0,+\infty]$, defined by
$$
G^{(0)}_j(\zeta):=\left\{\begin{array}{ll}\ds \int_{(B_N^{n-1}
\times (-N,N)) \setminus C_{1,N}} g_j (D \zeta )\, dx & \text{ if
}\zeta \in
W^{1,p}((B_N^{n-1} \times (-N,N)) \setminus C_{1,N} ;\Rb^m),\\
&\\+\infty & \text{ otherwise\,,}
\end{array}\right.
$$
$\G$-converges, with respect to the $L^p$-convergence, to
$$
G^{(0)}(\zeta):=\left\{\begin{array}{ll}\ds \int_{(B_N^{n-1}
\times (-N,N)) \setminus C_{1,N}} g (D \zeta)\, dx & \text{if }
\zeta \in
W^{1,p}((B_N^{n-1} \times (-N,N)) \setminus C_{1,N} ;\Rb^m),\\&\\
+\infty & \text{ otherwise\,.}
\end{array}\right.
$$
\end{proposition}

\noindent {\it Proof. }The result is an immediate consequence of
the pointwise convergence of the sequence of quasiconvex functions
$g_j$ towards $g$ together with Proposition 12.8 in
\cite{BD}.\hfill$\Box$

\bigskip

\begin{proposition}[Limit space]\label{bla3} Let
\begin{equation}\label{ass0}
\ell=\lim_{j\to +\infty} {r_j\over \d_j}= 0\,,\qquad
0<R^{(0)}=\lim_{j \to +\infty} \frac{r_j^{n-p}}{\e_j^{n-1} \d_j} <
+\infty
\end{equation}
and let $\zeta_{\g,j} \in Y^\g_j(z)$ such that, for every fixed
$\g>0$,
\begin{equation}\label{boundgj}
\sup_{j \in \Nb}\int_{(B_{\g N_j}^{n-1} \times I_j) \setminus
C_{1,\g N_j}} g_j(D \zeta_{\g,j})\, dx \le c\,.
\end{equation}
Then, there exists a sequence $\tilde\zeta_j \in
W^{1,p}_{\loc}(\Rb^n\setminus C_{1,\infty}; \Rb^m)$ such that
$$
\tilde\zeta_j= \zeta_{\g,j}\quad \text{ in } \quad (B^{n-1}_{\g
N_j} \times I_j) \setminus C_{1, \g N_j}
$$
and such that, up to subsequences, it converges weakly to $\zeta$
in $W^{1,p}_{\rm loc}(\Rb^n\setminus C_{1,\infty};\Rb^m)$.
Moreover, the function $\zeta$ satisfies the following properties
\begin{equation}\label{defzeta}
\left\{
\begin{array}{l}
D \zeta \in L^p(\Rb^n\setminus C_{1,\infty};\Rb^{m \times n}),\\
\\
\zeta-z \in L^{p}(0,+\infty;L^{p^*}(\Rb^{n-1};\Rb^m)),\\
\\
\zeta \in L^{p}(-\infty,0;L^{p^*}(\Rb^{n-1};\Rb^m))\,.
\end{array}
\right.
\end{equation}
\end{proposition}

\noindent {\it Proof.} By (\ref{boundgj}), (\ref{pgrowthWj}) and
(\ref{ass0}), we deduce that, for every fixed $\g>0$,
\begin{equation}\label{boundgrad0}
\sup_{j \in \Nb}\int_{(B_{\g N_j}^{n-1} \times I_j) \setminus
C_{1,\g N_j}} |D \zeta_{\g,j}|^p\, dx \le c\,.
\end{equation}
Let us first extend $\zeta_{\g,j}$ by reflection
\begin{equation}\label{barzetagj}
\bar\zeta_{\g,j}(x)=\left\{
\begin{array}{lll}
\zeta_{\g,j}\Bigl(x_\a, 2{\d_j\over r_j}-x_n\Bigr)& \text{if} &
x_\a\in B^{n-1}_{\g N_j} \text{ and } x_n\in (\d_j/r_j, 2 \d_j/r_j),\\
\\
\zeta_{\g,j}(x) & \text{if}& x\in (B^{n-1}_{\g N_j} \times
I_j)\setminus C_{1,\g N_j},\\
\\
\zeta_{\g,j}\Bigl(x_\a, -2{\d_j\over r_j}-x_n\Bigr)& \text{if}&
x_\a\in B^{n-1}_{\g N_j}\text{ and } x_n\in (-2\d_j/r_j, - \d_j/r_j)\\
\end{array}  \right.
\end{equation}
and then, we extend it by $(2\d_j/r_j)$-periodicity in the $x_n$
direction. The resulting sequence, still denoted by
$\bar\zeta_{\g,j}$, is defined in $\big(B^{n-1}_{\g N_j} \times
\Rb\big)\setminus C_{1,\g N_j}$. Hence, we define on
$\Rb^n\setminus C_{1,\infty}$,
\begin{equation}\label{barzetaj}
\bar\zeta_{j}(x):=\left\{
\begin{array}{lll}
z & \text{in} &\big( \Rb^{n-1}\setminus B^{n-1}_{\g
N_j}\big)\times (0,+\infty),\\
\\
\bar\zeta_{\g,j}(x)& \text{in} & \big(B^{n-1}_{\g N_j}\times
\Rb\big)\setminus C_{1,\g N_j},\\
\\
0 & \text{in} &\big( \Rb^{n-1}\setminus B^{n-1}_{\g
N_j}\big)\times (-\infty, 0)\,. \\
\end{array}  \right.
\end{equation}
Let us now introduce the cut-off functions $\phi_j \in
\C_c^{\infty}((-2\d_j/r_j, 2 \d_j/r_j);[0,1])$ such that
$\phi_j(x_n)=1$ if $|x_n|\le \d_j/r_j$, $\phi_j(x_n)=0$ if
$|x_n|\ge 2\d_j/r_j$ and $|D_n \phi_j |\le c(r_j/\d_j)$. Then, we
introduce our last sequence,
$$
\tilde\zeta_j(x_\a,x_n):= \left\{ \begin{array}{lll} \phi_j(x_n)
\bar\zeta_j (x_\a,x_n) + (1 - \phi_j(x_n)) z & \text{ if } &
(x_\a,x_n) \in \Rb^{n-1} \times (0,+\infty),\\
&&\\
\phi_j(x_n) \bar\zeta_j(x_\a,x_n) & \text{ if } & (x_\a,x_n) \in
\Rb^{n-1} \times (-\infty,0).
\end{array}\right.
$$
Note that
\begin{equation}\label{coincide}
\tilde\zeta_j= \zeta_{\g,j}\quad \text{ in } \quad (B^{n-1}_{\g
N_j}\times I_j )\setminus C_{1, \g N_j}\,.
\end{equation}
Moreover, by (\ref{boundgrad0})-(\ref{coincide}) we have that
\begin{equation}\label{bgradalphatilde}
\sup_{j \in \Nb}\int_{\Rb^n\setminus C_{1,\infty}}
|D_\a\tilde\zeta_j|^p\dx\le c\,,
\end{equation}
while, for every $(a,b)\subset \Rb$, with $a<b$, we have
\begin{equation}\label{bgradntilde}
\int_{\big(\Rb^{n-1}\times (a,b) \big) \setminus C_{1,\infty}}
|D_n\tilde\zeta_j|^p\dx\le c\,,
\end{equation}
for $j$ large enough and $c$ independent of $(a,b)$. Reasoning as
in Proposition \ref{bla1}, with $(0,+\infty)$ and $(-\infty,0)$ in
place of $(0,1)$ and $(-1,0)$, respectively, we can conclude that
there exist $\zeta_1\in L^p(0,+\infty;L^{p^*}(\Rb^{n-1};\Rb^m))$
and $\zeta_2\in L^p(-\infty,0;L^{p^*}(\Rb^{n-1};\Rb^m)) $ such
that, up to subsequences,
$$
\tilde\zeta_j-z\rightharpoonup \zeta_1 \quad\text{in}\quad
L^p(0,+\infty;L^{p^*}(\Rb^{n-1};\Rb^m))
$$
and
$$
\tilde\zeta_j \rightharpoonup \zeta_2 \quad\text{in}\quad
L^p(-\infty,0;L^{p^*}(\Rb^{n-1};\Rb^m))\,.
$$
Moreover, by (\ref{bgradalphatilde}) and (\ref{bgradntilde}), we
have that, up to subsequences, $\tilde\zeta_j$ converges weakly to
$\zeta$ in $W^{1,p}_{\rm loc}(\Rb^n\setminus C_{1,\infty};\Rb^m)$
where
$$
\zeta=\left\{
\begin{array}{lll}
\zeta_1+z & \text{in} & \Rb^{n-1}\times (0,+\infty)\\
\zeta_2 & \text{in} & (\Rb^{n-1}\times (-\infty, 0)) \cup
(B_1^{n-1}
\times\{0\})\,.\\
\end{array}  \right.
$$
In particular, for any compact set $K\subset \Rb^n\setminus
C_{1,\infty}$, we have that
$$
\int_K |D\zeta|^p\dx\le \liminf_{j\to +\infty} \int_K
|D\tilde\zeta_j|^p\dx \le c
$$
for some constant $c$ independent of $K$; hence, we get that
$D\zeta\in L^p(\Rb^n\setminus C_{1,\infty}; \Rb^{m\times n})$
which concludes the description of the limit function $\zeta$.
\hfill$\Box$
\bigskip
\begin{proposition}[Representation formula]
We have
\begin{eqnarray*}
\varphi^{(0)}(z) & = & \inf \Bigg\{\int_{\Rb^n\setminus
C_{1,\infty} } g (D \zeta) \, dx :\;\zeta \in W^{1,p}_{\rm
loc}(\Rb^n \setminus C_{1,\infty};\Rb^m), \, D \zeta \in L^p(\Rb^n
\setminus C_{1,\infty};\Rb^{m\times n}),\\&& \hspace{1.5cm}
\,\zeta-z \in L^{p}(0,+\infty;L^{p^*}(\Rb^{n-1};\Rb^m)) \text{ and
} \zeta \in L^{p}(-\infty,0;L^{p^*}(\Rb^{n-1};\Rb^m))\Bigg\}
\end{eqnarray*}
for every $z\in \Rb^m$.
\end{proposition}

\noindent {\it Proof. } We define
\begin{eqnarray*}
\psi^{(0)}(z) & := & \inf \Bigg\{\int_{\Rb^n \setminus
C_{1,\infty}} g (D \zeta) \, dx :\;\zeta \in W^{1,p}_{\loc}(\Rb^n
\setminus C_{1,\infty};\Rb^m), \, D \zeta \in L^p(\Rb^n \setminus
C_{1,\infty};\Rb^{m\times n}),\\
&& \hspace{1.5cm} \,\zeta-z \in
L^{p}(0,+\infty;L^{p^*}(\Rb^{n-1};\Rb^m)) \text{ and } \zeta \in
L^{p}(-\infty,0;L^{p^*}(\Rb^{n-1};\Rb^m))\Bigg\}
\end{eqnarray*}
and let us prove that $\varphi^{(0)}(z)=\psi^{(0)}(z)$ for every
$z\in \Rb^m$.

By definition of $\varphi_{\g,j}^{(0)}$ (see (\ref{phigj12})), for
every fixed $\eta>0$, there exists $\zeta_{\g,j} \in Y_j^\g(z)$
such that
\begin{equation}\label{liminf-}
\int_{(B_{\g N_j}^{n-1} \times I_j) \setminus C_{1,\g N_j}} g_j(D
\zeta_{\g,j})\, dx \leq \varphi_{\g,j}^{(0)}(z) + \eta\,;
\end{equation}
hence, by Proposition \ref{midi12} (i), (\ref{boundgj}) is
satisfied. Then by Propositions \ref{lzero} and \ref{bla3} we get
that
\begin{eqnarray}\label{liminf-1}
\lim_{j \to +\infty} \varphi_{\g,j}^{(0)}(z) + \eta & \geq &
\liminf_{j \to +\infty} \int_{(B_{\g N_j}^{n-1} \times I_j )
\setminus C_{1,\g N_j}}
g_j(D \tilde \zeta_j)\, dx\nonumber\\
 & \geq &\nonumber \liminf_{j
\to +\infty} \int_{\big(B_N^{n-1} \times (-N,N)\big)\setminus
C_{1,N}} g_j(D \tilde
\zeta_j)\, dx\\
&\ge&\int_{\big(B_N^{n-1} \times (-N,N)\big)\setminus C_{1,N}} g(D
\zeta)\, dx
\end{eqnarray}
for some fixed $N>1$, where $\zeta$ satisfies (\ref{defzeta}). Thus,
passing to the limit in (\ref{liminf-1}) as $N \to +\infty$, $\g \to
0^+$ and $\eta \to 0^+$, it follows that
$$
\varphi^{(0)}(z) \geq \int_{\Rb^n\setminus C_{1,\infty}}
g(D\zeta)\, dx \geq \psi^{(0)}(z)\,.
$$
Let us prove the converse inequality. For any fixed $\eta>0$, let
$\zeta\in W^{1,p}_{\rm loc}(\Rb^n \setminus C_{1,\infty};\Rb^m)$
be as in (\ref{defzeta}) and satisfying
\begin{equation}\label{inbis}
\int_{\Rb^n\setminus C_{1,\infty} }g (D\zeta )\, dx \leq
\psi^{(0)}(z) + \eta.
\end{equation}
For every $j\in \Nb$ and $\g>0$, we consider a cut-off function
$\theta_{\g,j} \in \C^\infty_c(B_{\g N_j}^{n-1};[0,1])$ such that
$\theta_{\g, j}=1$ in $B_{(\g N_j)/2}^{n-1}$, $|D_\a\theta_{\g,j}|
\leq c/\g N_j$ and we define $\zeta_{\g,j} \in Y_j^\g(z)$ by
$$
\zeta_{\g ,j}:= \left\{ \begin{array}{lll} \theta _{\g ,j}
(x_\alpha  )\zeta + (1 - \theta _{\g ,j} (x_\alpha  ))z &
\text{in} & (B_{\g N_j}^{n - 1})^{+(\d_j/r_j)}\\
&&\\
\theta _{\g ,j} (x_\alpha )\zeta & \text{in} & (B_{\g
N_j}^{n-1})^{-(\d_j/r_j)} \cup (B_1^{n-1} \times
\{0\})\,.  \\
\end{array} \right.
$$
Consequently, $\zeta_{\g,j}$ is an admissible test function for
(\ref{phigj12}) and we get that
$$\varphi_{\g,j}^{(0)}(z)  \leq \int_{(B^{n-1}_{\g N_j} \times I_j )
\setminus C_{1,\g N_j}}g_j(D\zeta_{\g,j} )\, dx.
$$
The same kind of computations as those already employed in the
proof of Lemma \ref{reprell} now with $g_j$ in place of $g$ and
with other obvious replacements (see (\ref{est1})-(\ref{est3}))
gives
$$\lim_{j \to +\infty} \varphi_{\g,j}^{(0)}(z) \leq \limsup_{j \to
+\infty} \int_{(B^{n-1}_{\g N_j} \times I_j ) \setminus C_{1,\g
N_j}} g_j (D\zeta)\,dx + o(1)\,,\qquad \text{as}\quad \g\to 0^+\,.
$$
On the other hand, Fatou's Lemma and (\ref{pgrowthWj}) imply
$$\limsup_{j \to +\infty} \int_{(B^{n-1}_{\g N_j} \times I_j )
\setminus C_{1,\g N_j}} g_j (D\zeta)\,dx \leq \int_{\Rb^n
\setminus C_{1,\infty}} g (D\zeta)\,dx + o(1)\,,\qquad
\text{as}\quad \g\to 0^+.$$ Hence by (\ref{inbis}), passing to the
limit as $\g\to 0^+$, we get that
$$\varphi^{(0)}(z) \leq \psi^{(0)}(z) +\eta$$
and by the arbitrariness of $\eta$, the thesis.\hfill$\Box$\\

\begin{rmk}\label{Nadia}{\rm As already recalled, in \cite{Ans} it is
proved that if $\d_j=1$ or $\d_j=\e_j$ then the critical size
$r_j$ of the contact zones is of order $\e_j^{(n-1)/(n-p)}$ or
$\e_j^{n/(n-p)}$, respectively; moreover, the interfacial energy
density is described by the following formula
\begin{eqnarray*}
\varphi(z)&=& \inf\Bigl\{\int_{\rr^n\setminus C_{1,\infty}} g (D
\zeta) \, dx :\;\zeta \in W^{1,p}_{\rm loc}(\Rb^n \setminus
C_{1,\infty};\Rb^m)\,
 \\
&& \qquad\qquad\qquad\qquad\quad\qquad \zeta-z\in
W^{1,p}(\rr_+^n;\rr^m),\ \zeta\in W^{1,p}(\rr_-^n;\rr^m)\Bigr\}
\end{eqnarray*}
where $\rr_{+}^n= \rr^{n-1}\times (0,+\infty)$, $\rr_{-}^n=
\rr^{n-1}\times (-\infty,0)$ (see \cite{Ans} Section 7, the case
$p=q$, with $\rho_{\e_j}=r_j$, $W_p=U_p=W$, $\widehat{W}_p=
\widehat{U}_p=g$ and $\rr_{+,-}^n\cup B^{n-1}_1(0)=\rr^n\setminus
C_{1,\infty}$).

We want to point out that from the analysis we carried on in the
case $\ell=0$ and in particular from
$$
0<R^{(0)}=\lim_{j\to+\infty}\frac{r_j^{n-p}}{\d_j\e_j^{n-1}}
$$
we recovered both the critical sizes founded in \cite{Ans} and
correspondent to the two cases $\d_j=1$ and $\d_j=\e_j$.

Moreover we want to show that $\varphi = \varphi^{(0)}$. We have
to check only the inequality $\varphi\leq \varphi^{(0)}$, the
other one being obvious.

For any fixed $\eta>0$ let $\zeta\in W^{1,p}_{\rm loc}(\Rb^n
\setminus C_{1,\infty};\Rb^m)$ be such that $\zeta-z \in
L^{p}(0,+\infty;L^{p^*}(\Rb^{n-1};\Rb^m))$, $\zeta \in
L^{p}(-\infty,0;L^{p^*}(\Rb^{n-1};\Rb^m))$, $D \zeta \in L^p(\Rb^n
\setminus C_{1,\infty};\Rb^{m\times n})$ and
\begin{equation}\label{minimizing}
\int_{\rr^n\setminus C_{1,\infty}} g (D \zeta) \dx\le
\varphi^{(0)}(z)+\eta\,.
\end{equation}
For every $N>2$ we denote by $B_N$ the $n$-dimensional ball of
radius $N$ centered in zero and by $B_N^\pm$ the set of the points
$x\in B_N$ such that $\pm x_n>0$; we consider a cut-off function
$\theta_N \in \C^\infty_c(B_N;[0,1])$ such that $\theta_N=1$ in
$B_{N/2}$, $|D\theta_N| \leq c/N$ and we define
$$
\bar\zeta:= \left\{ \begin{array}{lll} \theta _N (\zeta - z) + z &
\text{in} & B_N^{+},\\
&&\\
 \theta _N \zeta  & \text{in} & B_N^{-}\cup
(B^{n-1}_1 \times \{0\})
\end{array} \right.
$$
so that $\bar\zeta\in W^{1,p}(B_N\setminus C_{1,N};\rr^m)$,
$\bar\zeta=z$ on $\partial B_N^{+}$ and $\bar\zeta=0$ on $\partial
B_N^{-}$. Hence,
$$
\int_{B_N\setminus C_{1,N}} g(D\bar\zeta)\dx=
\int_{B_{N/2}\setminus C_{1,N/2}} g(D\zeta)\dx +
\int_{(B_N\setminus B_{N/2})\setminus  C_{1,N}}
g(D\bar\zeta)\dx\,;
$$
in particular, by (\ref{pgrowthg}), we have
\begin{eqnarray*}
\int_{(B_N\setminus B_{N/2})\setminus
C_{1,N}}g(D\bar\zeta)\dx&\le& \beta \Bigl(\int_{B^+_N\setminus
B^+_{N/2}} |D\theta_N|^p |\zeta-z|^p\dx +
\int_{B^-_N\setminus B^-_{N/2}} |D\theta_N|^p |\zeta|^p\dx\\
&&+ \int_{(B_N\setminus B_{N/2})\setminus  C_{1,N}} |D\zeta|^p\dx
\Bigr)\\
&\le&{c\over N^p} \Bigl(\int_{B^+_N\setminus
B^+_{N/2}}|\zeta-z|^p\dx +
\int_{B^-_N\setminus B^-_{N/2}} |\zeta|^p\dx\Bigl)\\
&&+ \int_{(\rr^n\setminus B_{N/2})\setminus  C_{1,\infty}}
|D\zeta|^p\dx\,.
\end{eqnarray*}
Since $\zeta-z \in L^{p}(0,+\infty;L^{p^*}(\Rb^{n-1};\Rb^m))$,
$\zeta \in L^{p}(-\infty,0;L^{p^*}(\Rb^{n-1};\Rb^m))$ and $D \zeta
\in L^p(\Rb^n \setminus C_{1,\infty};\Rb^{m\times n})$, we can
easily conclude that
\begin{equation}\label{resto}
\lim_{N\to +\infty}\int_{(B_N\setminus B_{N/2})\setminus
C_{1,N}}g(D\bar\zeta)\dx=0\,.
\end{equation}
Hence, by (\ref{resto}), we deduce
\begin{eqnarray*}
\varphi^{(0)}(z)+\eta&\geq& \int_{\rr^n\setminus C_{1,\infty}} g
(D\zeta) \dx\geq \int_{B_{N/2}\setminus C_{1,N/2}} g(D\zeta)\dx\\
&=& \int_{B_N\setminus C_{1,N}} g(D\bar\zeta)\dx + o(1)\\
&\geq& \inf\Bigl\{\int_{B_N\setminus C_{1,N}} g (D \zeta) \, dx
:\;\zeta\in W^{1,p}(B_N\setminus C_{1,N};\rr^m)\,
 \\
&& \qquad\qquad\qquad\qquad\qquad \zeta=z \,\, {\rm on}\,\,
\partial B_N^{+},\ \zeta=0\,\,{\rm on}\,\, \partial
B_N^{-}\Bigr\}+ o(1)
\end{eqnarray*}
as $N\to +\infty$. Finally, passing to the limit as $N\to
+\infty$, by the arbitrariness of $\eta$, we get
$\varphi^{(0)}\geq \varphi $.

Note that the proof of the explicit formula for $\varphi$ in
\cite{Ans} relies on the fact that $\d_j$ is of order $\e_j$ or
bigger than it, while in Proposition \ref{bla3} and Proposition
7.12 we have to take into account that $\d_j\ll\e_j$. This is the
reason why our proof is different from the one of \cite{Ans} even
if, at the end, the two representation formulas turn out to
coincide. }\end{rmk} \vspace{0.2cm}

\section{Conclusions}

\noindent Our paper deals with the characterization of the
effective energy of weakly connected thin structures through a
periodically distributed contact zone. We highlight the presence
of three different regimes (depending on the mutual rate of
convergence of the radii of the connecting zones and the thickness
of the domain) and for each of them we derive the limit energy by
a $\G$-convergence procedure. For each regime an interfacial
energy term, depending on the jump of the deformation at the
interface, appears in the limit representing the asymptotic memory
of the sieve. We completely describe the interfacial energy
densities by nonlinear capacitary type formulas.

\vspace{1cm}

\textsc{Acknowledgments.} The authors wish to thank Andrea Braides
for suggesting the problem, for fruitful comments and stimulating
discussions.

The research of J.-F. Babadjian has been supported by the Marie
Curie Research Training Network MRTN-CT-2004-505226 `Multi-scale
modelling and characterisation for phase transformations in advanced
materials' (MULTIMAT).

\end{document}